\documentclass{amsart}
\usepackage[margin=2.7cm]{geometry}
\usepackage{systeme}
\usepackage{graphicx}
\usepackage[]{units}
\usepackage{color}
\usepackage{mathrsfs}
\usepackage{bm}
\usepackage[ansinew]{inputenc}
\usepackage{pgfplots}
 \usepackage[breaklinks]{hyperref}
\usepackage{float}
\usepackage{amsmath}
\usepackage{isomath}
\usepackage{amsthm}
\usepackage{amssymb}
\usepackage[foot]{amsaddr}
\usepackage{hyperref}
\usepackage{algorithm,algorithmic}
\usepackage{framed} % Framing content
\usepackage{multicol} % Multiple columns environment
\usepackage{lineno,hyperref}
\usepackage{amsthm}
\usepackage{algorithm}
\usepackage{algorithmic}
\usepackage{tabularx}
\usepackage{tabulary}
\usepackage{subfig}
\usepackage{framed} % Framing content
\usepackage{multicol} % Multiple columns environment
\usepackage{color}
\usepackage{resizegather}
\usepackage{siunitx}
\usepackage{tikz}

\makeatletter
\edef\orig@output{\the\output}
\output{\setbox\@cclv\vbox{\unvbox\@cclv\vspace{0pt plus 20pt}}\orig@output}
\makeatother

\providecommand{\norm}[1]{\lVert#1\rVert}
\DeclareMathOperator{\Tr}{tr}
\newcommand{\ltwonorm}[1]{\langle #1 \rangle_{L_2(\Omega)}}
\newcommand{\frob}[2]{\langle #1, #2 \rangle_F}
\newcommand{\dif}{\mbox{d}}

\usepackage{nomencl} % Nomenclature package
\usepackage{xpatch}

\makenomenclature
\makenomenclature 
\RequirePackage{ifthen}
\let\oldref\ref
\renewcommand{\ref}[1]{(\oldref{#1})}

\newcolumntype{Y}{>{\centering\arraybackslash}X}

\usepackage{amssymb}
\providecommand{\norm}[1]{\lVert#1\rVert}

\DeclareMathAlphabet\mathbfcal{OMS}{cmsy}{b}{n}

\xpatchcmd{\thenomenclature}{%
  \section*{\nomname}% Look for `\section*... etc.
}{% Replace it by 'nothing'
}{\typeout{Success}}{\typeout{Failure}}

\usepackage{ifthen}
\renewcommand{\nomgroup}[1]{%
  \ifthenelse{\equal{#1}{A}}{\item[\textbf{Abreviations}]}{%
    \ifthenelse{\equal{#1}{G}}{\item[\textbf{Symbols}]}{%
      \ifthenelse{\equal{#1}{C}}{\item[\textbf{Abbreviations}]}{%
        \ifthenelse{\equal{#1}{S}}{\item[\textbf{Subscripts}]}{%
          \ifthenelse{\equal{#1}{Z}}{\item[\textbf{Mathematical Symbols}]}{}
        }% matches mathematical symbols
      }% matches Subscripts
    }% matches Abbreviations
  }% matches Greek Symbols
}% matches Roman Symbols

\newcommand{\RA}{\color{black}}

\begin{document}
\nomenclature{$\text{ROM}$}{Reduced Order Model}
\nomenclature{$\text{POD}$}{Proper Orthogonal Decomposition}
\nomenclature{$\text{ALE}$}{Arbitrary Lagrangian Eulerian}
\nomenclature{$\text{D-EIM}$}{Discrete Empirical Interpolation Methods}
\nomenclature{$\text{RBF}$}{Radial Basis Function}

%%% Domain and notation
\nomenclature[G]{$\Omega (\bm \mu)$}{Parametrized Domain}
\nomenclature[G]{$\Gamma_D (\bm \mu)$}{Parametrized Dirichlet Boundary}
\nomenclature[G]{$\Gamma_N (\bm \mu)$}{Parametrized Neumann Boundary}
\nomenclature[G]{$\mathcal{N}(\bm{x},u(\bm{\mu}))$}{Generic parametrized differential operator}
\nomenclature[G]{$\Box_i$}{Variable at the center of the cell $i$}
\nomenclature[G]{$\Box_{if}$}{Variable at the center of the face $f$ of cell $i$}
\nomenclature[G]{$\Box_{\bm n}$}{Normal derivative}
\nomenclature[G]{$\langle \cdot, \cdot \rangle_F$}{Frobenious inner product}
\nomenclature[G]{$\mathcal T$}{FV tessellation}

%%% Unknown matrices and vectors
\nomenclature[G]{$\theta$}{Temperature field}
\nomenclature[G]{$\bm \theta_h$}{FOM solution for $\theta$}
\nomenclature[G]{$\bm \theta_{rb}$}{ROM solution for $\theta$}
\nomenclature[G]{$\bm{M'}$}{Modified mass matrix}
\nomenclature[G]{$\norm{\cdot}_\Omega$}{$L_2$ norm over the computational domain $\Omega$}
\nomenclature[G]{$\bm{A}$}{Discrete FOM Laplace operator}
\nomenclature[G]{$\bm{f}$}{Discrete FOM source term}
\nomenclature[G]{$\bm{A^r}$}{Discrete ROM Laplace operator}
\nomenclature[G]{$\bm{f^r}$}{Discrete ROM source term}
\nomenclature[G]{$\bm a^\theta$}{Reduced coefficients for $\theta$}

%%% Dimensions
\nomenclature[G]{$N_h$}{Number of DOFs for the FOM}
\nomenclature[G]{$N_A$}{Number of D-EIM modes for $\bm A$}
\nomenclature[G]{$N_f$}{Number of D-EIM modes for $\bm f$}
\nomenclature[G]{$N_r$}{Number of DOFs for the ROM}
\nomenclature[G]{$N_{train}$}{Number of training samples}
\nomenclature[G]{$N_{test}$}{Number of testing samples}
\nomenclature[G]{$N_{b}$}{Number of the RBF boundary points}

%%% Spaces and parametrization
\nomenclature[G]{$\mathcal{K}_train$}{Discrete training set}
\nomenclature[G]{$\mathcal{K}_test$}{Discrete testing set}
\nomenclature[G]{$\mathcal P$}{Parameter space}
\nomenclature[G]{$\bm \mu$}{Parameter vector}
\nomenclature[G]{$\bm \kappa_{i_{train}}$}{Sample point in the training set}
\nomenclature[G]{$\bm \kappa_{i_{test}}$}{Sample point in the testing set}

%%% Snapshots Matrices
\nomenclature[G]{$\bm{\mathcal{S}_\theta}$}{Snapshots matrix for Temperature}
\nomenclature[G]{$\bm{\mathcal{S}_A}$}{Snapshots matrix for $\bm{A}$ matrix}
\nomenclature[G]{$\bm{\mathcal{S}_f}$}{Snapshots matrix for $\bm f$ vector}
\nomenclature[G]{$\bm{\mathcal{S}_M}$}{Snapshots matrix for $\bm M$ matrix}

%%% BASES and D-EIM
\nomenclature[G]{$\bm{\varphi^\theta_i}$}{Basis function for $\theta$}
\nomenclature[G]{$\bm{\chi^A_k}$}{Matrix D-EIM basis for $\bm{A}$}
\nomenclature[G]{$\bm{\chi^A_f}$}{Vector D-EIM basis for $\bm{f}$}
\nomenclature[G]{$\bm{b^A}$}{D-EIM coefficients for $\bm{A}$}
\nomenclature[G]{$\bm{c^f}$}{D-EIM coefficients for $\bm{f}$}

\title[]{Efficient Geometrical Parametrization for Finite-Volume based Reduced Order Methods}

\author{Giovanni Stabile$^1$}
\author{Matteo Zancanaro$^1$}
\author{Gianluigi Rozza$^1$}
\address{$^1$ mathLab, Mathematics Area, SISSA, via Bonomea 265, I-34136 Trieste, Italy}

\maketitle

\begin{abstract}
In this work, we present an approach for the efficient treatment of parametrized geometries in the context of POD-Galerkin reduced order methods based on Finite Volume full order approximations. On the contrary to what is normally done in the framework of finite element reduced order methods, different geometries are not mapped to a common reference domain: the method relies on basis functions defined on an average deformed configuration and makes use of the Discrete Empirical Interpolation Method (D-EIM) to handle together non-affinity of the parametrization and non-linearities. In the first numerical example, different mesh motion strategies, based on a Laplacian smoothing technique and on a Radial Basis Function approach, are analyzed and compared on a heat transfer problem. Particular attention is devoted to the role of the non-orthogonal correction. In the second numerical example the methodology is tested on a geometrically parametrized incompressible Navier--Stokes problem. In this case, the reduced order model is constructed following the same segregated approach used at the full order level.
\end{abstract}

\section{Introduction and Motivation}\label{sec:intro}
In several cases there is the need to repeatedly solve Partial Differential Equations (PDEs) on different parametric domains. Such situations occur for example in the case of geometric design optimization where different parametric domains are in need of being tested in order to determine the configuration which is maximizing or minimizing a certain quantity of interest. 

Finite Volume (FV) approximations, due to the fact that they enforce conservation at local level, are particularly suited to be employed for the discretization of systems of Partial Differential Equations (PDEs) generated by conservation laws. Conservation laws model a large range of different physical problems such as fluid dynamics, heat transfer, solid mechanics, etc. In particular,  the FV method is particularly used and convenient to model hyperbolic conservation laws that are common in fluid dynamics problems.

For this reason, the finite volume method is widespread in many different engineering fields such as mechanical, aerospace, civil, nuclear, \dots and non-engineering ones such as meteorology, environmental marine sciences, medicine, etc. As also for other full order discretization techniques (Finite Elements, Finite Differences, Spectral Elements) in several cases (e.g. uncertainty quantification, inverse problems, optimization, real-time control, etc) the numerical simulation of the governing equations using standard techniques, becomes not affordable. This fact is particularly evident for the case of geometrical optimization \cite{Lombardi_2012} on which this work is focused. A viable way to reduce the computational burden is given by the reduced order modeling methodology; for a comprehensive review of the reduced order modeling methodology the interested reader may see \cite{QuaMaNe15,HeRoSta16,BeOhPaRoUr17}. 

Reduced Order Models (ROMs) and in particular projection-based reduced order models, which are the focus of this work, have been applied to a large variety of different mathematical problems based on linear elliptic equations \cite{Rozza2008229}, linear parabolic equations \cite{Grepl2005} and also non linear problems \cite{Veroy2003,Grepl2007}. 

Some of the ingredients and the issues encountered to generate a projection-based reduced order model are common to all the full order discretization technique but some aspects require particular attention if we change full order discretization. The usage of the reduced basis method for shape optimization in a finite element setting has been widely studied and the state of the art counts already several scientific contributions \cite{Rozza:148535,Rozza_2013,Rozza_2011,Iapichino_2016}. 
%On the other hand, geometrical parametrization using the reduced basis method and the finite volume method has not been fully studied and needs therefore further development. 

Summarizing pioneering and more recent works dealing with the coupling of the finite volume methods with projection based reduced order methods it is worth mentioning \cite{Haasdonk2008,HOR08} that present one of the first contributions dealing with FV and the reduced basis (RB) method and \cite{Drohmann2012} dealing with the empirical interpolation method in a FV setting. For more recent contributions it is worth mentioning \cite{Stabile2017, stabile_stabilized} where pressure stabilization techniques normally employed in a finite element setting have been adapted to a FV one and \cite{Carlberg2018} that proposes a new structure preserving strategy specifically tailored for a FV setting, \cite{saddam2018,Lorenzi2016151,GeStaRoBlu2018} that extend FV-POD-Galerkin ROMs to a turbulence setting and to thermal mixing problems. Focusing on works that are specifically dealing with geometrically parametrized problems in a finite volume setting we mention \cite{LeGresley2001, Zimmermann2010} that focus the attention on inviscid Euler equations, \cite{Washabaugh2016} focusing on turbulent compressible Navier-Stokes equations and \cite{Zahr2014} that deals with PDE-constrained optimization problems.

The main focus of this work is the development of computational strategies for the geometrical parametrization of reduced order methods starting from a finite volume full order discretization. As highlighted in \autoref{subsec:geom_par}, this issue, respect to a finite element setting, requires particular attention. The strategy usually employed to deal with geometrical parametrization with reduced order models is to write all the equations in a reference domain making use of a map $\mathcal{M}(\bm \mu):\Omega(\bm \mu)\to \tilde\Omega$ which transforms the equations from the physical domain $\Omega(\bm \mu)$ to the reference one $\tilde \Omega$. This approach has been adopted successfully in several cases starting from finite element discretization. The affinity of the differential operators can be reconstructed by domain decomposition and piecewise affine reference mappings \cite{Rozza2008229,ballarin2015supremizer,Milani2008} or by the usage of empirical interpolation techniques \cite{Jggli2014,Lassila2010,Rozza2009}. The first approach cannot be easily transferred to a finite volume setting (see \autoref{subsec:geom_par}) while the second approach has been used in some of the previously mentioned works dealing with geometrical parametrization in a finite volume environment. In the present work, all the possible domain configurations are solved in the physical domain and the reduced basis functions are computed by means of an inner product computed on an average deformed configuration similarly to what done in \cite{LeGresley2001}. 

The online-offline splitting is guaranteed by the usage of a discrete empirical interpolation procedure at both matrix and vector level. The resolution of the equations is carried out on the physical domain where the mesh is deformed by the usage of an Arbitrary Lagrangian Eulerian Framework \cite{Hirt_1997}.

The manuscript is organized as follows: \autoref{sec:FOM} briefly introduces the parametrized mathematical problem, \autoref{sec:FVM} presents the full order finite volume formulation and its main features with a particular attention on the different mesh motion strategies. In \autoref{sec:rom} we describe the methodologies used to construct the reduced order model with emphasis on the empirical interpolation techniques used to ensure an efficient offline-online splitting. In \autoref{sec:num_ex} we test the proposed methodologies on a geometrically parametrized heat transfer problem and the attention here is focused on the different mesh motion strategy. In \autoref{sec:num_ex2} the methodology is tested instead on a steady incompressible Navier-Stokes setting. Finally in \autoref{sec:outlooks} some conclusions and perspectives for future works are given. 

\section{The parametrized mathematical problem}\label{sec:FOM}
The interest of this work is on generic stationary linear or nonlinear PDEs describing conservation laws whose domain is parametrized by a parameter vector $\bm{\mu} \in \mathcal{P} \subset \mathbb{R}^k$, where $\mathcal{P}$ is a $k$-dimensional parameter space. Given a generic scalar or vectorial quantity of interests $u(x,\bm{\mu})$ the abstract problem reads:
\begin{equation}\label{eq:abstract_problem}
\begin{cases}
\mathcal{N}(x,u(\bm{\mu})) = 0 &\mbox{ in } \Omega(\bm{\mu}),\\
u(x,\bm\mu) = g_D(x,\bm{\mu}) &\mbox{ on } {\Gamma_D(\bm{\mu})},\\
u_{\bm{n}}(x,\bm{\mu}) = g_N(x,\bm{\mu}) &\mbox{ on } {\Gamma_N(\bm{\mu})},
\end{cases}
\end{equation}
where $\mathcal{N}(x,\bm{\mu})$ is a generic linear or nonlinear differential operator, $\Omega(\bm{\mu}) \subset \mathbb{R}^d $ is a bounded parametrized domain, $\Gamma_D(\bm{\mu})$ and $\Gamma_N(\bm{\mu})$ are its parametrized Dirichlet and Neumann boundaries respectively with the relative boundary conditions $g_D(x,\bm{\mu})$ and $g_N(x,\bm{\mu})$, while the subscript $\Box_{\bm{n}}$ indicates the normal derivative. Reduced Basis (RB) methods are based on the assumption that the response of such parametrized PDE is approximated by a reduced number of dominant modes and therefore it is possible to project the original equations onto a low dimensional subspace in order to reduced the computational cost. In this work this low-dimensional space is obtained using a Proper Orthogonal Decomposition (POD) approach. More details on this aspect are reported in \autoref{sec:rom}.

\section{The finite volume discretization and the geometrical parametrization}\label{sec:FVM}
In this section we recall the main features of the full order discretization method used to discretize the mathematical problem and the techniques used to deal with the parametrization of the geometry. We stress the attention on key aspects of this full order discretization technique that will play a crucial role in order to ensure an accurate and efficient model reduction strategy for parametrized geometries. 

The abstract problem of \autoref{eq:abstract_problem} is discretized using a finite volume approach. The procedure used here is recalled for a simple heat transfer problem but the procedure is general and can be easily extended also to other type of equations. On the contrary respect to projection-based (FEM, SEM) full order discretization techniques, which start from a weak formulation of the governing equations, the starting point of the finite volume method is the integral form of the equations written in conservative form. For a simple steady heat transfer problem on a parametrized domain, which is described by the Poisson equation, the formulation reads:

\begin{equation}\label{eq:poisson_problem}
\begin{cases}
\int_{\tilde{\Omega}(\bm{\mu})}\mbox{div}(\alpha_\theta \nabla \theta) \mbox{d}\omega = \int_{\tilde{\Omega}(\bm{\mu})} f \mbox{d}\omega &\forall \tilde{\Omega}({\bm{\mu}}) \in {\Omega}({\bm{\mu}}),\\
\theta(x,\bm{\mu}) = \theta_D(x,\bm{\mu}) &\mbox{ on } {\Gamma_D(\bm{\mu})},\\
\theta_{\bm{n}}(x,\bm{\mu}) = \theta_N(\bm{\mu},x) &\mbox{ on } {\Gamma_N(\bm{\mu})},
\end{cases}
\end{equation}
where $\theta$ is the temperature field, $\alpha_\theta$ is the thermal diffusivity, $f$ is a generic forcing term, $\tilde{\Omega}({\bm{\mu}})$ is an arbitrary control volume inside $\Omega(\bm\mu)$, $\theta_D$ and $\theta_N$ are Dirichlet and Neumann boundary condition for the temperature field while the other terms assume the same meaning of those reported in the abstract problem formulation (\autoref{eq:abstract_problem}).
\subsection{The finite volume discretization}\label{subsec:fvm_disc}
In order to solve the problem of \autoref{eq:poisson_problem}, the domain $\Omega(\bm{\mu})$ is subdivided into a tessellation $\mathcal{T}(\bm{\mu}) = \{ \Omega_i(\bm{\mu}) \}_{i=1}^{N_{h}}$ composed by a set of $N_h$ convex and non overlapping polygonals in $2D$ or polyhedrons in $3D$ (finite volumes) such that $\Omega(\bm{\mu}) = \bigcup_{i=1}^{N_{h}} \Omega_i(\bm{\mu})$ and $\Omega_i(\bm{\mu}) \bigcap \Omega_j(\bm{\mu}) = \emptyset \mbox{ for } i\neq j$. Once the tessellation is set, the divergence term is transformed into a surface integral using the Gauss' theorem and the integral terms of \autoref{eq:poisson_problem} can be discretized as follows:
\begin{gather}\label{eq:laplace}
\sum_{i=1}^{N_{h}}\int_{\Omega_i(\bm{\mu})} \mbox{div}(\alpha_\theta \nabla \theta) \mbox{d}\omega = 
\sum_{i=1}^{N_{h}}\int_{\partial\Omega_i(\bm{\mu})} \bm{n} \cdot (\alpha_\theta \nabla \theta) \mbox{d}s = 
\sum_{i=1}^{N_{h}}\sum_{f=1}^{N_{f}} \alpha_{\theta_{if}} \bm{S_{if}}  \cdot (\nabla \theta)_{if}, \\ \nonumber
\int_{\tilde \Omega(\bm{\mu})} f \mbox{d}\omega = 
\sum_{i=1}^{N_{h}}\int_{\Omega_i(\bm{\mu})}f \mbox{d}\omega = 
\sum_{i=1}^{N_{h}} f_i V_i.
\end{gather}

\begin{figure}
\centering
\includegraphics[height=0.3\textwidth]{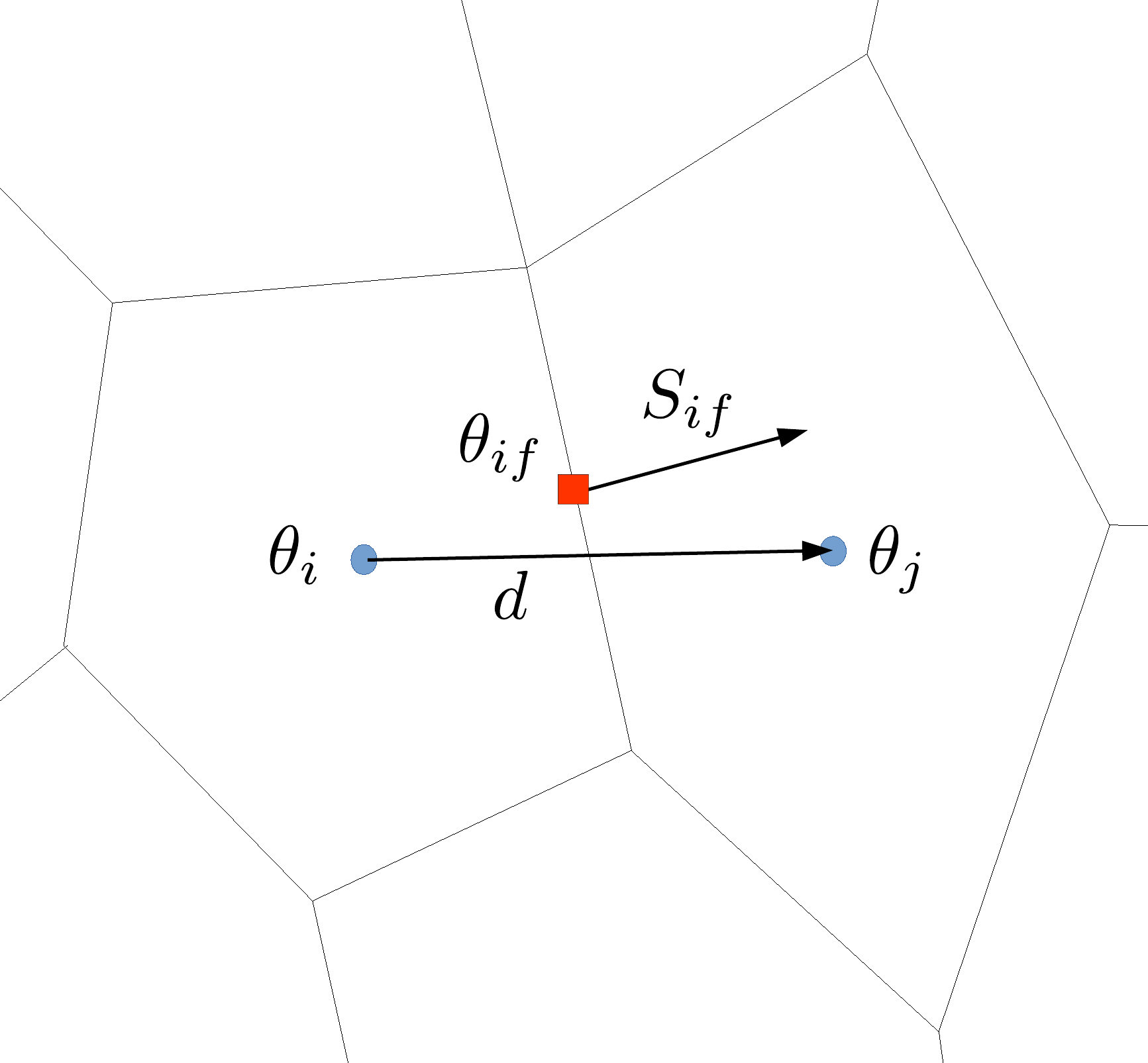}
\caption{Example of a non-orthogonal mesh for a $2$-dimensional problem.}
\label{fig:cell_ortho}
\end{figure}

In the above expressions the subscript $\Box_i$ indicates the value at the centre of the $i$-cell, the subscript $\Box_{if}$ indicates the value of the variables at the centre of the $f$-face relative to the $i$-cell, $\bm n$ represents the unit normal vector respect to the control volume surface, $V$ the cell volumes and $\bm S$ the surface vectors. One of the key points of the finite volume method, especially for advection dominated problem, is the evaluation of the conservative variables at the centre of the faces starting from their values at the centre of the cells. Several numerical schemes are in fact possible such as central differencing, first-order upwind, second-order upwind, Rusanov, MUSCL, etc. For the particular case the term $\bm{S_{if}}\cdot(\nabla \theta)_{if}$, which represents the flux of the gradient of the temperature field passing trough the face $f$, can be evaluated using both an explicit or implicit scheme with or without orthogonal correction. In case of orthogonal meshes, the flux relative to the gradient term can be estimated as:
\begin{equation}\label{eq:gradient_comp}
\bm{S_{if}} \cdot (\nabla \theta )_{if} = \bm{S_{if}} \cdot \frac{\theta_i - \theta_j}{|\bm{d}|},
\end{equation}
where by orthogonal meshes we mean those cases where the surface vector $\bm{S_{if}}$ (see \autoref{fig:cell_ortho}, relative to the face which is dividing two cells, is parallel to the distance vector $\bm{d}$ which is connecting the cell centres. Unfortunately, orthogonal meshes, especially considering complex geometries, are not common and in many cases even not realizable. In the present case, even starting from a fully orthogonal grid, after the mesh motion, we will end up with a non-orthogonal grid. For this reason, in order to have accurate results, \autoref{eq:gradient_comp} needs to take into account also a non-orthogonal correction. In the present work, following what presented in \cite{Jasak1996}, when a non-orthogonal correction is considered, the product $\bm{S_{if}} \cdot (\nabla \theta )_{if}$ is split into two parts:
\begin{equation}\label{eq:nonortho}
\bm{S_{if}} \cdot (\nabla \theta )_{if} = |\bm{\Delta_{if}}| \cdot \frac{\theta_i - \theta_j}{|\bm{d}|} + \bm{k_{if}} \cdot (\nabla \theta )_{if},
\end{equation}
where $\bm{S_{if}} = \bm{\Delta_{if}} + \bm{k_{if}}$ is decomposed into two vectors with $\bm{\Delta_{if}}$ parallel to $\bm{d}$ while $\bm{k_{if}}$ can be selected using different approaches such as minimum correction approach, orthogonal correction approach, over-relaxed approach, etc. We refer to \cite{Jasak1996} for more details. Therefore, the total contribution is split into the first orthogonal contribution and into a non orthogonal one. The term $(\nabla \theta )_{if}$ which represents the value of the gradient at the centre of the face is calculated explicitly via interpolation after the calculation of the value of the gradient at the centres of the neighboring cells $(\nabla \theta )_{i}$ and $(\nabla \theta )_{j}$. Then the value at the centre of the face $(\nabla \theta )_{if}$ is evaluated using an interpolation scheme:
\begin{equation}
(\nabla \theta )_{if} = w_x (\nabla \theta )_{i} + (1-w_x)(\nabla \theta )_{j},
\end{equation}
where $w_x$ is a proper weight that depends on the employed interpolation scheme. Moreover, in some cases, in order to ensure that the discretized diffusion term preserves its bounded behavior also after the correction, the non-orthogonal contribution can be limited. Therefore, when using a correction for the computation of the flux of the gradient one needs to resolve the discrete equations using an iterative scheme because the term $(\nabla \theta)_{if}$ is based on the previous converged value of the conservative variable $\theta$. This additional term act like a source term in the discretized equations.
It is beyond the scope of this work to go into the details concerning the different methodologies used to deal with non-orthogonal meshes but it was important to briefly recall them since they play an important role also at the reduced order level. As we will see in the numerical examples different non-orthogonal correction approaches will lead to different dimensions of the computational stencil that we need to store during the D-EIM procedure and will also lead to different values of the error respect to the full order model. In the numerical examples we will examine two different cases, one without non-orthogonal correction and one with a bounded orthogonal correction. 

Once that the discretization of all the terms inside the continuous formulation is described it is finally possible to recast  \autoref{eq:laplace} in matrix form and to transform it into the following system of linear equations:
\begin{equation}\label{eq:FOM}
\bm{A}(\bm{\mu})\bm{\theta}(\bm{\mu})=\bm{f}(\bm{\mu}).
\end{equation}
where $\bm{A}(\bm{\mu}) \in \mathbb{R}^{N_h \times N_h}$ is a matrix which represents the discretized Laplace operator, $\bm{\theta} \in \mathbb{R}^{N_h}$ is a vector of unknown values, $\bm{f}(\bm{\mu}) \in \mathbb{R}^{N_h}$ is a source term vector which accounts for physical source terms, for implicit non-orthogonal corrections and for boundary conditions. 
\subsection{Geometrical parametrization}\label{subsec:geom_par}
Once the finite volume discretization is defined one important aspect is the selection of the strategy used to solve a geometrically parametrized problem. The usual strategy applied in the reduced order modeling community, and more in general for a finite element discretization, is to construct a map $\mathcal{M}(\bm{\mu}):\Omega(\bm{\mu})\to\hat\Omega$ which maps all the possible domain configurations to a reference one $\hat\Omega$. This is also the procedure employed in a finite volume setting, for a linear mapping, in \cite{Drohmann2009}. Unfortunately this operation is affordable only for simple cases and for certain types of schemes and creates the additional complexities:

\begin{itemize}
\item In a finite volume code, respect to a finite element one, also at the element level, integrals are not written in a reference domain setting. In finite elements, in order to perform the integration of the shape functions, all the elements are mapped to a reference domain. In the finite volume method such operation is not performed and therefore to write integrals in a reference setting, making use of a mapping, would be much more intrusive. 
\item For a general nonlinear map, also in the finite element case, the equations written into a reference domain, transform a pure diffusion equation into a diffusion-advection equation \cite{Drohmann2009}.
\item The finite volume method, with respect to the finite element one, involves also surface integrals. Even though it is possible to map infinitesimally small distances, areas and volumes from a real to a reference domain, and to write the continuous equations in a reference framework, at the discrete level there are additional source terms created by non-orthogonality (we refer to \autoref{subsec:fvm_disc} for this issue) and by non-linear interpolation schemes (i.e. upwinding schemes). The finite volume method involves the computation of the values at centre of the faces starting from the values at the centre of the cells. Unless one is using a central differencing scheme, where the value at the centre of the faces depends only on its distance respect to the centre of the cells to write this operation in a reference domanin becomes not trivial.
\item One of the great advantages of the finite volume consists into the possibility of using a mesh made of polyhedra with an arbitrary number of faces. Writing equations in a reference domain using polyhedra with arbitrary shape may be extremely challenging and would be realistically possible only for tetrahedral elements. This would strongly limit the versatility of the method.
\end{itemize}
For the above reasons, in the present work we opted to work with the equations directly on the physical domain and, for each new value of the input parameters, to deform the mesh and to solve the full order problem on the real physical domain. One of the drawbacks of this approach is given by the fact that an additional problem related to the motion of the grid points must be solved. Moreover, all the different snapshots are defined in a different domain and therefore an adjustment is needed. These aspects will be examined in details in the next subsections.

The proposed approach is not the only one possible to avoid the usage of a reference domain approach. In \cite{1807.07753v2,1807.07790v2,KaBaRO18} the authors make use of immersed methods during the full order stage to deal with the geometrical parametrization and to avoid the mapping to a reference domain.
\subsection{Mesh Motion Strategies}
As mentioned in the previous section, in this work we rely on an Arbitrary Lagrangian-Eulerian (ALE) formulation \cite{Hirt_1997} for the treatment of all the possible different geometric configurations. Therefore, the shape of the computational domain changes according to the parameter values. Working with such an approach, we have to distinguish between the motion of boundary points, which is a given data defined by the parameter values and by the geometrical parametrization, and the motion of the internal points that depends on the employed mesh motion strategy. Dealing with a full order finite volume discretization, we highlight here, as seen in the description of the full order discretization, that the accuracy of the solution is strongly dependent on the quality of the mesh. For this reason, in the selection of the mesh motion strategy it is crucial to select an algorithm which produces a deformed mesh that has a non-orthogonality value as small as possible. In this work we analyze two different strategies to solve the mesh motion problem and we highlight their features in view of the construction of a reduced order model. Here the two methods are briefly recalled and only the aspects relevant for reduced order modeling perspectives are detailed. The interested reader may refer to \cite{jasak2006automatic,deBoer2007} and to references therein for more details. 
\subsubsection{Laplace equation with variable diffusivity}

The first type of approach is based on the solution of a Laplace equation with a variable diffusivity:
\begin{gather}
\mbox{div} (\gamma \nabla \bm{s}) = 0 \mbox{ in } \Omega, \\
\bm{s} = \bm{s}_D \mbox{ on } \Gamma_D,
\end{gather}
where $\bm{s}$ is the displacement field of the grid points and $\gamma$ is a diffusivity field. For the diffusivity field several options are possible such as a constant value, a value proportional to the inverse of the distance respect to moving boundaries, etc. In this work it is chosen to use a diffusivity field which is equal to the inverse of the square of the distance respect to the moving boundaries $r$:
\begin{equation}
\gamma = r^{-2}.
\end{equation}
This value of the diffusivity coefficient proved to produce a smooth mesh motion and therefore a good quality of the deformed mesh in cases with moderate mesh deformations \cite{jasak2006automatic}.
The computation of the displacement of the grid points, at the full order level, requires the resolution of a linear sparse system that has the same dimension of the full order problem multiplied by the number of physical dimensions of the problem $d=2,3$ (i.e. $N_h \cdot d$). The problem is in fact discretized as:
\begin{equation}\label{eq:mesh_lapl}
\bm A_D s = \bm b_D (\mu),
\end{equation}
where $\bm A_D \in \mathbb{R}^{N_h\cdot d \times N_h\cdot d}$ is the discretized Laplace operator and the term $\bm b_D \in \mathbb{R}^{N_h\cdot d}$ which depends on the boundary conditions $\bm s_D$ is the source term. 
Therefore, also during the online stage, when a new parameter value needs to be simulated, one has to solve a full order mesh motion problem. For general deformation mappings it is usually not possible to recover an affine decomposition of the boundary condition $\bm s_D$ and therefore also for the mesh motion problem, in order to obtain an efficient offline-online splitting, we have to rely on hyper-reduction techniques. More details on this aspect are reported in the numerical example sections.

\subsubsection{Radial Basis functions}\label{subsec:RBF}
The second examined approach consists in a Radial Basis Function interpolation approach. In this approach the displacements of all the mesh points not belonging to the moving boundaries are approximated by a Radial Basis Function (RBF) interpolant function:
\begin{equation}
\bm{s}(\bm{x}) = \sum_{i=1}^{N_{b}} \beta_i \xi (\norm{\bm{x}-\bm{x}_{b_i}}) + q(\bm{x}),
\end{equation}
where $\bm{x}_{b_i}$ are the coordinates of points for which we know the boundary displacements, $N_b$ is the number of points on the boundary, $\xi$ is a given basis function, $q(\bm{x})$ is a polynomial. The coefficients $\beta_i$ and the polynomial $q(\bm{x})$ are obtained by the imposition of interpolation conditions
\begin{equation}
\bm{s}(\bm{x}_{b_i}) = \bm{d}_{b_i},
\end{equation}
where $\bm{d}_{b_i}$ is the displacement value at the boundary points and by the additional requirement:
\begin{equation}
\sum_{i=1}^{N_b}\beta_i q(\bm{x}_{b_i}) = 0.
\end{equation}
In the present case, we select basis functions for which it is possible to use linear polynomials $q(\bm{x})$. For more informations concerning the selection of the order of polynomials, respect to the type of basis functions, the reader may see \cite{Beckert2001}. Finally the values of the coefficients $\beta_i$ and the coefficients $\delta_i$ of the linear polynomial $q$ can be obtained solving the linear problem:
\begin{equation}\label{eq:RBF_system}
\begin{bmatrix}
\bm{d}_b \\ 0
\end{bmatrix}
=
\begin{bmatrix}
\bm{M}_{b,b} & P_b \\ P_b^T & 0  
\end{bmatrix}
\begin{bmatrix}
\bm{\beta} \\ \bm{\delta}
\end{bmatrix},
\end{equation}
where $\bm{M}_{b,b} \in \mathbb{R}^{N_b\times N_b}$ is a matrix containing the evaluation of the basis function $\xi_{b_i b_j} = \xi (\norm{\bm{x}_{b_i}-\bm{x}_{b_j}})$ and $\bm{P}_b \in \mathbb{R}^{N_b\times (d+1)}$ is a matrix where $d$ is the spatial dimension. Each row of this matrix, that contains the coordinates of the boundary points, is given by $(\bm P_b)_{i\bullet} = \begin{bmatrix}1 & \bm{x_{b_i}}\end{bmatrix}$. Once the system of \autoref{eq:RBF_system} is solved one can obtain the displacement of all the internal points using the RBF interpolation:
\begin{equation}
\bm d_{in_i} = \bm s(\bm x_{in_i}),
\end{equation}
where $\bm x_{in_i}$ are the coordinates of the internal grid points. The computation of the displacement of the grid points entails the resolution of a dense system of equations that has dimension $N_b + d + 1$. Usually, the number of boundary points $N_b$ is much smaller respect to the number of grid points $N_h$. Moreover, as it will be shown in the numerical examples, usually it is sufficient to select only a subset of the boundary points and therefore the dimension of this system is further reduced. Making a comparison with the Laplace equation approach, the RBF methodology entails the resolution of a smaller but dense system of equations. In the numerical examples it will also be done a comparison in terms of mesh qualities after the mesh deformation. 

\section{The POD-Galerkin reduced order model}\label{sec:rom}

\autoref{eq:FOM} entails the resolution of a possibly large system of equations. The resolution of this system of equations, in a multi-query setting, may become unfeasible. In this work, in order to alleviate the computational burden, we rely on a projection based ROM and in particular on a POD-Galerkin reduced order model. Projection-based reduced order models exploit the fact that, in most of the cases, the solution manifold lies in a low dimensional space and can be therefore approximated by a linear combination of a reduced number of properly selected global basis functions:
\begin{equation}
\theta(\bm{\mu},\bm x) = \sum_{i=1}^{N_r} {a^\theta_i}(\bm{\mu})\varphi^\theta_i(\bm{x}),
\end{equation}
where $\varphi^\theta_i(\bm{x})$ are parameter independent basis functions and $a^\theta_i(\bm{\mu})$ are parameter dependent coefficients. There exist different approaches to construct the set of basis functions such as the greedy approach, the proper orthogonal decomposition, the proper generalized decomposition, etc \cite{BeOhPaRoUr17,QuaMaNe15,HeRoSta16,ChiHuRo16,Chinesta2011,Dumon20111387}. We decided here to rely on a POD approach. Given a parameter vector $\bm{\mu} \in \mathcal{P}$ we select a finite dimensional training set $\mathcal{K}_{train} \subset \mathcal{P}$ and for each of the possible combinations of the parameter values we solve a full order problem. The full order order problem is then solved for each $\bm \kappa_{j_{train}} \in \mathcal{K}_{train} = \{\bm \kappa_{1_{train}}, \dots, \bm \kappa_{N_{train}}\}$ where $\bm \kappa_{j_{train}}$ is the $j$-sample belonging to the finite dimensional training set $\mathcal{K}_{train}$. Each sample $\bm \kappa_{j_{train}}$ corresponds to a certain full order solution $\bm \theta_j$ and the ensemble of whole the full order solutions gives the snapshots matrix:
\begin{equation}
\bm{\mathcal{S}_\theta} = [\bm{\theta}(\bm \kappa_{1_{train}}), \dots,\bm{\theta}(\bm \kappa_{N_{train}})] \in \mathbb{R}^{N_h \times N_{train}}.
\end{equation}
Once the snapshots matrix is set it is possible to apply POD in order to generate a reduced basis space to be used for the projection of the governing equations. Given a general scalar or vectorial parametrized function $\bm{u}(\bm{\mu}):\Omega \to \mathbb{R}^d$, with a certain number of realizations $\bm{u}_1,\dots, \bm{u}_{N_{train}}$, the POD problem consists in finding, for each value of the dimension of POD space $N_{POD} = 1,\dots,N_{train}$, the scalar coefficients $a_1^1,\dots,a_1^{N_{train}},\dots,a_{N_{train}}^1,\dots,a_{N_{train}}^{N_{train}}$ and functions $\bm{\varphi}_1,\dots,\bm{\varphi}_{N_{train}}$ that minimize the quantity:
%%      POD ENERGY - eq:pod_energy
%%
\begin{gather}\label{eq:pod_energy}
E_{N_{POD}} = \sum_{i=1}^{N_{train}}||\bm{u}_i-\sum_{k=1}^{N_{POD}}a_i^k \bm{\varphi_k}||^2_{L_2(\Omega)} \hspace{0.5cm}\forall\mbox{ } N_{POD} = 1,\dots,N\\
\mbox{ with } \ltwonorm{\bm{\varphi}_i,\bm{\varphi}_j} = \delta_{ij} \mbox{\hspace{0.5cm}} \forall \mbox{ } i,j = 1,\dots,N_{train} .
\end{gather}
It is well known \cite{Kunisch2002492} that the minimization of \autoref{eq:pod_energy} reduces to the following eigenvalue problem:
%%%      APPROXIMATION MATRIX FORM - eq:min1
%%%
\begin{gather}\label{eq:eigenvalue_problem}
\bm{{C}^u}\bm{Q}^u = \bm{Q^u}\bm{\lambda^u} ,\\
{C}^u_{ij} = \ltwonorm{\bm{u}_i,\bm{u}_j} \mbox{\hspace{0.5cm} for } i,j = 1,\dots,N_{train},
\end{gather}
%%%
%%%
where $\bm{Q}^u \in \mathbb{R}^{N_{train} \times N_{train}}$ is the matrix containing the eigenvectors while $\bm{\lambda^u} \in \mathbb{R}^{N_{train} \times N_{train}}$ is the matrix containing the eigenvalues.
We remind here that the interest is into parametrized geometries and that we do not map all the possible configurations to a reference one. For this reason, in this case, there is an additional difficulty given by the fact that the different snapshots ``live'' in different domains. The method used here consists into the usage of the method of snapshots where the correlation matrix of \autoref{eq:eigenvalue_problem} is constructed using an inner product with a modified mass matrix $\bm{M}'$ which is referred to an ensemble average of all the different mesh configurations:
\begin{equation}\label{mass matrix}
C^\theta_{ij} = \bm{\theta}_i^T \bm{M}'\bm{\theta}_j.
\end{equation}
In the case of snapshots defined all in the same domain, the above expression reduces to the standard $L_2$ inner product. Operating with this procedure, during the offline stage, we need to additionally compute the mass matrix associated with the ensemble average of the parameter values inside the training set:
\begin{equation}
\bm{M'}(\overline{\bm \kappa})\quad \mbox{with}\quad \overline{\bm \kappa} = \frac{1}{N_{train}}\sum_{i = 1}^{N_{train}} \bm \kappa_{i}.
\end{equation}

Once the modified correlation matrix is assembled and the eigenvalue problem is solved, one can compute the basis functions in a standard way:
\begin{equation}
\bm{\varphi^\theta_i} = \frac{1}{N_{train}\lambda_i^\theta}\sum_{j=1}^{N_{train}} \bm{\theta}_j Q^\theta_{ij}.
\end{equation}
Based on the eigenvalue decay of the POD eigenvalues, we can select only the first $N_r$ basis functions and assemble the matrix:
\begin{equation}
\bm{L} = \left[ \bm{\varphi}^\theta_1, \dots, \bm{\varphi}^\theta_{N_r}\right] \in
\mathbb{R}^{N_h\times N_r}.
\end{equation}
Once the bases are set, it is possible to approximate the solution vector with:
\begin{equation}
\bm \theta \approx \bm{\theta}_{rb} = \sum_{i=1}^{N_r} a^\theta_i(\bm{\mu})\bm{\varphi^\theta}_i(\bm{x}),
\end{equation}
and finally the original problem can be projected onto the POD space giving rise to the following reduced problem:
\begin{equation}\label{eq:reduced_problem}
\bm{L}^T \bm{A}(\bm{\mu}) \bm{L} \bm{a^\theta} = \bm{L}^T \bm{f}(\bm{\mu}),
\end{equation}
that can be rewritten as
\begin{equation}
\bm{A^r}(\bm \mu) \bm{a^\theta} = \bm{f^r}(\bm{\mu}),
\end{equation}
where $\bm {A^r} \in \mathbb{R}^{N_r \times N_r}$ is the discrete reduced parametrized differential operator, $\bm{a^\theta} \in \mathbb{R}^{N_r}$ is the unknown vector of reduced coefficients and $\bm{f^r}(\bm \mu) \in \mathbb{R}^{N_r}$ is the reduced source term. Before solving the reduced problem for a new value of the parameter one has to solve the mesh motion problem, to assemble the discretized operator $\bm{A}(\bm{\mu})$ and the source term $\bm{f}(\bm{\mu})$ and finally to perform the projection. The resolution of the mesh problem has the same dimension of the full order problem, as the assembly of the discretized operator. In order to ensure an efficient offline-online splitting it is crucial to perform full-order operations only during the offline stage. The approach employed here to ensure an efficient offline-online splitting is to rely on the discrete variant of the empirical interpolation method (D-EIM) \cite{BARRAULT2004667,Chaturantabut2010,Drohmann2012}. Other approaches are also possible and among them we list the gappy-POD \cite{Willcox2006} or the GNAT method \cite{Carlberg2013,Carlberg2010}. Since the mass matrix of \autoref{mass matrix} is computed on an ensemble average of all the deformed meshes, this deformed configuration can be considered as a "reference" geometry that is used for the computation of the reduced basis functions. Therefore, the present methodology shares some common features with a standard reference-domain approach.\footnote{ At the continuous level, there are not particular advantages with respect to a standard reference-domain formulation. As highlighted in \autoref{sec:FVM} the main benefits are on the implementation side when the starting full order discretization is a finite volume one.}

\subsection{Discrete Empirical Interpolation Method}
As illustrated in the previous section, given the fact that it is possible to efficiently solve the mesh motion problem, also during the online stage (more details on this issue are reported in \autoref{subsec:lapl_red}), \autoref{eq:reduced_problem} still needs the assembly of a full order differential operator. In fact, for a new value of the parameter, it is not possible to recover an affine decomposition of the differential operator $\bm{A}(\bm{\mu})$ and of the source term vector $\bm{f}(\bm{\mu})$. The approach used here consists into an approximate affine expansion of the differential operator $\bm{A}(\bm{\mu})$ and of the source term vector $\bm{f}(\bm{\mu})$ as:
\begin{gather}
\bm{A}(\bm{\mu}) = \sum_{k=1}^{N_A} b^A_k(\bm{\mu})\bm{\chi}^A_k \mbox{ , } \bm{f}(\bm{\mu}) = \sum_{k=1}^{N_f} c^f_k(\bm{\mu})\bm{\chi}^f_k ,
\end{gather}
where $b^A_k(\bm{\mu})$, $c^f_k(\bm{\mu})$ and $\bm{\chi}^A_k \in \mathbb{R}^{N_h\times N_h}$, $\bm{\chi}^f_k \in \mathbb{R}^{N_h}$ are parameter dependent coefficients and parameter independent basis functions for the discretized differential operator $\bm{A}$ and for the source term vector $\bm{f}$ respectively. 

For the computation of the basis functions $\bm{\chi}^A_k$ and $\bm{\chi}^f_k$ different approaches are possible, here we used a matrix version of the snapshot POD method. The details of the computational procedure are reported in \autoref{alg:MPOD}. The starting point of the algorithm are the matrices $\bm{\mathcal{S}_A} = \{{\bm{A}(\bm \kappa_{i_{train}})}\}_{i=1}^{N_{train}} \in \mathbb{R}^{N_h \times N_h \times N_{train}}$ and vectors $\bm{\mathcal{S}_f} = \{{\bm{f}(\bm \kappa_{i_{train}})}\}_{i=1}^{N_{train}} \in \mathbb{R}^{N_h \times N_{train}}$. These need to be stored during the offline stage together with the full order model solutions. On the matrix and vector ``snapshots'' we apply then POD with the method of snapshots. As inner product to compute the correlation matrices we used the Frobenious inner product defined as $\langle \bm{A},\bm{B} \rangle_F = \Tr(\bm{A}^T \bm{B})$. 

\begin{algorithm}
\caption{THE POD algorithm for the D-EIM basis creation}
\label{alg:MPOD}
\hspace*{\algorithmicindent} \textbf{Input:} matrices $\bm{\mathcal{S}_A} = \{{\bm{A}(\bm \kappa_{i_{train}})}\}_{i=1}^{N_{train}}$ , vectors $\bm{\mathcal{S}_f} = \{{\bm{f}(\bm \kappa_{i_{train}})}\}_{i=1}^{N_{train}}$ , $N_A$, $N_f$\\
\hspace*{\algorithmicindent} \textbf{Output:} bases $\{\bm{\chi}_i^A\}_{i=1}^{N_A}$, $\{\bm{\chi}_i^f\}_{i=1}^{N_f}$
\begin{algorithmic}[1]
\STATE Compute the correlation matrices $C^A_{ij} = \frob{\bm{A}(\bm \kappa_{i_{train}})}{\bm{A}(\bm \kappa_{i_{train}})}$ and $C^f_{ij} = \frob{\bm{f}(\bm \kappa_{i_{train}})}{\bm{f}(\bm \kappa_{j_{train}})}$
\STATE Perform the eigenvalue decompositions $\bm{C^A}\bm{Q^A} = \bm{Q^A} \bm{\lambda^A}$ and $\bm{C^f}\bm{Q^f} = \bm{Q^f} \bm{\lambda^f}$
\STATE compute $\bm{\chi}_i^A = \frac{1}{N_{train} \lambda^A_i}\sum_{j=1}^{N_{train}}\bm{A}(\bm \kappa_{j_{train}})Q^A_{ij} \mbox{ }$ for $i=1:N_A$
\STATE compute $\bm{\chi}_i^f = \frac{1}{N_{train} \lambda^f_i}\sum_{j=1}^{N_{train}}\bm{f}(\bm \kappa_{j_{train}})Q^f_{ij} \mbox{ }$ for $i=1:N_f
$
\end{algorithmic}
\end{algorithm}

Once the basis functions are set, one needs to define a way to compute the coefficients vectors $\bm{b^A}(\bm{\mu}) \in \mathbb{R}^{N_A}$ and $\bm{c^f}(\bm{\mu}) \in \mathbb{R}^{N_f}$. To construct the above set of basis functions we rely on the standard discrete empirical interpolation procedure as introduced in \cite{Chaturantabut2010} for parametrized functions and extended in \cite{Negri2015} to parametrized differential operators:
\begin{gather}
\bm{b}^A(\bm{\mu}) = \bm{B}_A^{-1}\bm{A}_I(\bm{\mu}),\\
\bm{A}_{I_i}(\bm{\mu}) = \sum_{k=1}^{N_h} \bm{P}^T_{A_{k \bullet i}} \bm{A}(\bm{\mu})_{{\bullet k}},\\
\bm{c^f}(\bm{\mu}) = (\bm{P}_f^T\bm{U}_f)\bm{f}_I(\bm{\mu}),\\
\bm{f}_I(\bm{\mu}) = \bm{P}^T_f \bm{f}(\bm{\mu}),
\end{gather}
where the terms $\bm{B}_A \in \mathbb{R}^{N_h \times N_h \times N_A}$ , $\bm{P}_A \in \mathbb{R}^{N_h \times N_h \times N_A}$, $\bm{U}_f \in \mathbb{R}^ {N_h \times N_f}$ and $\bm{P}_f \in \mathbb{R}^{N_h \times N_f}$ are obtained following \autoref{alg:DEIM}. In the above expression we have introduced the $\Box_{\bullet}$ notation meaning that we let free the indices relative to the position of the $\bullet$ symbol. For example $\bm{P}_{A_{k\bullet i}}$ means we are considering the $k$-row of the $i$-layer of the $\bm{P}_A$ three-dimensional matrix. The $\bm{P}_A$ and $\bm{P}_f$ matrices assume the same meaning as those in \cite{Chaturantabut2010} with the difference that, working with parametrized differential operators, the $\bm{P}_A$ term is a three-dimensional matrix. Each $l$-layer of the $\bm{P}_A$ matrix is given by a matrix $\bm{E}^A_{l} \in \mathbb{R}^{N_h \times N_h}$ that is a matrix where the only non-zero element, which has unitary value, is located in the $(x,y)$ position identified by the magic point with coordinates $(I^A_{rl} \mbox{ } I^A_{cl})$  obtained at the $l$-iteration of the D-EIM \autoref{alg:DEIM}. Each $k$-column of $\bm{P}_f$ is given by the vector $\bm{E}^f_{k} \in \mathbb{R}^{N_h}$ which is a vector where the only non-zero element, which has also unitary value, is located at the position identified by the magic point $I^f_k$ obtained at the $k$-iteration of the D-EIM \autoref{alg:DEIM}.
\begin{algorithm}
\caption{The D-EIM procedure}
\label{alg:DEIM}
\hspace*{\algorithmicindent} \textbf{Input:} $\{\bm{\chi}_i^A\}_{i=1}^{N_A} \in \mathbb{R}^{N_h\times N_h\times N_A}$, $\{\bm{\chi}_i^f\}_{i=1}^{N_f} \in \mathbb{R}^{N_h \times N_f}$ \\
\hspace*{\algorithmicindent} \textbf{Output:} $\bm{I_A} \in \mathbb{R}^{N_h \times 2}$ , $\bm{I_f} \in \mathbb{R}^{N_h}$
\begin{algorithmic}[1]
  \STATE $\left[|\rho_A|,I^A_{r1},I^A_{c1}\right] = \max \{ | \bm{\chi}_1^A | \}$ , $\bm{I^A_1} = [I^A_{r1} \mbox{ } I^A_{c1}]$ and 
  $\left[|\rho^f|,I_1^f\right] = \max \{ | \bm{\chi}_1^f | \}$
  \STATE $\bm{E}^A_{1} = \bm{0}_{N_h,N_h}$, $\bm{E}^A_{1}(I^A_{r1},I^A_{c1}) = 1$, $\bm{E}^f_{1} = \bm{0}_{N_h,1}$, $\bm{E}^f_{1}(I^f_1) = 1$
  \STATE $\bm{U_A} = [\bm{\chi}_1^A] $ , $\bm{P_A} = [\bm{E}^A_{I_1}]$ , $\bm{I_A} = [\bm{I^A_{1}}]$ and 
  $\bm{U_f} = [\bm{\chi}_1^f] $ , $\bm{P_f} = [\bm{E}^f_{1}]$ , $\bm{I_f} = [I^f_{1}]$
  \FOR {$l = 2 \mbox{ to } N_A$ and  $k = 2 \mbox{ to } N_f $} 
  \STATE {Solve $(\bm{P_A^T U_A)\bm{c_A} = P_A^T \bm{\chi^A_l}}$ for $\bm{c_A}$} and 
  {Solve $(\bm{P_f^T U_f)\bm{c_f} = P_f^T \bm{\chi^f_k}}$ for $\bm{c_f}$}.
  \STATE {$\bm{r_A} = \bm{\chi^A_l} - \bm{U_A c_A}$} and
  {$\bm{r_f} = \bm{\chi^f_k} - \bm{U_f c_f}$} 
  \STATE { $\left[|\rho^A|,I^A_{rl},I^A_{cl}\right] = \max \{ | \bm{r_A} | \}$} ,  $\bm{I^A_l} = [I^A_{rl} \mbox{ } I^A_{cl}]$ and 
  { $\left[|\rho^A|,I^f_{k}\right] = \max \{ | \bm{r_f} | \}$} 
  \STATE $\bm{E}^A_{l} = \bm{0}_{N_h,N_h}$, $\bm{E}^A_{l}(I^A_{rl},I^A_{cl}) = 1$, $\bm{E}^f_{k} = \bm{0}_{N_h,1}$, $\bm{E}^f_{k}(I^f_k) = 1$
  \STATE {$\bm{U_A} \gets \left[ \bm{U_A} \mbox{ } \bm{\chi^A_l} \right]$ } , 
  {$\bm{P_A} \gets \left[ \bm{P_A} \mbox{  } \bm{E}^A_{l} \right]$ } , {$\bm{I_A} \gets 
  \begin{bmatrix}
  \bm{I_A}\\
  \bm{I^A_l}
  \end{bmatrix}
  $},\\
  {$\bm{U_f} \gets \left[ \bm{U_f} \mbox{ } \bm{\chi^f_k} \right]$ } , 
  {$\bm{P_f} \gets \left[ \bm{P_f} \mbox{ } \bm{E}^f_{k} \right]$ } , {$\bm{I_f} \gets 
  \begin{bmatrix}
  \bm{I_f}\\
  {I^f_k}
  \end{bmatrix}
  $} \\
  \ENDFOR
\end{algorithmic}
\end{algorithm}

Finally we can write the D-EIM approximation as:
\begin{gather}
\bm{A}(\bm{\mu}) \approx \sum_{k=1} \bm{D}_{A_{\bullet \bullet k}} \bm{A}_{I_k}(\bm{\mu}), \\ 
\bm{D}_{A_{\bullet \bullet i}} = \sum_{k=1}^{N_A}\bm{U}_{A_{\bullet \bullet k}} \bm{B}^{-1}_{A_{ki}}, \\
\bm{B}_{A_{ij}} = \sum_{k=1}^{N_h} \bm{P}^T_{A_{k \bullet i}} \bm{U}_{A_{\bullet k j}}, \\
\bm{f}(\bm{\mu}) \approx \bm{B_f} \bm{f}_I(\bm{\mu}),\\
\bm{B_{f} = \bm{U}_f(\bm{P}_f^T \bm{U}_f)^{-1} \in \mathbb{R}^{N_h\times N_f}},\\
\bm{f}_I(\bm{\mu}) = \bm{P}_f^T \bm{f}(\bm{\mu}) \in \mathbb{R}^{N_f} ,
\end{gather}
where the D-EIM approximation for the vector quantities is reported with the original notation presented in \cite{Chaturantabut2010} while the one for D-EIM approximation of the discretized differential operator $\bm{A}(\bm{\mu})$, since we are working with three-dimensional sparse matrices, and we want to rely on sparse matrices linear algebra operations, the notation is slightly modified to account for three-dimensional matrices. In the above formulations we have introduced the additional variables $\bm{D_A} \in \mathbb{R}^{N_h \times N_h \times N_A}$.

In the end, the idea is to precompute during the offline stage the terms $\bm{A}^r_k \in \mathbb{R}^{N_r \times N_r}$ and  $\bm{f}^r_k \in \mathbb{R}^{N_r}$ and to express reduced operator and vector as:
\begin{gather}
\bm{A}^r(\bm{\mu}) = \sum_{k=1}^{N_A} \bm{A}^r_k \bm{A}_{I_k}(\bm{\mu}), \quad \bm{f}^r(\bm{\mu}) = \sum_{k=1}^{N_f} \bm{f}^r_k \bm{f}_{I_k}(\bm{\mu}),
% \bm{A}^r_i = \bm{L}^T \bm{\chi}^A_i\bm{L}, \quad \bm{f}^r_i = \bm{L}^T \bm{\chi}^f_i.
\end{gather}
where:
\begin{gather}
\bm{A}_k^r = \bm{L}^T \bm{D}_{\bullet\bullet k} \bm{L}, \quad \bm{f}_k^r(\bm{\mu})= \bm{L}^T \bm{B}_{\bullet k}.
\end{gather}
The terms $\bm{A}^r_k$ and $\bm{f}^r_k$ are precomputed during the offline stage and $\bm{A}_{I}(\bm{\mu})$ and $\bm{f}_{I}(\bm{\mu})$ are vector of coefficients that must be computed during the online stage. 

In the above formulations it is important to notice that the $\bm{A}_I$ and $\bm{f}_I$ vector of coefficients correspond to point-wise evaluations of the differential operator $\bm{A}(\bm{\mu})$ and of the source term vector $\bm{f}$ in the locations corresponding to the magic points inside the matrix $\bm{I}_A$ and the vector $\bm{I}_f$ reported in \autoref{alg:DEIM}.

\begin{figure}
\centering
\begin{minipage}{0.8\textwidth}
\centering
\begin{minipage}{0.48\textwidth}
\centering
\includegraphics[height=0.6\textwidth]{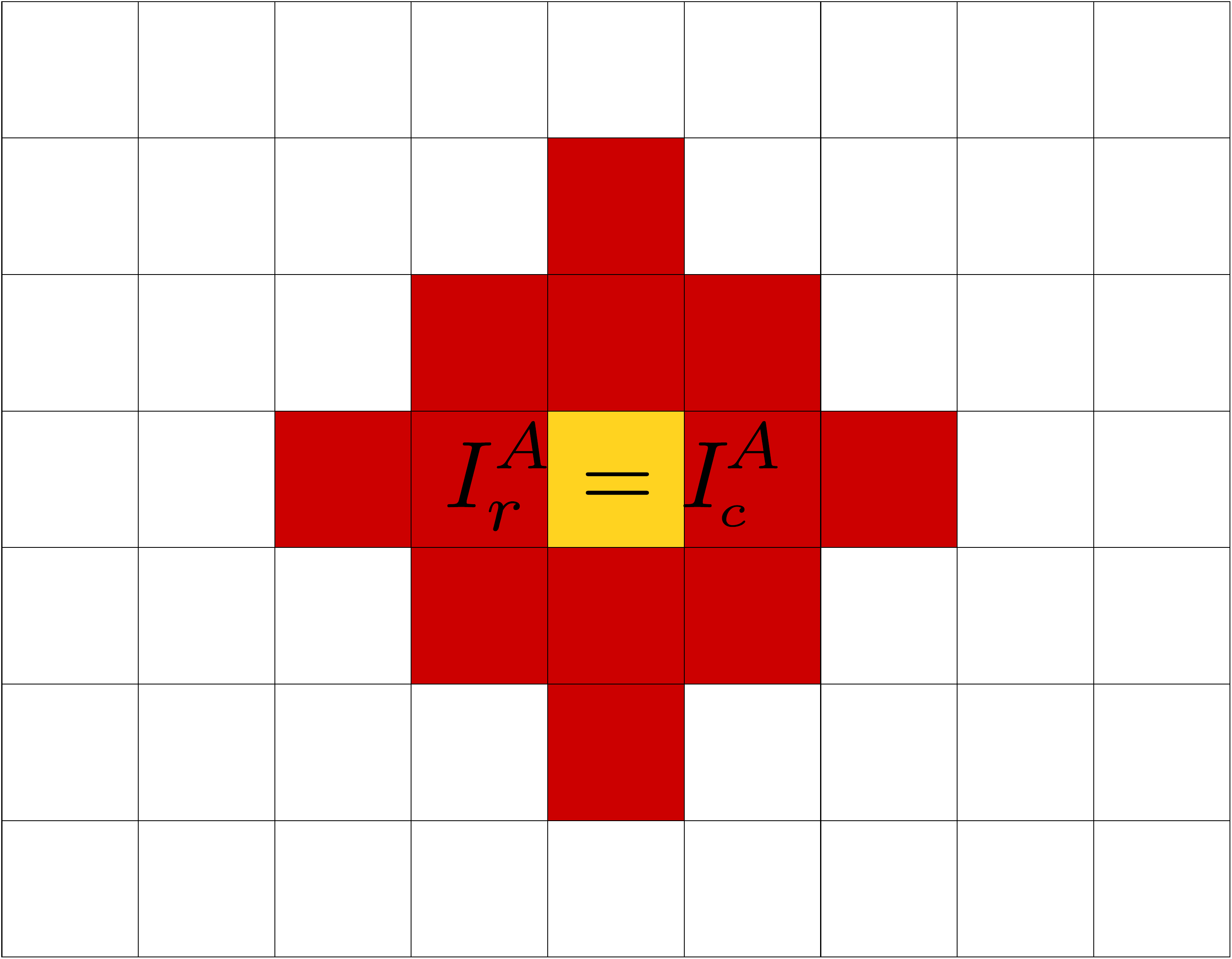}
\end{minipage} 
\begin{minipage}{0.48\textwidth}
\centering
\includegraphics[height=0.6\textwidth]{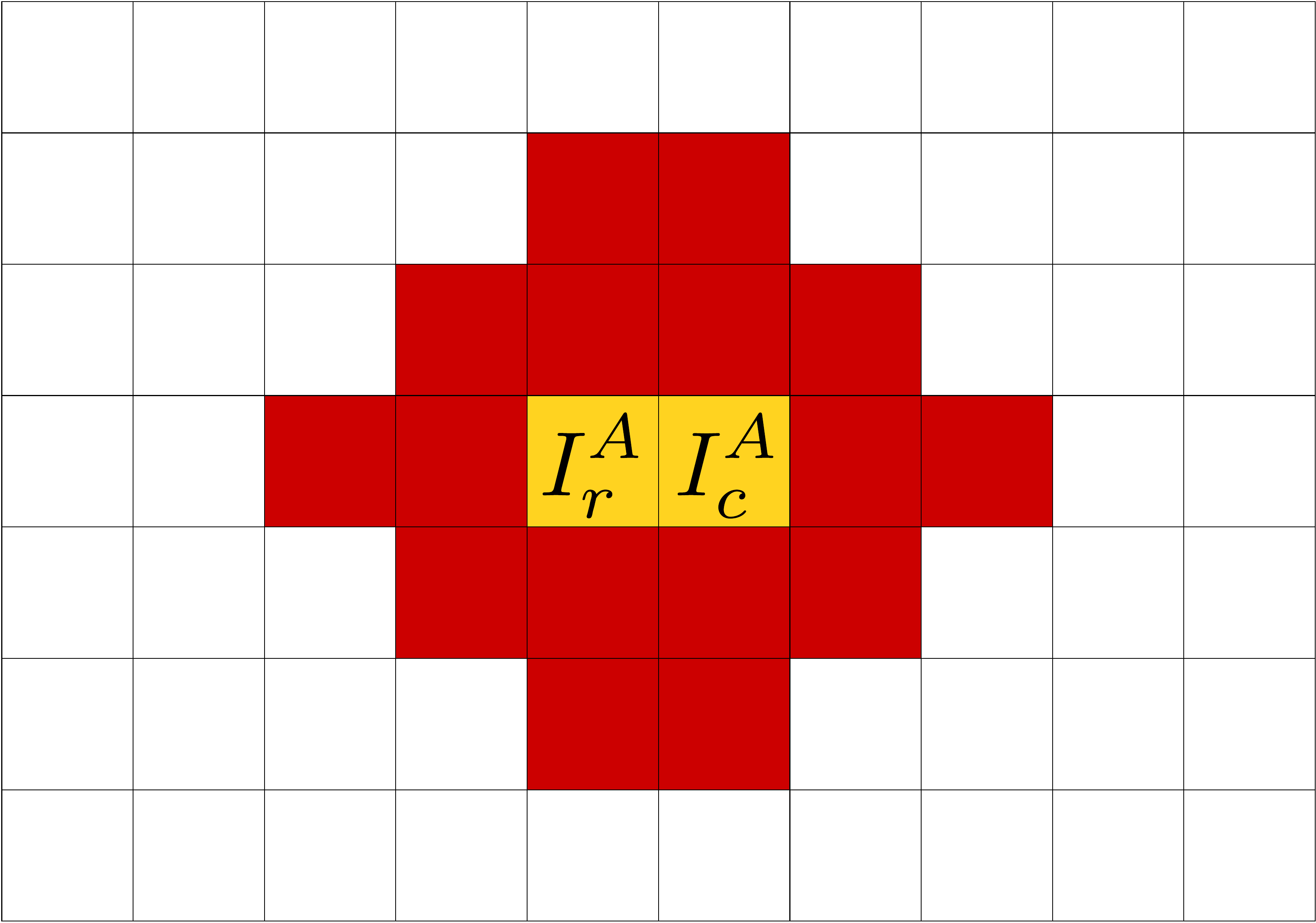}
\end{minipage} 
\end{minipage} 
\caption{Example of the computational stencil needed for the evaluation of the D-EIM coefficients during the online stage. On the left there is an example of a point on the diagonal of the discretized differential operator, while on the right there is an example of a point outside the diagonal of the differential operator. In the image two layers of cells are considered. The yellow cells depict the ones identified during the D-EIM procedure by the indices $I^A_r$ and $I^A_c$ while the red ones depict the additional cells required for the differential operator evaluation.}
\label{fig:comp_mol}
\end{figure}

Therefore, the only additional complexity, during the online stage, corresponds into point-wise evaluation of the source term vector and of the discretized differential operator. In order to make the procedure efficient and independent with respect to mesh size, this point-wise evaluation must involve operations only on a subset of the entire domain and must be independent with respect to the total number of degrees of freedom of the full order problem. The finite volume method satisfies this property, in fact, the point evaluation of an operator can be computed also on a subset of the whole domain. Each subset includes the cells relative to the point evaluation and a certain number of neighboring layers of cells. The required number of neighboring layers of cells depends on the type of scheme employed for the computation of the operator. For example a central differencing scheme, without orthogonal correction, requires only one layer of cells while for second order schemes or a central differencing scheme with orthogonal correction at least two layers of cells are necessary. The process for the construction of the subsets is exemplified in \autoref{fig:comp_mol}. In the figure are reported two examples of the required computational stencil for both a point on the diagonal of the discretized differential operator and a point outside of the diagonal of the discretized differential operator. 

\section{A parametrized heat transfer problem}\label{sec:num_ex}
\begin{figure}[H]
\centering
\begin{minipage}{\textwidth}
\centering
\begin{minipage}{0.24\textwidth}
\centering
\includegraphics[width=\textwidth]{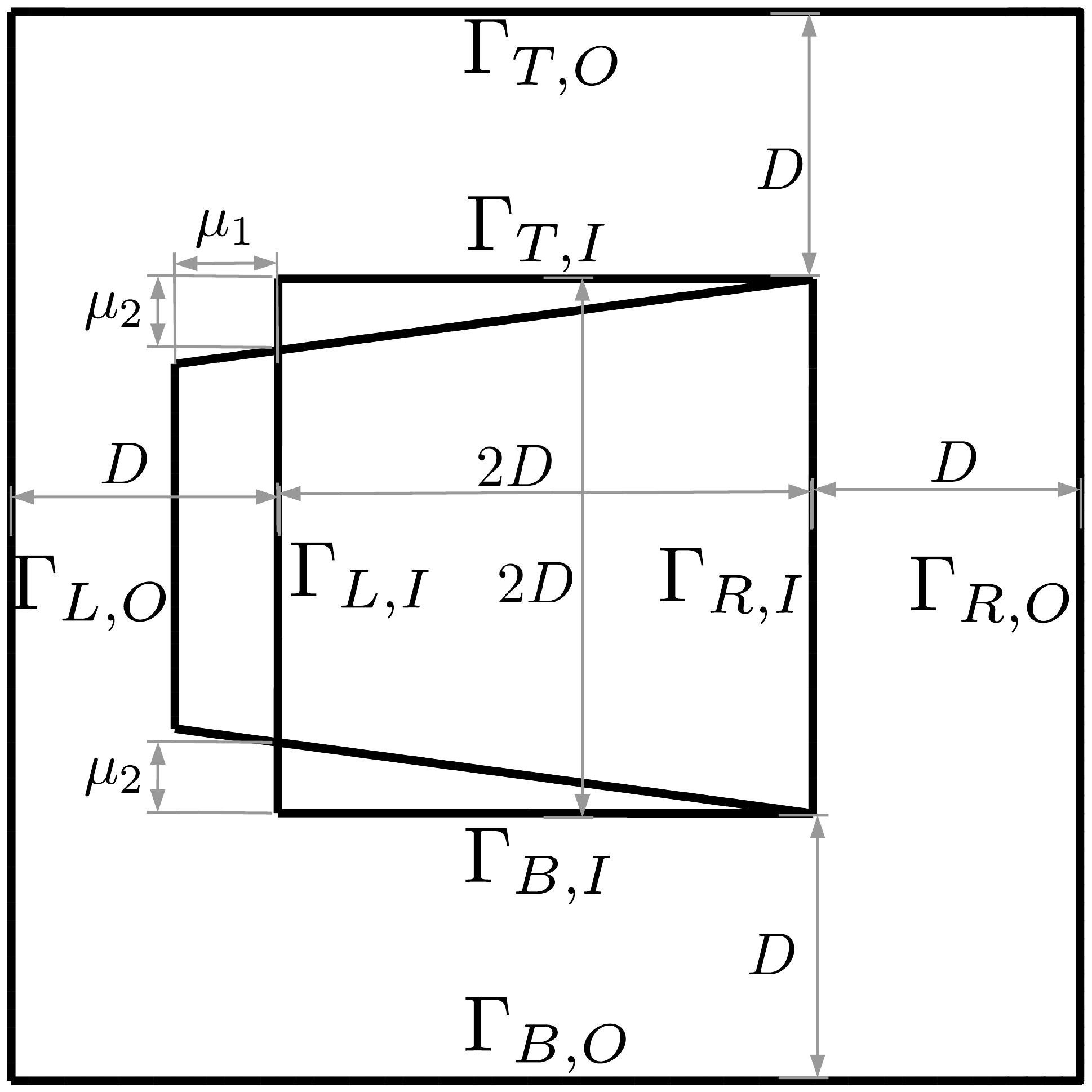}
\end{minipage} 
\begin{minipage}{0.24\textwidth}
\centering
\includegraphics[width=\textwidth]{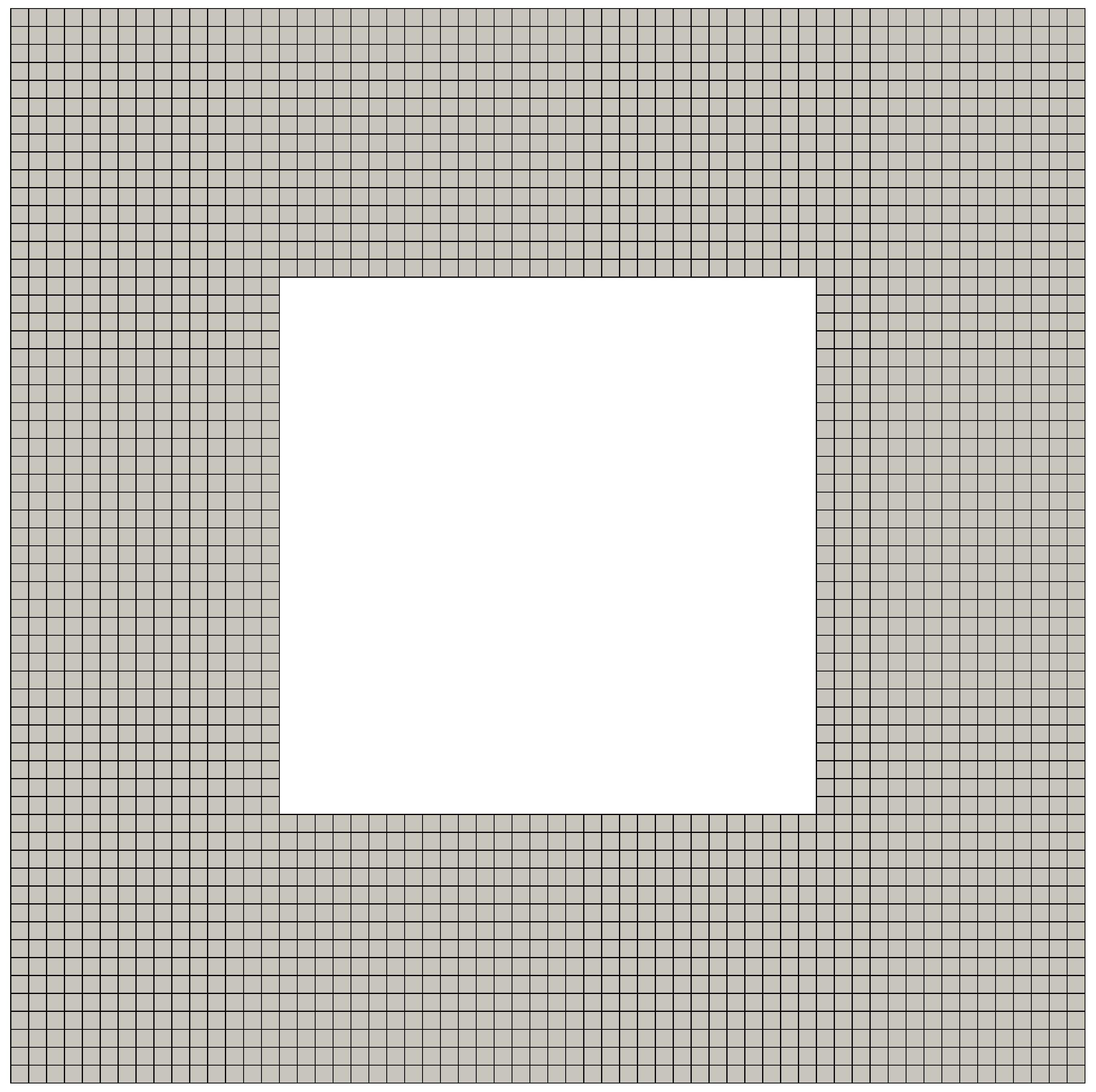}
\end{minipage} 
\begin{minipage}{0.24\textwidth}
\centering
\includegraphics[width=\textwidth]{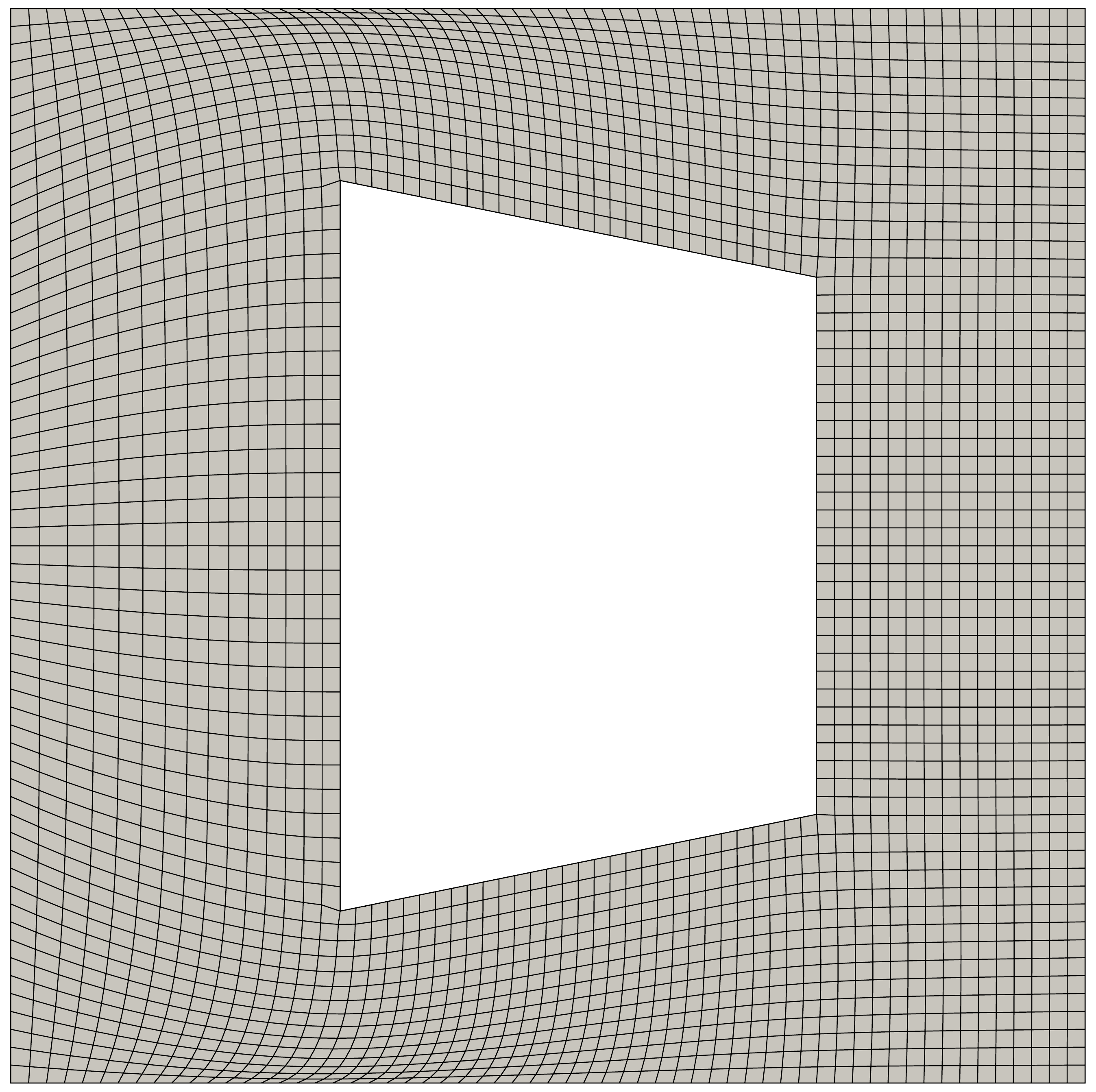}
\end{minipage} 
\begin{minipage}{0.24\textwidth}
\centering
\includegraphics[width=\textwidth]{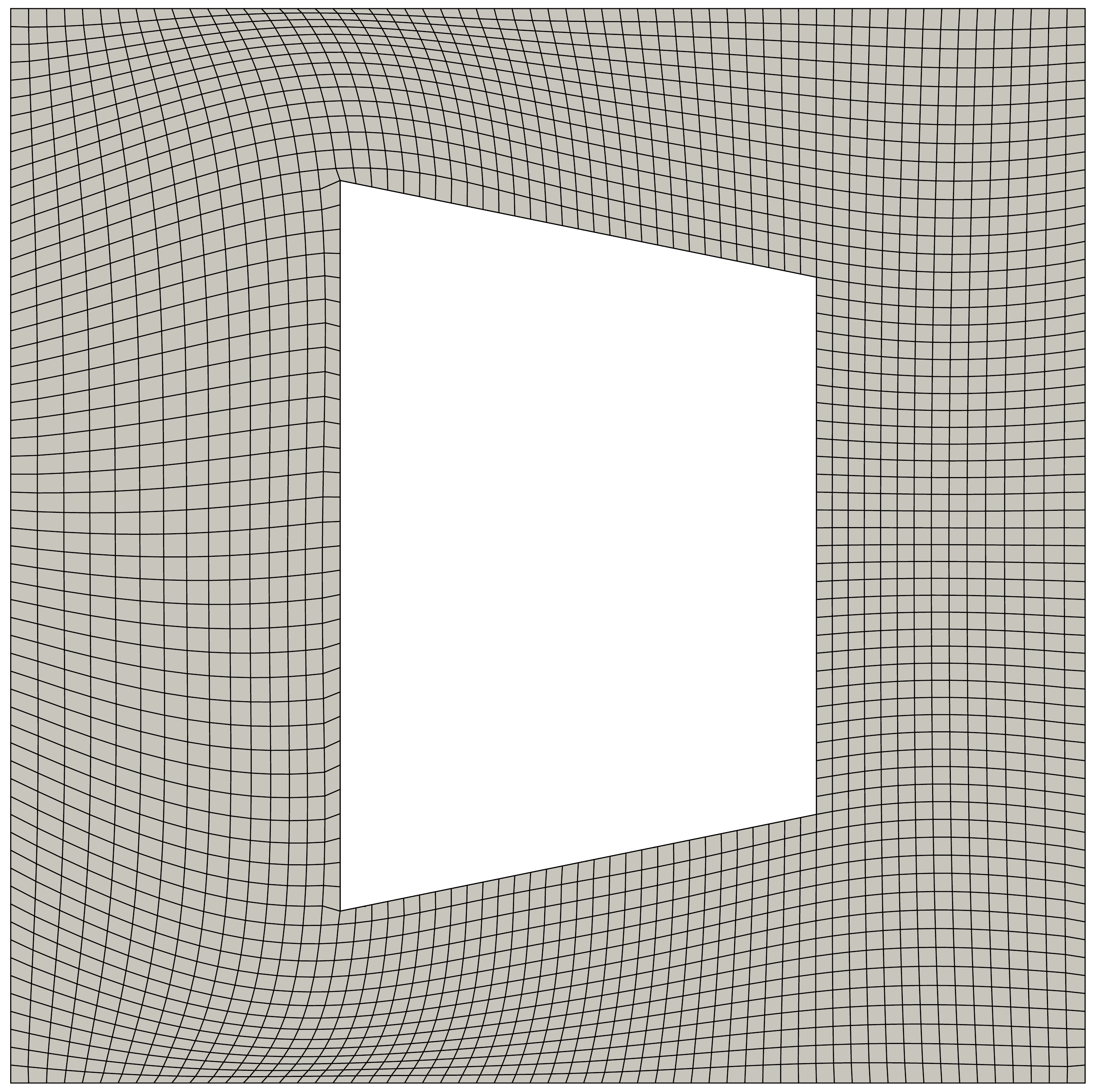}
\end{minipage} 
\end{minipage} 
\caption{From left to the right we can observe the domain together with the dimensions ($D = 1$) and the parameters scheme, the tessellation in its undeformed configuration, the tessellation deformed using a Laplacian smoothing approach and the tessellation deformed using a radial basis function approach.}
\label{fig:domains}
\end{figure}

In this section we present the proposed methodology on a parametrized heat transfer problem. Particular attention is posed onto the selection of the best mesh motion algorithm especially in view of model reduction purposes. The domain is parametrized by the parameter vector $\bm{\mu} = (\bm{\mu}_1,\bm{\mu}_2) \in [\bm{\mu}_{1,min},\bm{\mu}_{1,max}] \times[\bm{\mu}_{2,min},\bm{\mu}_{2,max}] \subset \mathbb{R}^2 $. The parameters scheme is depicted on the left of \autoref{fig:domains}. In the same figure it is reported also the chosen tessellation in its undeformed configuration which accounts for $2700$ quadrilateral cells. In the same figure it is also depicted the mesh configuration using both the Laplacian smoothing and the radial basis function approach. The numerical tests are conducted using the in-house open source library \emph{ITHACA-FV} (In real Time Highly Advanced Computational Applications for Finite Volumes) \cite{RoSta17} which is a computational library based on the finite volume solver OpenFOAM 6.0 \cite{OF}. The computational domain is given by a square block with a hole inside. The boundary $\Gamma = \{\Gamma_{T,O},\Gamma_{R,O},\Gamma_{B,O},\Gamma_{L,O},\Gamma_{T,I},\Gamma_{R,I},\Gamma_{B,I},\Gamma_{L,I}\}$ is subdivided into eight different parts which consists of the top, right, bottom and left sides of the outer and inner polygons. The boundary conditions are set according to \autoref{tab:BC}.

\begin{table}
\centering
\begin{tabular}{ l | c | c | c | c | c | c | c | c }
      & $\Gamma_{L,O}$ & $\Gamma_{T,O}$ & $\Gamma_{R,O}$ & $\Gamma_{B,O}$ & $\Gamma_{L,I}$ & $\Gamma_{T,I}$ & $\Gamma_{R,I}$ & $\Gamma_{B,I}$ \\
      \hline
  B.C. & $\theta_{\bm n} = 10$ & $\theta = 10$ & $\theta_{\bm n} = - 10$ & $\theta = 10$ & $\theta = -10$ & $\theta = 0$ & $\theta = 0$ & $\theta = 0$ \\
\end{tabular}
\caption{Imposed values of the Dirichlet and Neumann boundary conditions for the different boundaries of the domain}\label{tab:BC}%
\end{table}

Before going deeply into numerical solutions, let us report here the strong form of the steady state heat transfer equations and describe briefly what is the expected behavior of the physical solution for the problem we took into account:
\begin{equation}\label{eq:poisson_problem_2}
\begin{cases}
\mbox{div}(\alpha_\theta \nabla \theta) = f &\mbox{ in } {\Omega}({\bm{\mu}}),\\
\theta(x,\bm{\mu}) = \theta_D(x,\bm{\mu}) &\mbox{ on } {\Gamma_D(\bm{\mu})},\\
\theta_{\bm{n}}(x,\bm{\mu}) = \theta_N(\bm{\mu},x) &\mbox{ on } {\Gamma_N(\bm{\mu})},
\end{cases}
\end{equation}
where $\alpha_\theta$ is the diffusivity of the homogeneous isotropic medium and has been fixed equal to $1$ while a null forcing term $f$ has been applied to the problem. 

This is the most classical choice when dealing with elliptic problems. The Laplace operator only acts like a smoothing-diffusing term. For this reason no peaks or discontinuities are expected into the domain, boundary temperatures should be spread allover $\Omega$ since the forcing term is adding no contribution. According with the \textit{maximum principle}, the maximum (and the minimum as well) of the solution should be located on the boundaries since $f=0$.

\subsection{The mesh motion strategies}\label{subsec:mms}
\begin{figure}
\begin{minipage}{\textwidth}
\centering
\begin{minipage}{0.4\textwidth}
\centering
\includegraphics[width = \textwidth]{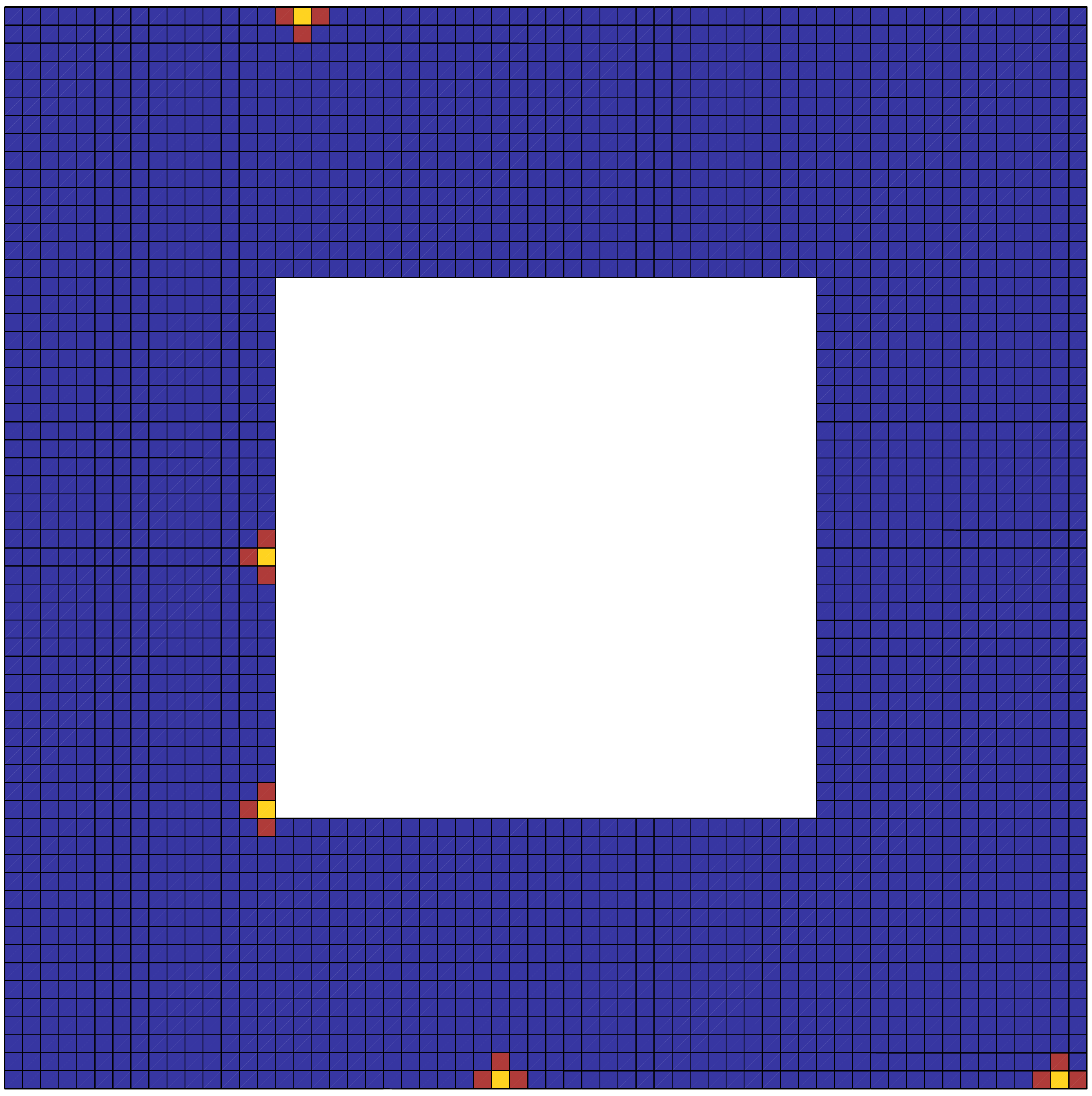}
\end{minipage}
\begin{minipage}{0.4\textwidth}
\centering
\includegraphics[width=\textwidth]{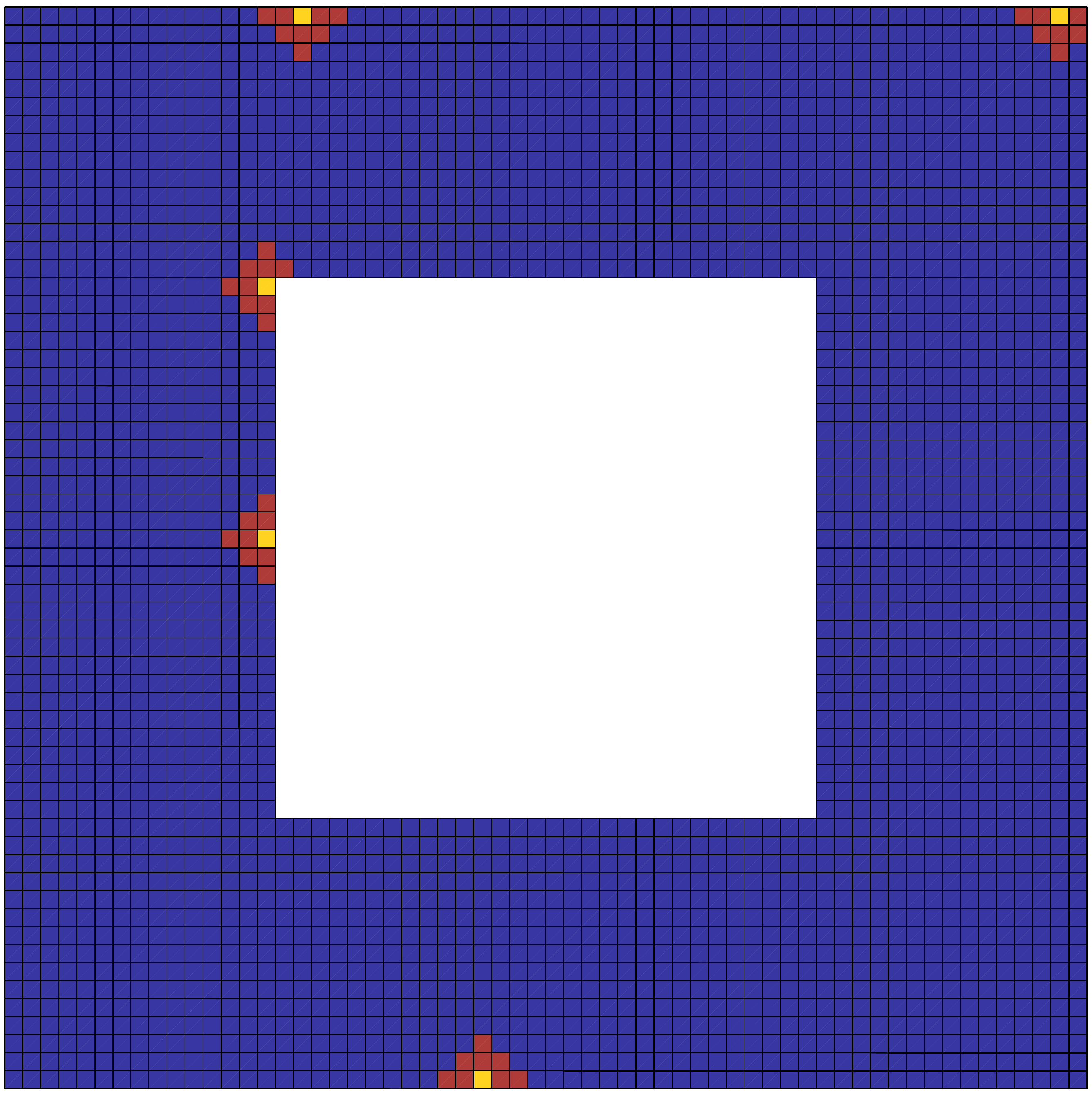}
\end{minipage}
\end{minipage}
\caption{The first $5$ magic points identified during the D-EIM procedure for the source term vector in case of a linear scheme without orthogonal correction and in the case of a linear scheme with orthogonal correction.}
\label{fig:DEIM_mesh}
\end{figure}

In the numerical example we examine different mesh motion strategies starting from different full order schemes. As reported in the previous section, and as exemplified in \autoref{fig:DEIM_mesh}, the selection of the full order numerical scheme used to deal with non-orthogonality affects the number of layers of additional cells around the ones identified by the magic points during the D-EIM procedure. In the picture we report, for a case with just $5$ D-EIM modes, the location of the cells identified by the magic points (yellow cells), together with the additional cells of the computational stencils necessary for the point-wise evaluation of the differential operators during the online stage (red cells). The identified cells are relative to the magic points identified for the source term vector $\bm{f}$. 

\begin{figure}
\begin{minipage}{\textwidth}
\begin{minipage}{0.2\textwidth}
\centering
LAPL.
\includegraphics[width = \textwidth]{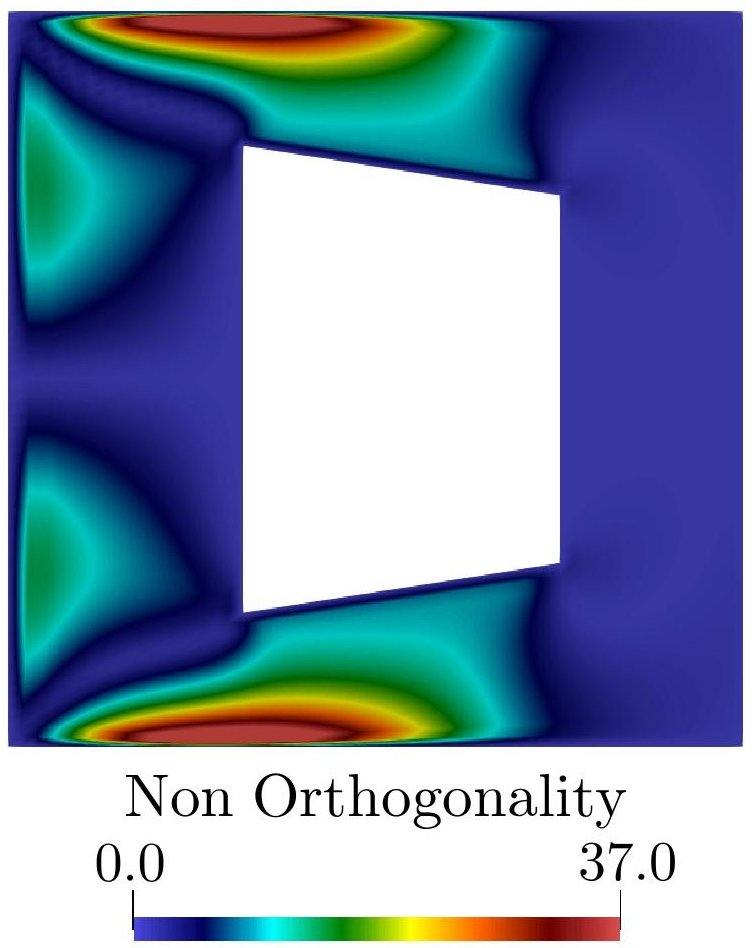}
\end{minipage}
\begin{minipage}{0.2\textwidth}
\centering
RBF
\includegraphics[width = \textwidth]{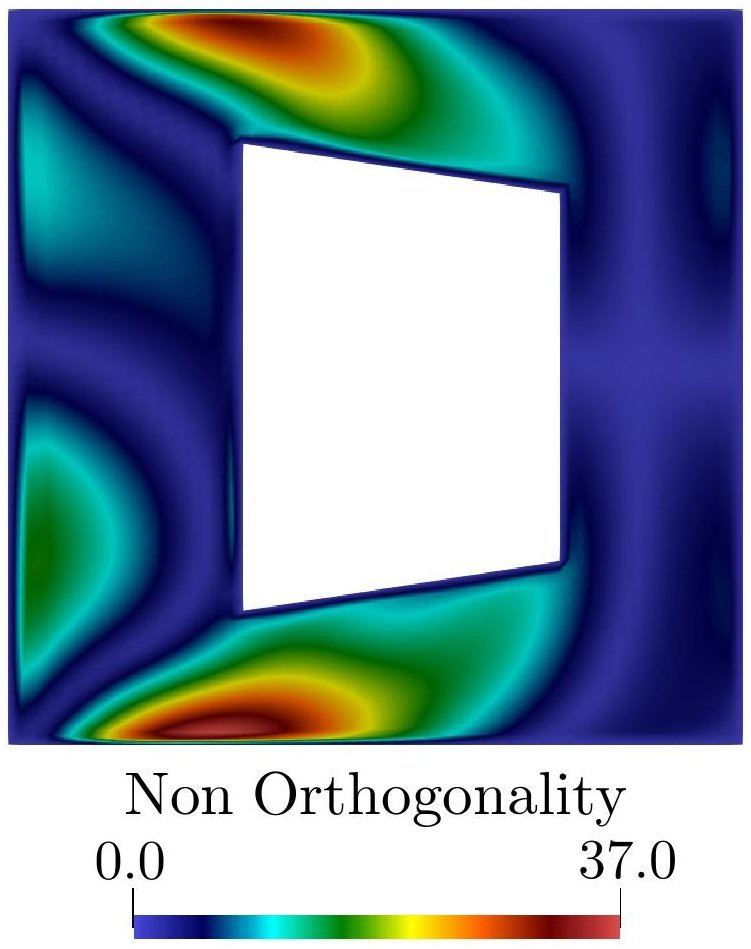}
\end{minipage}
\begin{minipage}{0.2\textwidth}
\centering
Control Points
\vspace{0.9cm}
\includegraphics[width = \textwidth]{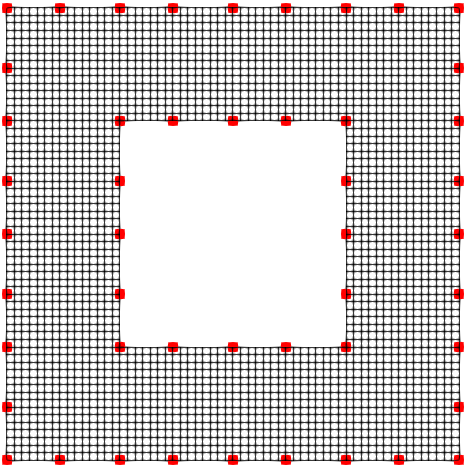}
\end{minipage}
\begin{minipage}{0.35\textwidth}
\resizebox{\textwidth}{!}{%
\begin{tabular}{l | r | r}
  & Max Non-Ortho. & Average Non-Ortho.\\
  \hline
  LAPL. & $\ang{29.774}$ & $\ang{7.724}$ \\
  \hline
  RBF & $\ang{25.518}$ & $\ang{7.531}$ 
\end{tabular}}
\end{minipage}
\end{minipage}
\caption{Non-Orthogonality comparison between a Laplacian smoothing and a RBF approach. The two images on the left show, for one selected sample value $\bm{\kappa*} = (0.181184,0.288162)$, the value of the non-orthogonality inside the domain. In the third image on the right are reported the control points identified by the automatic coarsening procedure. In the table it is reported the average value, over $100$ sample points, for the maximum value of the non-orthogonality and for the average value of the non-orthogonality.}
\label{Fig:non_ortho}
\end{figure}

\begin{figure}
\includegraphics[width = \textwidth]{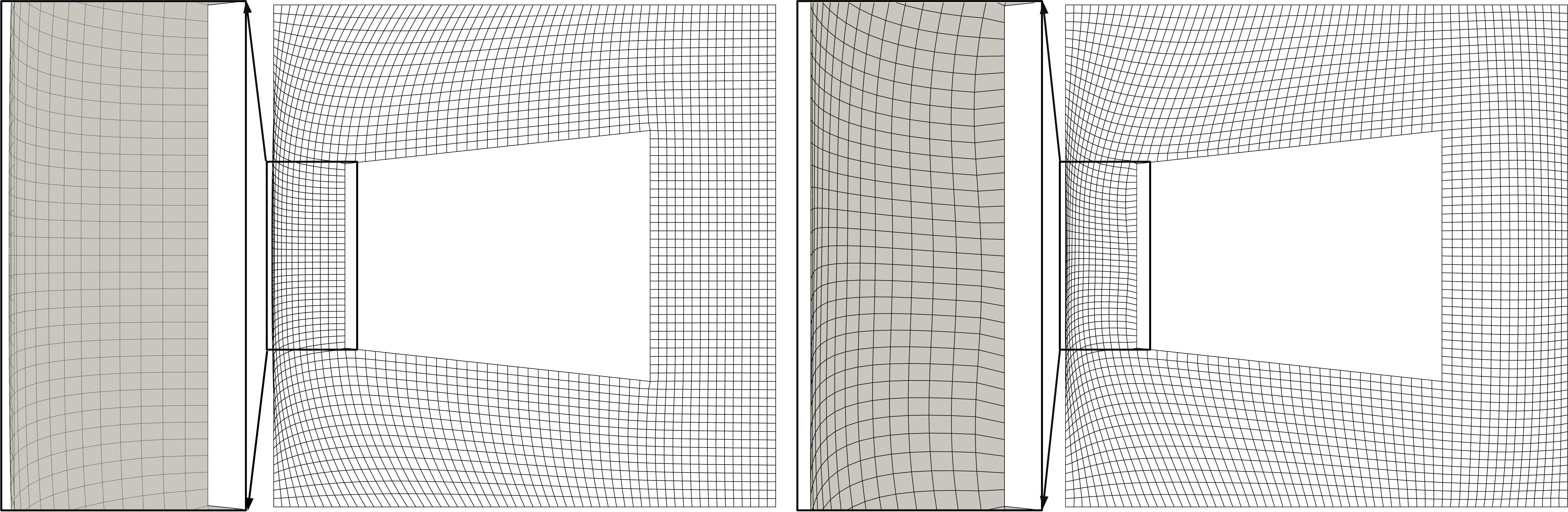}
\caption{Comparison of the two mesh deformation strategies for large parametric geometrical variations. On the left we observe the Laplacian smoothing approach that, as highlighted in the zoom, produces a deformed mesh that has cells with negative volumes (cells close to the left boundary). On the right we have the RBF approach.}
\label{Fig:large_def}
\end{figure}

Before analyzing the efficiency and the approximation properties of the reduced order model we analyze the performances of the mesh motion strategies at the full order level. As extensively illustrated in the previous section, when using ALE approaches to deal with parametrized geometries, the accuracy of the results depends strongly on the quality of the deformed mesh. In a finite volume context, one of the most important indicators to certify the quality of a mesh consists into the non-orthogonality factor. We analyze the parametrized geometry of \autoref{fig:domains} and we deform the mesh according to the discrete set of parameter samples $\mathcal{K}^1_{train} = \{\bm{\kappa}_{i_{train}}\}_{i=1}^{N_{train}} \in [-0.32,0.32] \times [-0.32,0.32]$. We set $N_{train}=100$ and parameter samples are selected randomly inside the parameter space. The results of such test are reported qualitatively and quantitatively in the images and in the table of \autoref{Fig:non_ortho}. As can be seen from the image and the table, the Laplacian smoothing approach, compared with the RBF strategy, produces a deformed mesh with higher values of the non-orthogonality. This fact is particularly evident in the part of the domain close to the upper and lower boundaries. Moreover, changing the dimension of the parameter space used to determine the training set to $\mathcal{K}^2_{train} \in [-0.45,0.45] \times [-0.45,0.45]$, as shown in \autoref{Fig:large_def}, for some values of the samples, the mesh motion strategy using a Laplacian smoothing approach produces a mesh where some of the cells have negative volumes. In these cases of course the simulation results become unreliable. The RBF mesh motion strategy requires the selection of $56$ control points on the moving boundaries and requires therefore, as illustrated in \autoref{subsec:RBF}, the resolution of relatively small linear system (especially if we consider the dimension of the full order problem $N_h = 2700$). The selection of the control points has been done automatically fixing a coarsening ratio which defines the ratio between the total number of points on the boundaries and the number of control points. The location of the control points identified by this automatic procedure is depicted in \autoref{Fig:non_ortho}. The radial basis functions are given by Gaussian functions with a radius $d_{RBF} = 0.6$.

\subsubsection{Reduction of the mesh motion problem in the case of Laplacian smoothing technique}\label{subsec:lapl_red}
The discretized operator $\bm A_D$ of \autoref{eq:mesh_lapl} is always assembled in the same physical domain (i.e. the undeformed one) and has a diffusivity value $\gamma$ which is changing in space but that is always constant with respect to the parameter values. Therefore, it can be assembled only once during the offline stage and projected onto the reduced basis space of the geometrical deformation field. For the source term $\bm b_D(\bm \mu)$, in the general case it is usually not possible to recover an affine decomposition and we have therefore to rely on an approximate affine expansion using the Discrete Empirical Interpolation Method:
\begin{equation}
\bm b_D (\bm \mu) \approx \sum_{i=1}^{N_{b_d}} a^D_i (\bm \mu) \chi^D_i.
\end{equation}
Once this approximation is done, it is straightforward to obtain a reduced order model of the mesh motion problem. For this particular case, since the motion of the moving boundaries depends linearly with respect to the input parameters, two D-EIM modes for the source term $\bm b_d$ were sufficient to obtain accurate results. The average value of the L2 relative error between the full and the reduced order solutions over the entire testing set is equal to $1.98\cdot 10^{-7}$. This results was obtained using only 2 POD modes for the mesh motion field and 2 D-EIM basis function for the source term vector.

%Comparing the computational time necessary to solve the mesh motion problem, for the whole set of samples, using the two different approaches, we notice comparable results. The Laplacian smoothing approach which entails the resolution of a larger but sparse system, requires a total time of $1.49 \mbox{s}$ while the RBF approach, that on the other side entails the resolution of a smaller dense system, requires $1.07 \mbox{s}$. Both cases run serially on a desktop machine with an Intel\textsuperscript{\tiny\textregistered} Core(TM) i7-7700 CPU @ 3.60GHz CPU.

% \begin{table}
% \centering
% \begin{tabular}{l | r | r}
%   & Max Non-Ortho. & Mean Non-Ortho.\\
%   \hline
%   RBF & $\ang{25.518}$ & $\ang{7.531}$ \\
%   \hline
%   LAPL. & $\ang{29.774}$ & $\ang{7.724}$ \\
% \end{tabular}
% \caption{Non-Orthogonality}
% \label{tab:non_ortho}
% \end{table}

\subsection{The reduced order model}

\begin{figure}[H]
\centering
\begin{minipage}{\textwidth}
\centering
\begin{minipage}{0.49\textwidth}
\centering
\includegraphics[width=\textwidth]{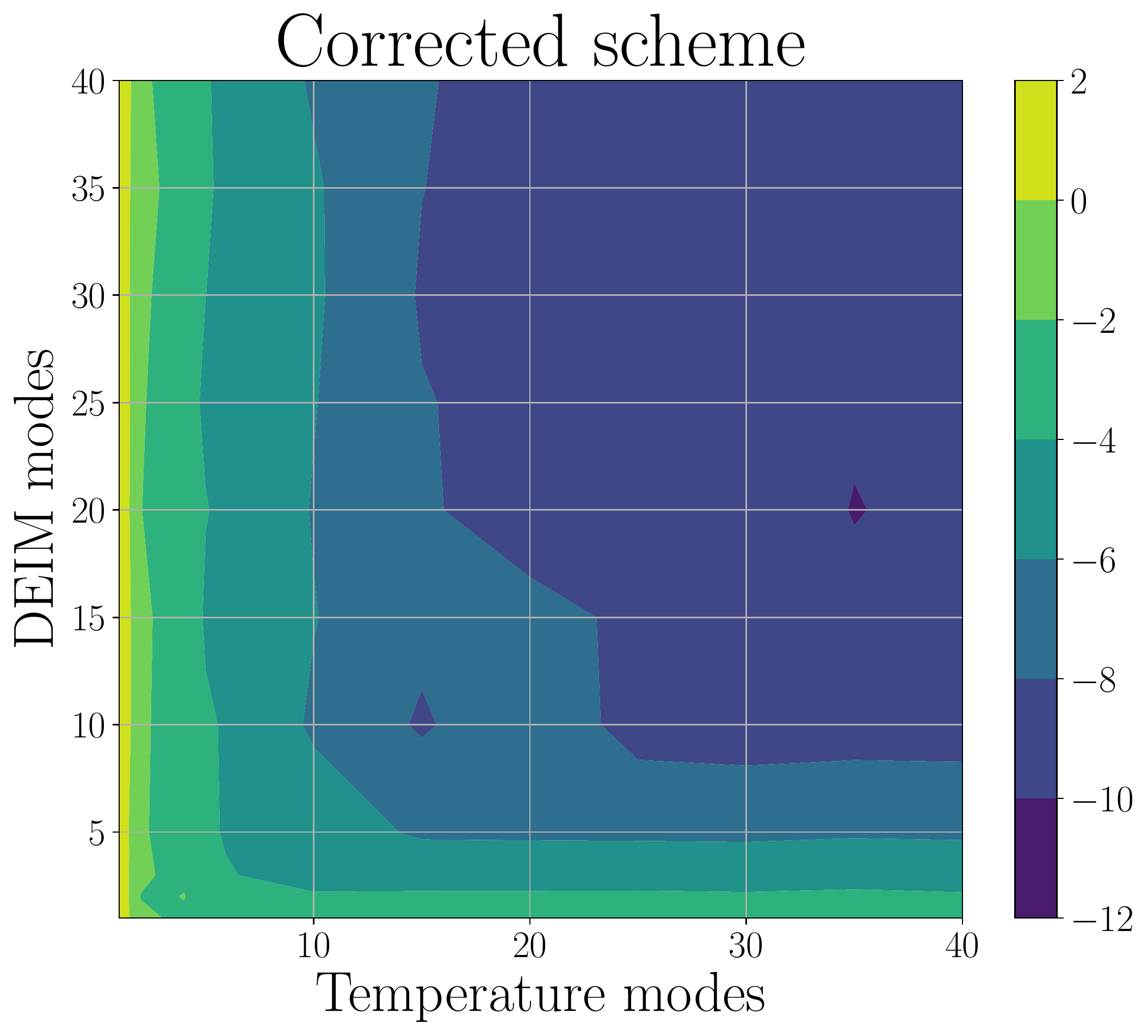}
\end{minipage}
\begin{minipage}{0.49\textwidth}
\centering
\includegraphics[width=\textwidth]{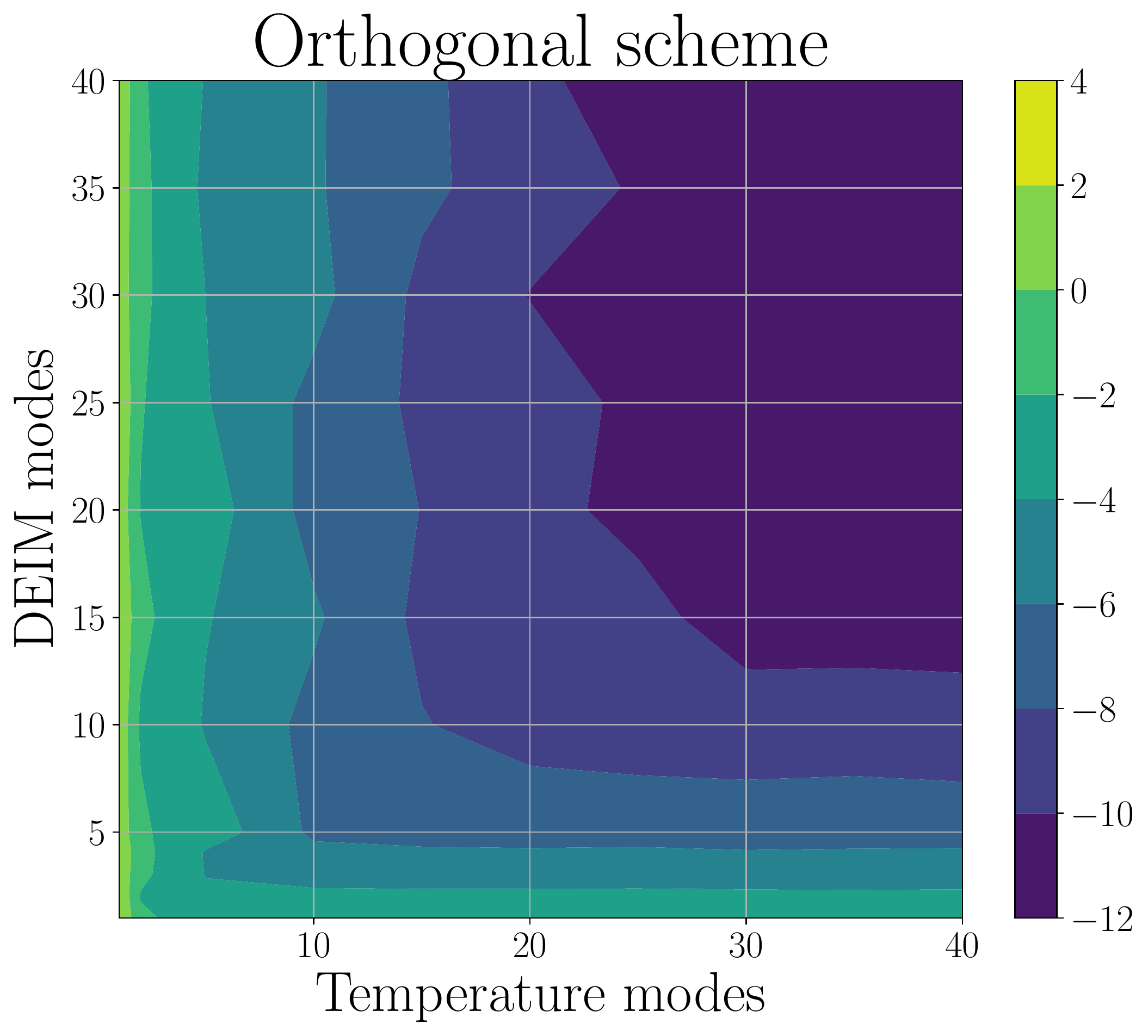}
\end{minipage} 
\end{minipage} 
\caption{The average of the $L^2$ relative error norm of the temperature field for the Laplacian smoothing mesh motion strategy. The logarithm of the $L^2$ relative error norm is plotted against the number of modes used for the temperature field and the number of D-EIM modes used to approximate the discretized differential operator $\bm A$ and the source term $\bm f$ ($N_{A} = N_f$). The plots are reported with (left) and without (right) non-orthogonal correction.}
\label{fig:Lapl_comp}
\end{figure}

\begin{figure}[H]
\centering
\begin{minipage}{\textwidth}
\centering
\begin{minipage}{0.49\textwidth}
\centering
\includegraphics[width=\textwidth]{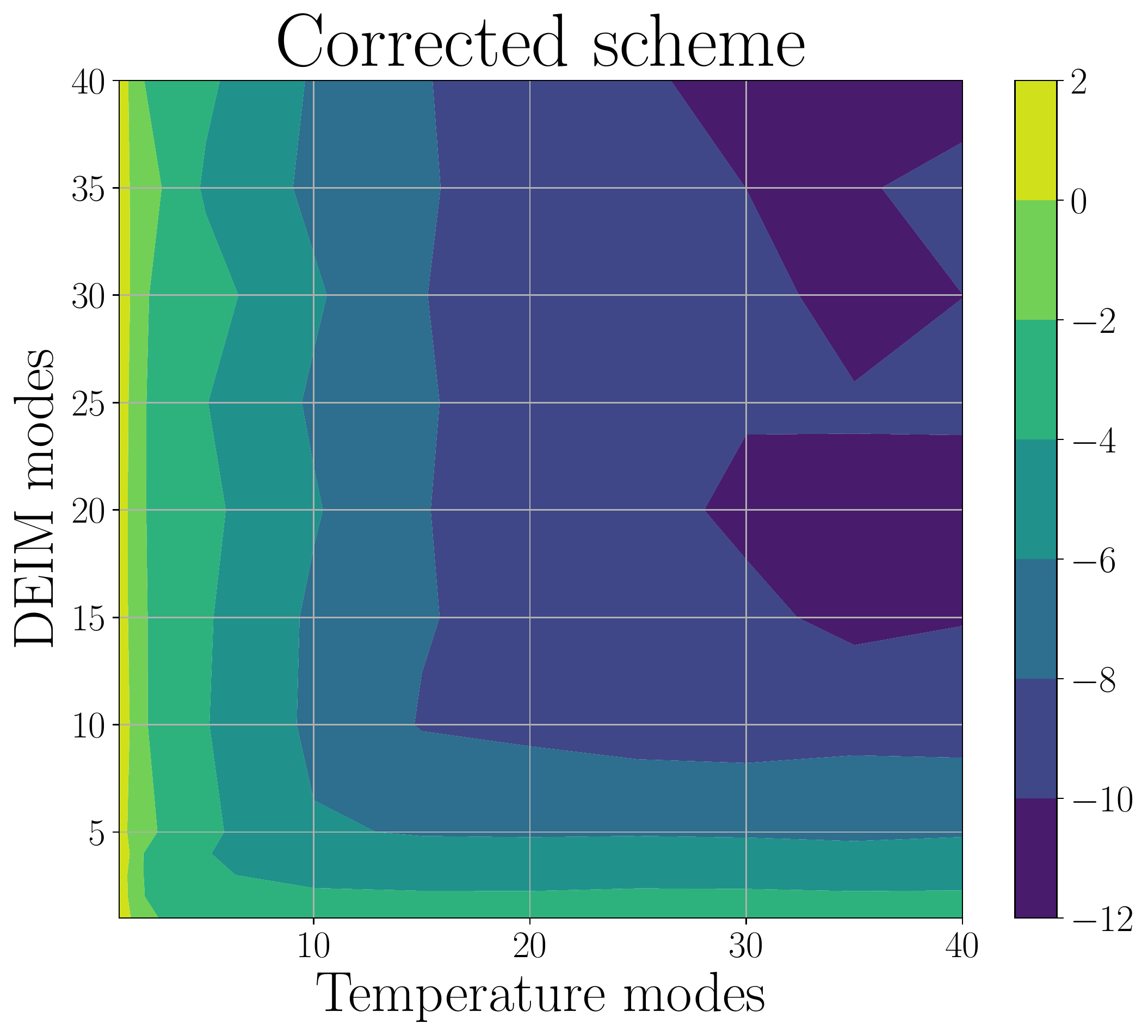}
\end{minipage} 
\begin{minipage}{0.49\textwidth}
\centering
\includegraphics[width=\textwidth]{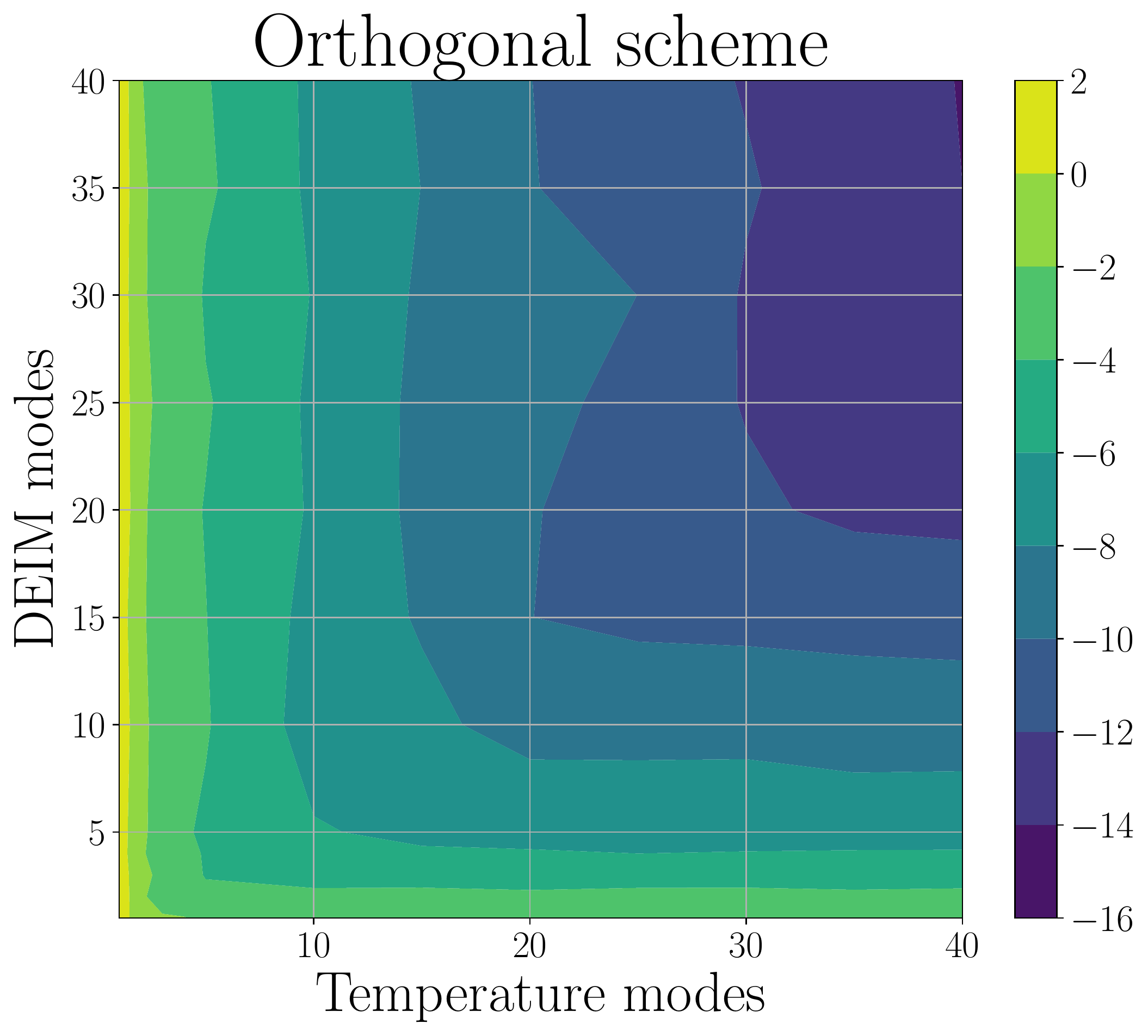}
\end{minipage} 
\end{minipage} 
\caption{The average of the $L^2$ relative error norm of the temperature field for the RBF mesh motion strategy. The logarithm of the $L^2$ relative error norm is plotted against the number of modes used for the temperature field and the number of D-EIM modes used to approximate the discretized differential operator $\bm A$ and the source term $\bm f$ (with $N_A = N_f$). The plots are reported with (left) and without (right) non-orthogonal correction.}
\label{fig:RBF_comp}
\end{figure}

\begin{figure}[H]
\centering
\begin{minipage}{\textwidth}
\centering
\begin{minipage}{0.49\textwidth}
\centering
\includegraphics[width=\textwidth]{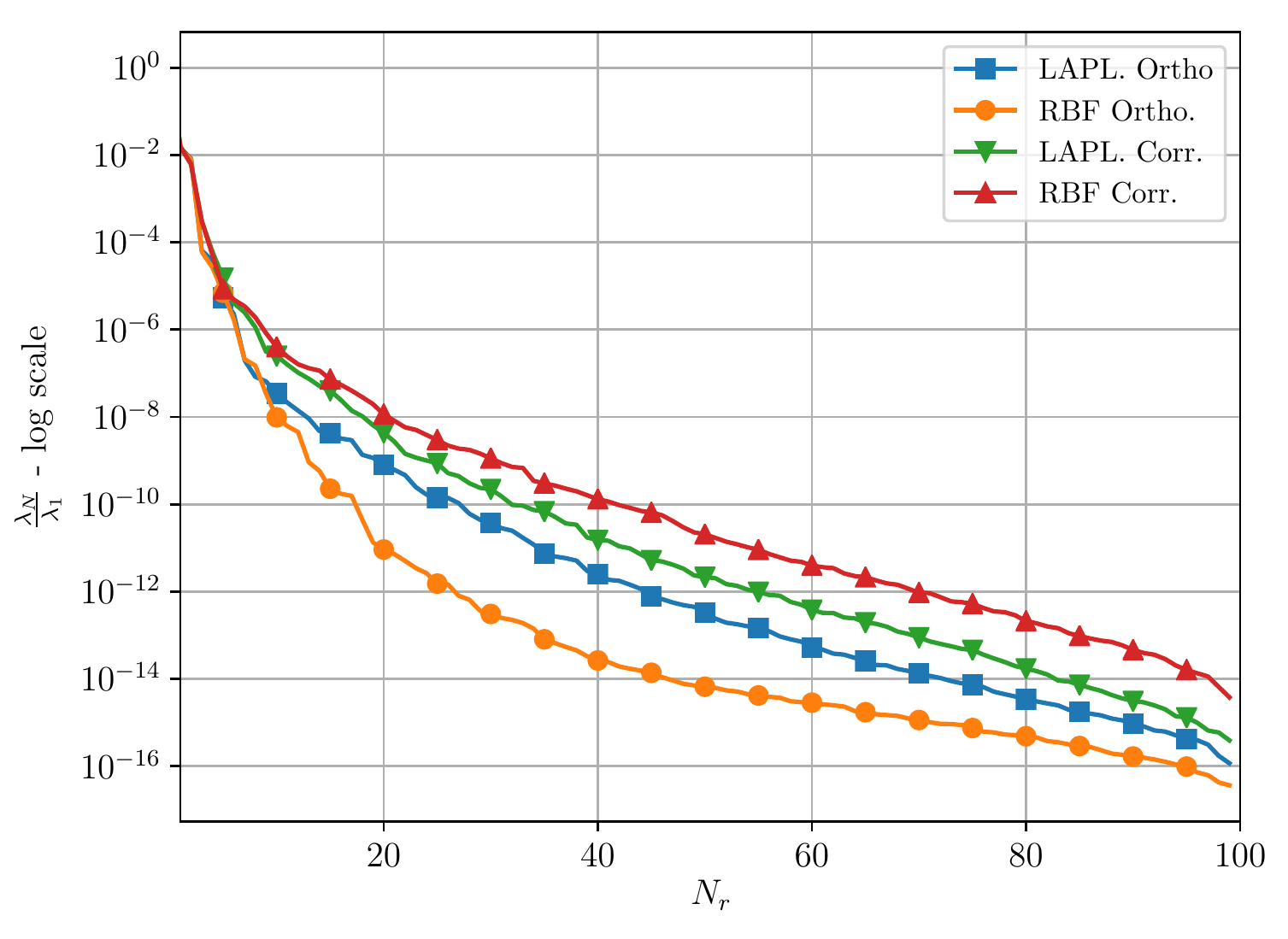}
\end{minipage} 
\end{minipage} 
\caption{Plot reporting the eigenvalue decay relative to the POD procedure used to compute the modes for the temperature field.}
\label{fig:T_eigenvalue}
\end{figure}

\begin{figure}[H]
\centering
\begin{minipage}{\textwidth}
\centering
\begin{minipage}{0.49\textwidth}
\centering
\includegraphics[width=\textwidth]{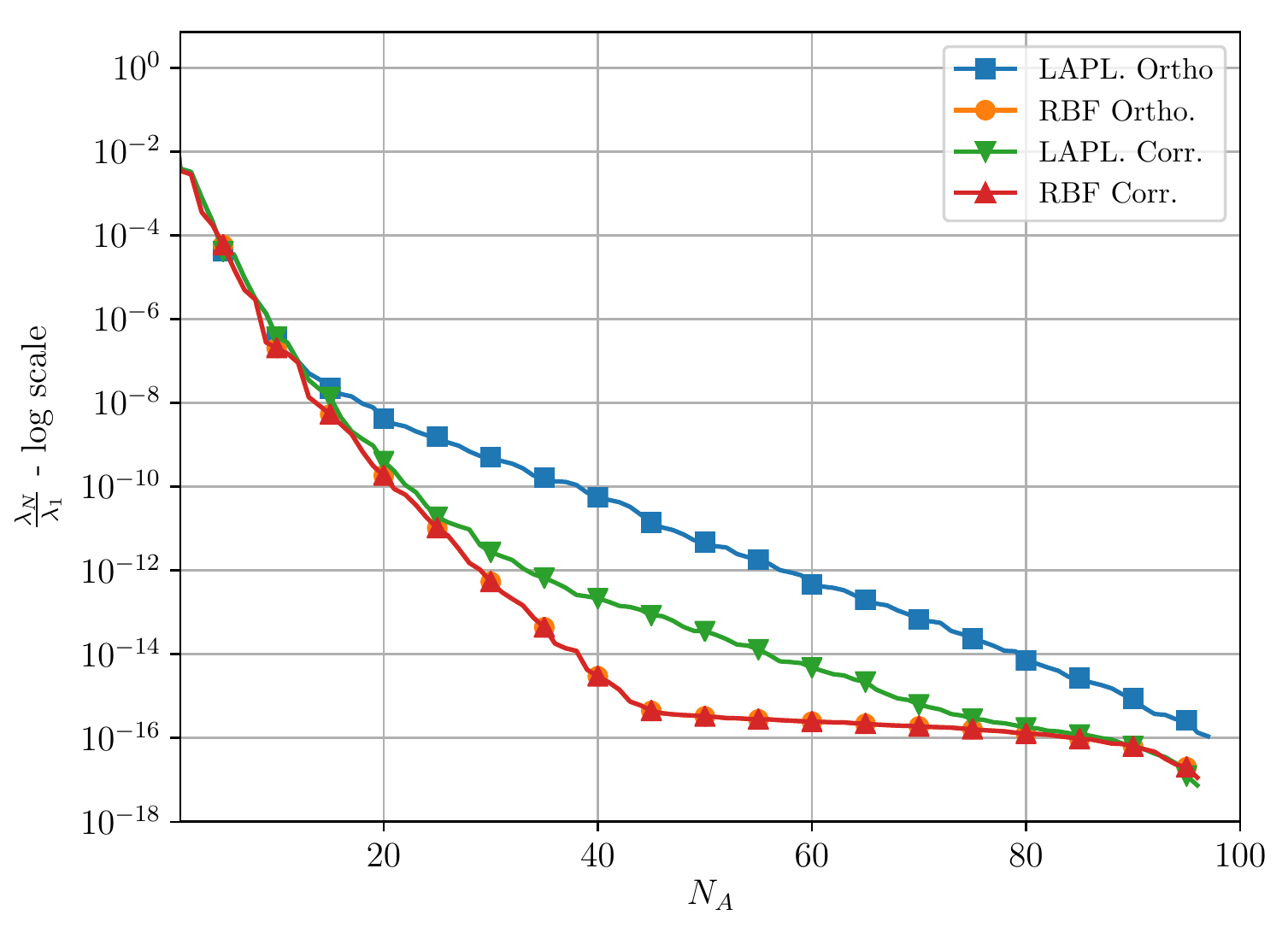}
\end{minipage} 
\begin{minipage}{0.49\textwidth}
\centering
\includegraphics[width=\textwidth]{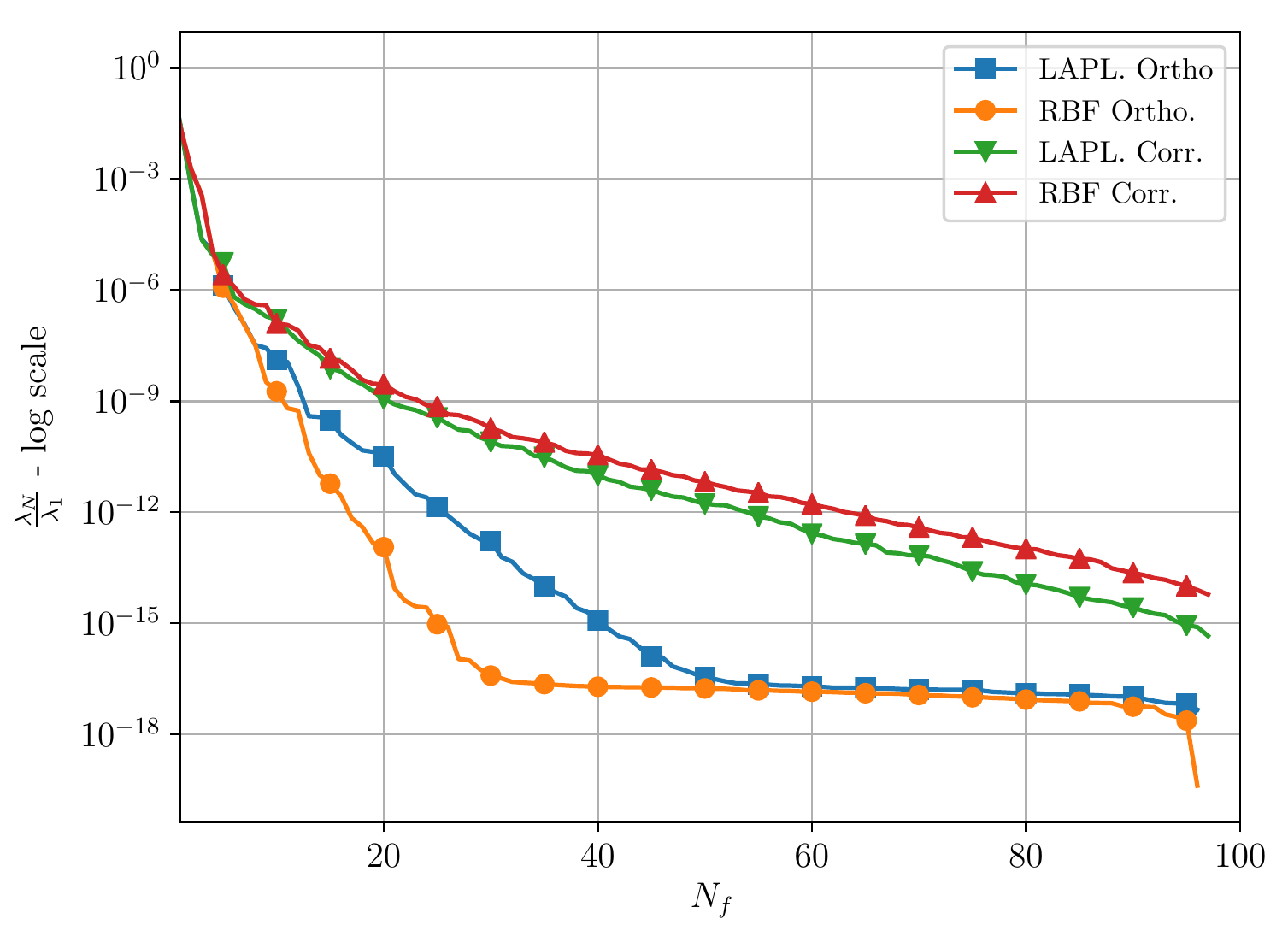}
\end{minipage} 
\end{minipage} 
\caption{Plot reporting the eigenvalue decay relative to the POD procedure used to compute the D-EIM modes for the discretized differential operator $\bm A$ (on the left) and the source term $\bm f$ (on the right).}
\label{fig:AB_eigenvalue}
\end{figure}

\begin{figure}[H]
\centering
\begin{minipage}{\textwidth}
\centering
\begin{minipage}{0.49\textwidth}
\centering
\includegraphics[width=\textwidth]{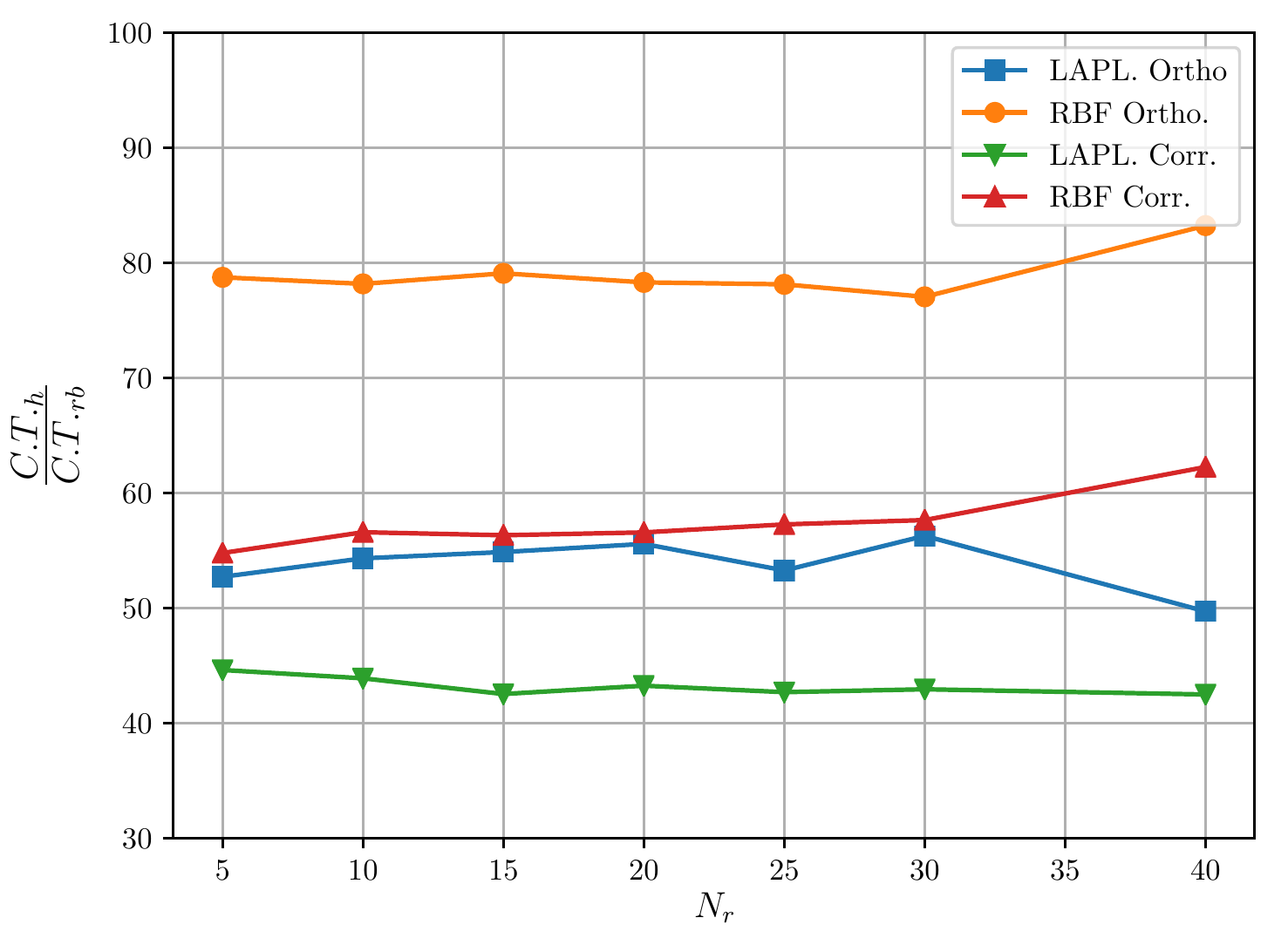}
\end{minipage} 
\begin{minipage}{0.49\textwidth}
\centering
\includegraphics[width=\textwidth]{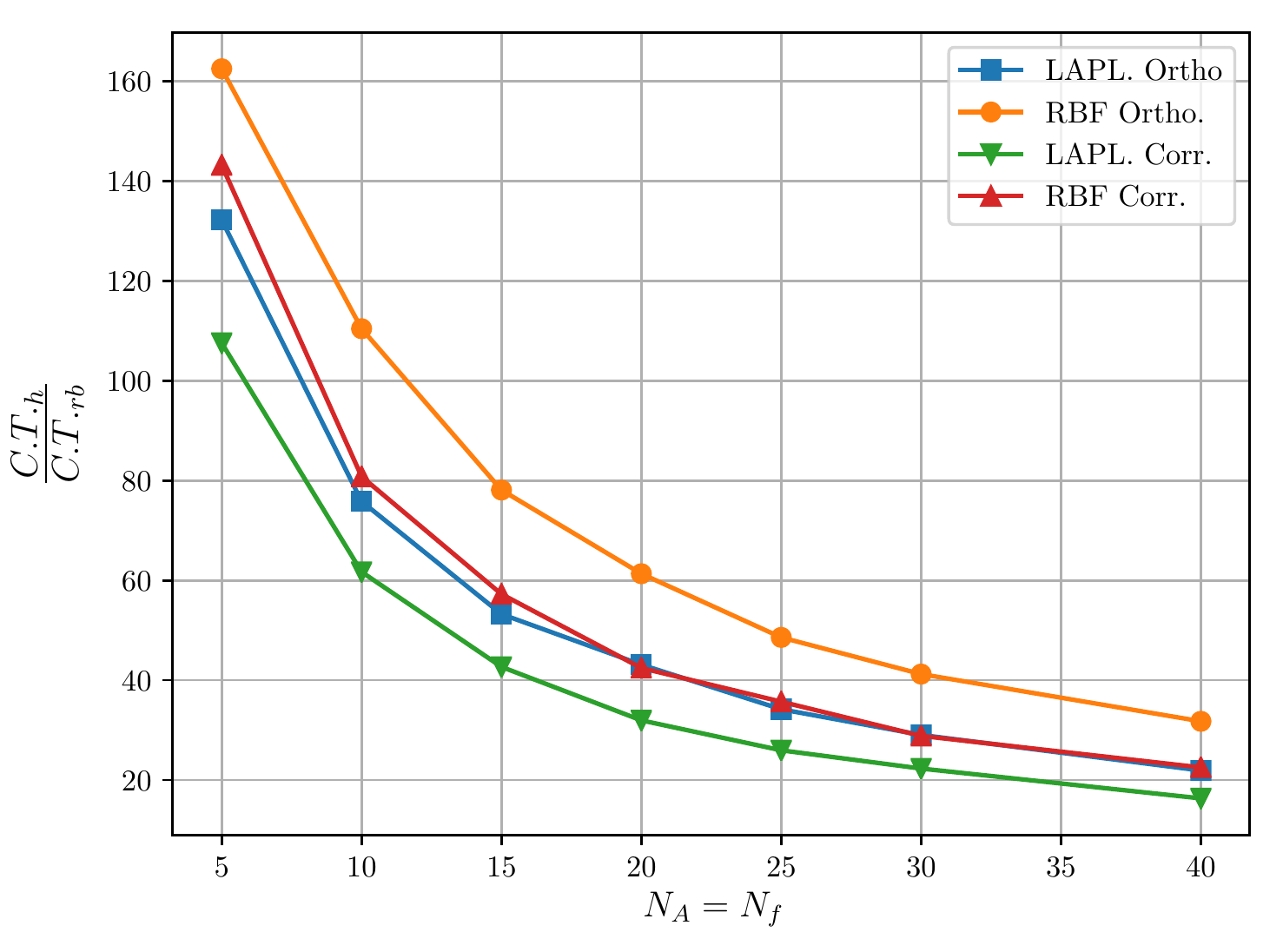}
\end{minipage} 
\end{minipage} 
\caption{Comparison of the computational time, in terms of speedup, for the different strategies. On the left we have a comparison keeping constant the number of D-EIM modes used to approximate the discretized differential operator $\bm A$ and the source term $\bm f$ ($N_{A} = N_f = 15$ and changing the number of temperature modes $N_r$. On the right we have a comparison keeping constant the number of temperature modes and changing the number of D-EIM modes used to approximate the discretized differential operator $\bm A$ and the source term $\bm f$.}
\label{fig:computational_time}
\end{figure}

Once we have determined the characteristic of the analyzed mesh motion strategies at the full order level, we examine the characteristics of a reduced order model generated using the different approaches. We will compare the results using two different schemes with and without orthogonal correction. The reduced order model is constructed using the discrete set $\mathcal{K}^1_{train}$ as defined in \autoref{subsec:mms} and the results are tested on a different testing sets $\mathcal{K}_{test} = \{\bm{\kappa}_{i_{test}}\}_{i=1}^{N_{test}} \in [-0.28,0.28] \times [-0.28,0.28]$ with $N_{test}=100$. As mentioned in \autoref{sec:FOM}, in order to improve the accuracy of the results, it is possible to employ, at the full order level, a non-orthogonal correction strategy. In what follows we will indicate with $\bm \theta_h$ the temperature field obtained using the full order model and with $\bm \theta_{rb}$ the temperature field obtained using the reduced order model. The accuracy of the ROM is measured using the following expression to identify the approximation error:
\begin{equation}\label{eq:appr_error}
\frac{1}{N_{test}} \sum_{i=1}^{N_{test}} \frac{|| \bm \theta_{h,i} - \bm \theta_{rb,i} ||_{\Omega}}{|| \bm \theta_{h,i}||_\Omega},
\end{equation}
where $N_{test}$ is the number of the parameter samples in the testing space and $||\cdot||_{\Omega}$ denotes the $L_2$ norm over the computational domain $\Omega$. In \autoref{fig:Lapl_comp} we report a comparison of the approximation error for the Laplacian mesh motion strategy, with and without orthogonal correction, while changing the dimension of the reduced basis space and the number of D-EIM modes used for the approximation of the discretized differential operator $\bm A$ and the source term $\bm f$. In this case we chose to use the same number of D-EIM modes $N_A = N_f$ for both $\bm A$ and $\bm f$. In \autoref{fig:RBF_comp} we report a similar plot but this time applying the RBF mesh motion procedure. In \autoref{fig:T_eigenvalue} and \autoref{fig:AB_eigenvalue} we report the eigenvalue decay relative to the POD procedure used to construct the modes for the temperature field and for the discretized differential operator $\bm A$ and source term $\bm f$. To measure the efficiency of the different model reduction strategies we measured also the computational speedup changing the number of D-EIM modes and the number of temperature modes. This comparison is depicted in \autoref{fig:computational_time}. 

From the numerical example we can draw several conclusions. We can note that the RBF mesh motion approach is generally performing slightly better respect to the Laplacian smoothing approach. Moreover the orthogonal non-corrected scheme, produces a ROM with better approximation properties respect to the corrected approach\footnote{The approximation error is measured using equation \autoref{eq:appr_error} with respect to the FOM and not to the analytical solution. A corrected scheme usually produces numerical results with a smaller numerical error and therefore a comparison computed using the analytical solution might produced a different trend in the results.}. 
%\autoref{fig:RBF_comp} shows a remarkable difference between the orthogonal approach and the corrected approach. In the first case, the approximation error reduces monotonically as the number of D-EIM modes increases, while in the second case, for both mesh motion strategies, using more than $10$ D-EIM modes we observe that the approximation error reaches a plateau and stops improving. We also observe that for the particular case, especially for the orthogonal scheme, it is more effective to increase the number of D-EIM modes rather than the number of POD modes used to approximate the temperature field. 
From \autoref{fig:computational_time} one observes two different evidences. The first one is that the RBF approach produces a ROM with the better computational speedups and the same is true for the orthogonal approach. For what concern the orthogonal approach this fact is due to the smaller computational stencil required by the methodology without correction. The second evidence is that an increase into the number of D-EIM modes produces a remarkable decrease into the computational speedup while the increase into the number of POD modes is negligible.

\begin{figure}[H]
\centering
\begin{minipage}{\textwidth}
\centering
\begin{minipage}{0.24\textwidth}
\centering
\footnotesize FOM Lapl. Ortho.\\
\includegraphics[width=\textwidth]{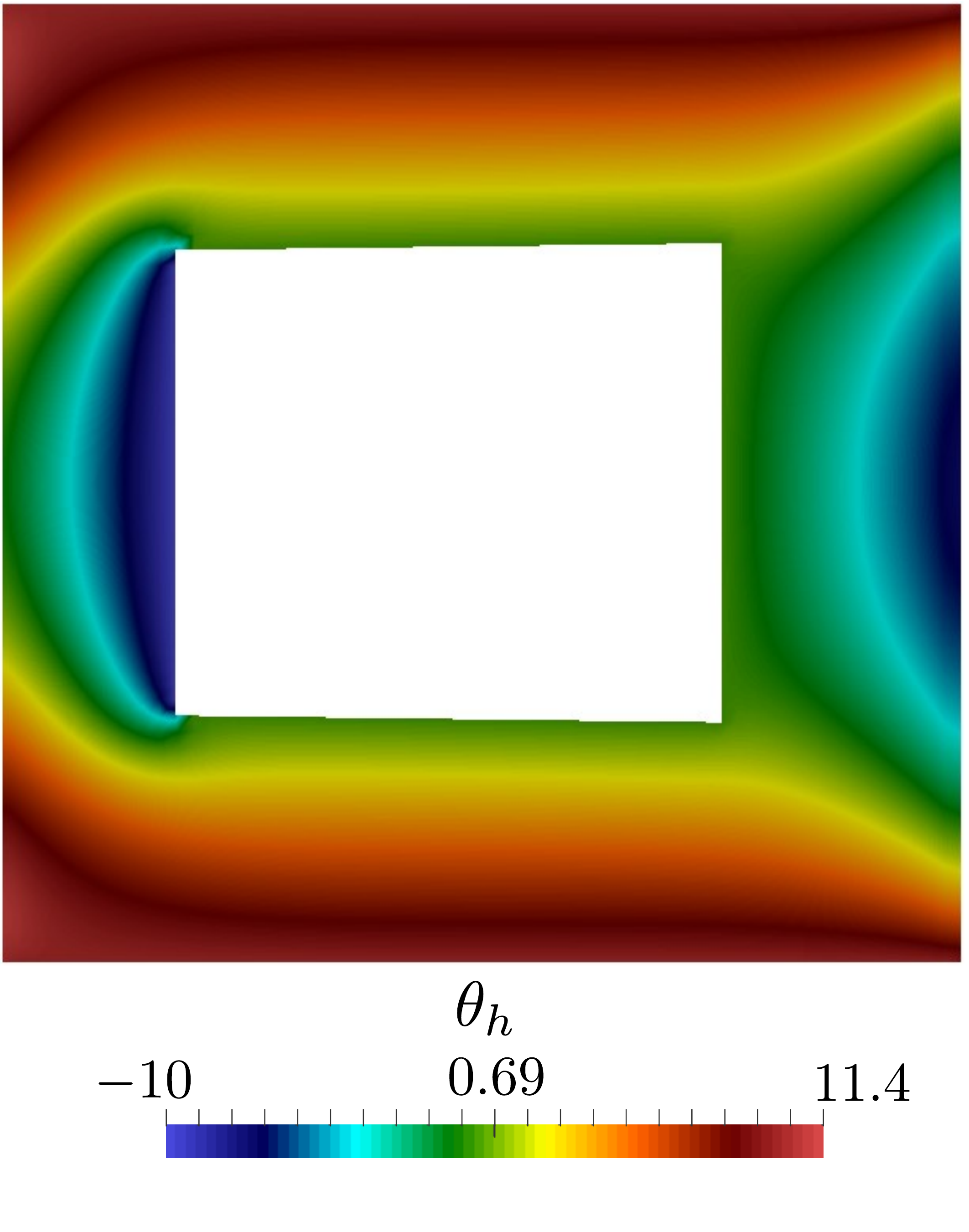}
\vspace{0.2cm}
\end{minipage} 
\begin{minipage}{0.24\textwidth}
\centering
\footnotesize ROM Lapl. Ortho. \\
\includegraphics[width=\textwidth]{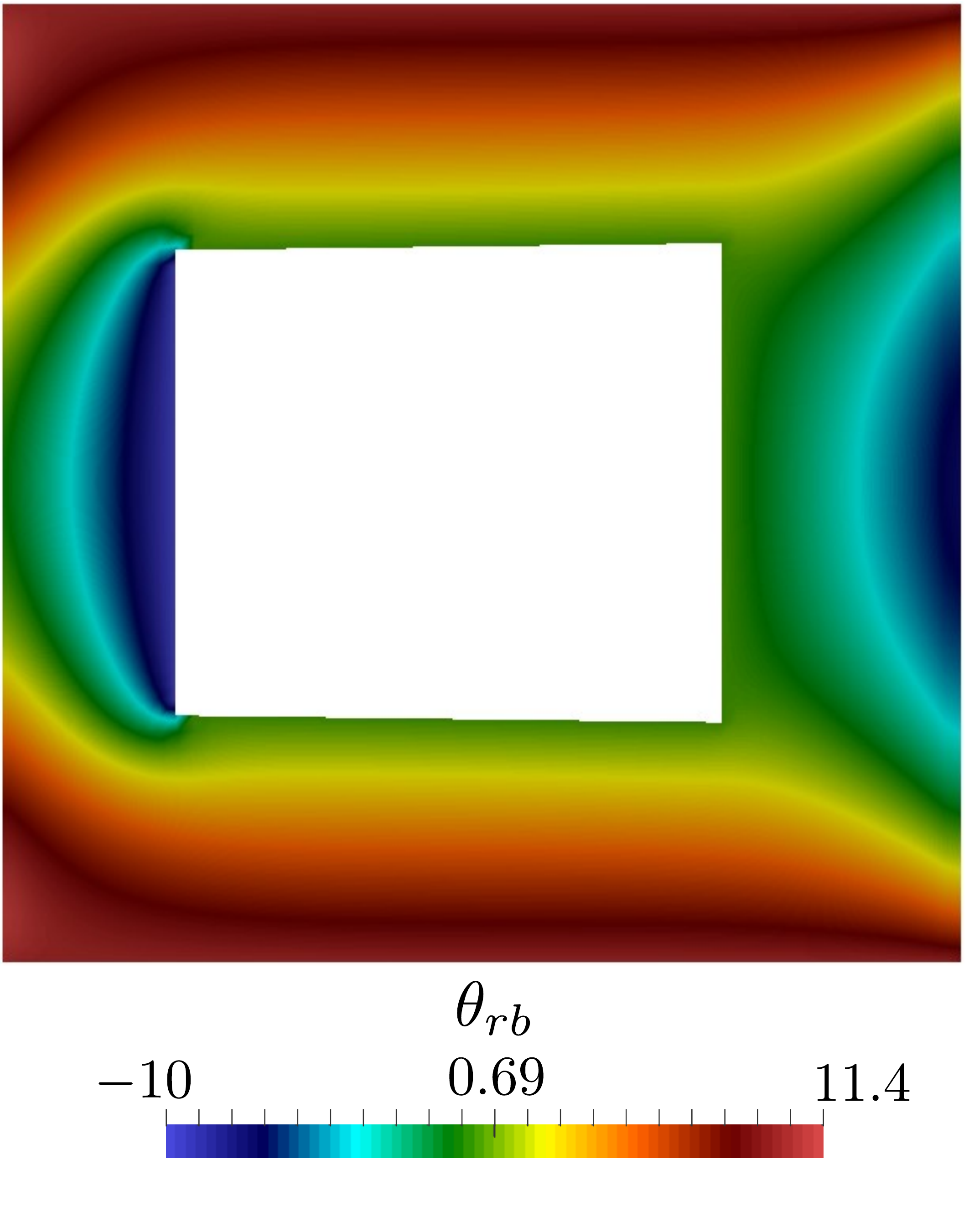}
\vspace{0.2cm}
\end{minipage} 
\begin{minipage}{0.24\textwidth}
\centering
\footnotesize FOM Lapl. Corr. \\
\includegraphics[width=\textwidth]{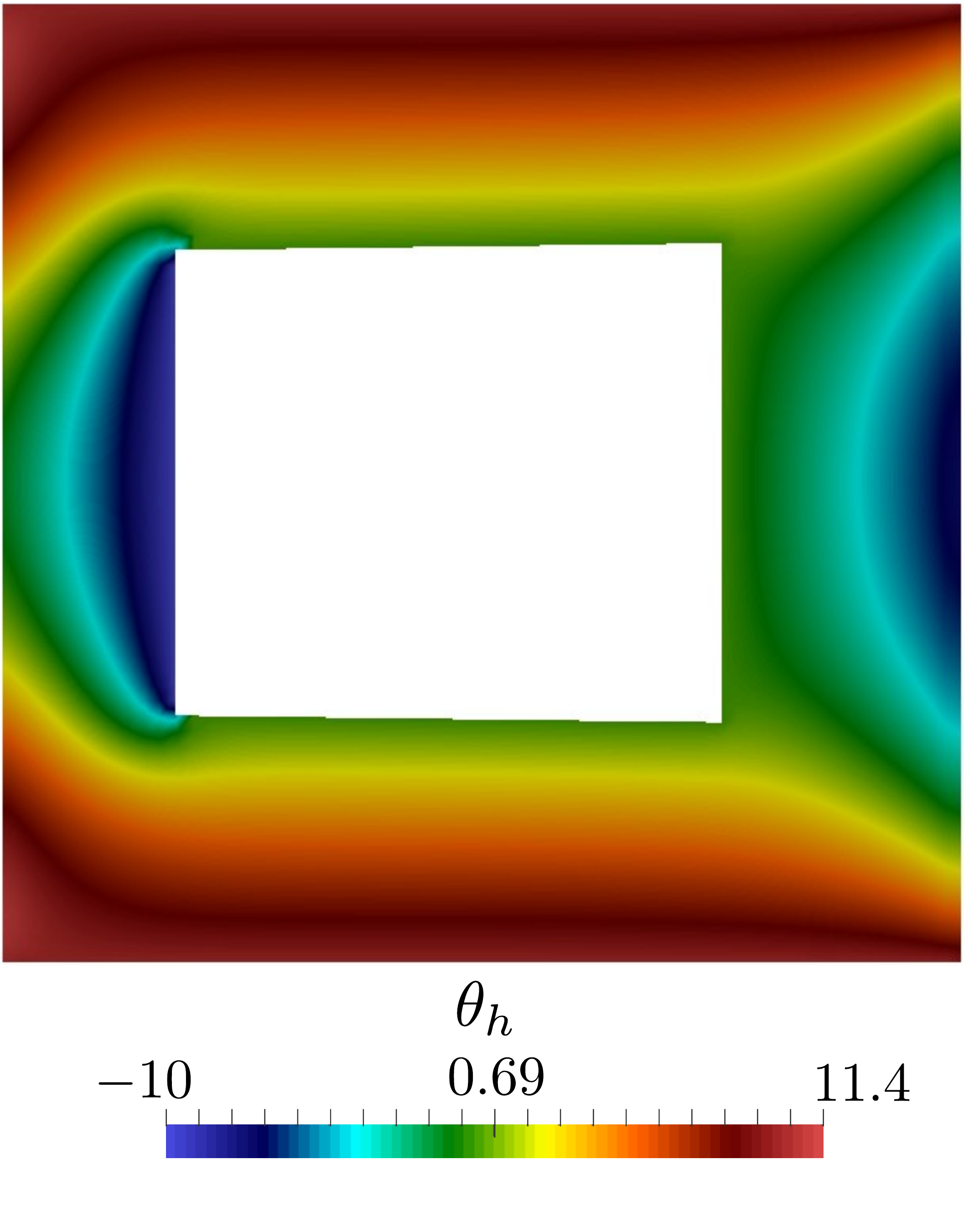}
\vspace{0.2cm}
\end{minipage} 
\begin{minipage}{0.24\textwidth}
\centering
\footnotesize ROM Lapl. Corr. \\
\includegraphics[width=\textwidth]{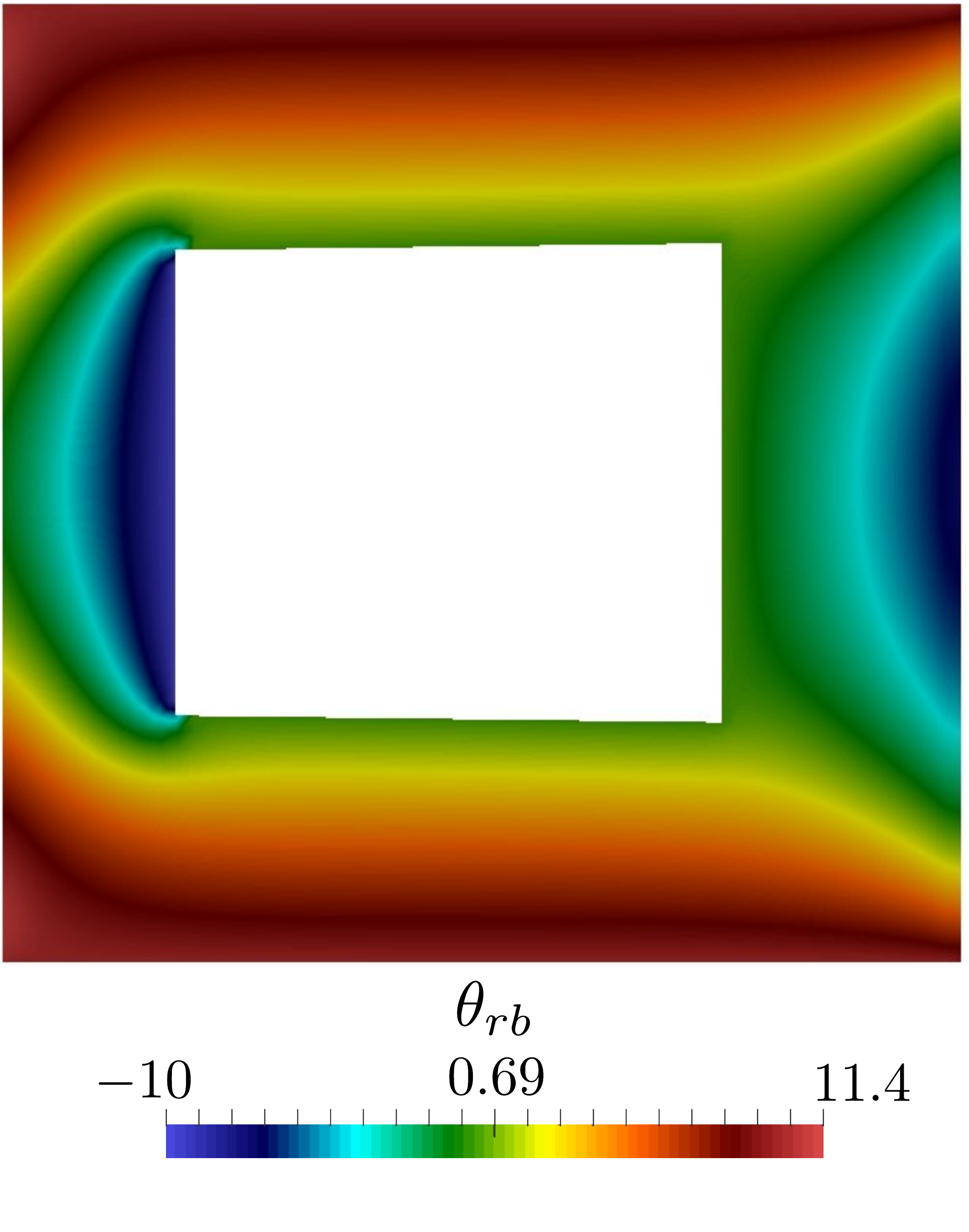}
\vspace{0.2cm}
\end{minipage} 
\begin{minipage}{0.24\textwidth}
\centering
\footnotesize FOM RBF. Ortho. \\
\includegraphics[width=\textwidth]{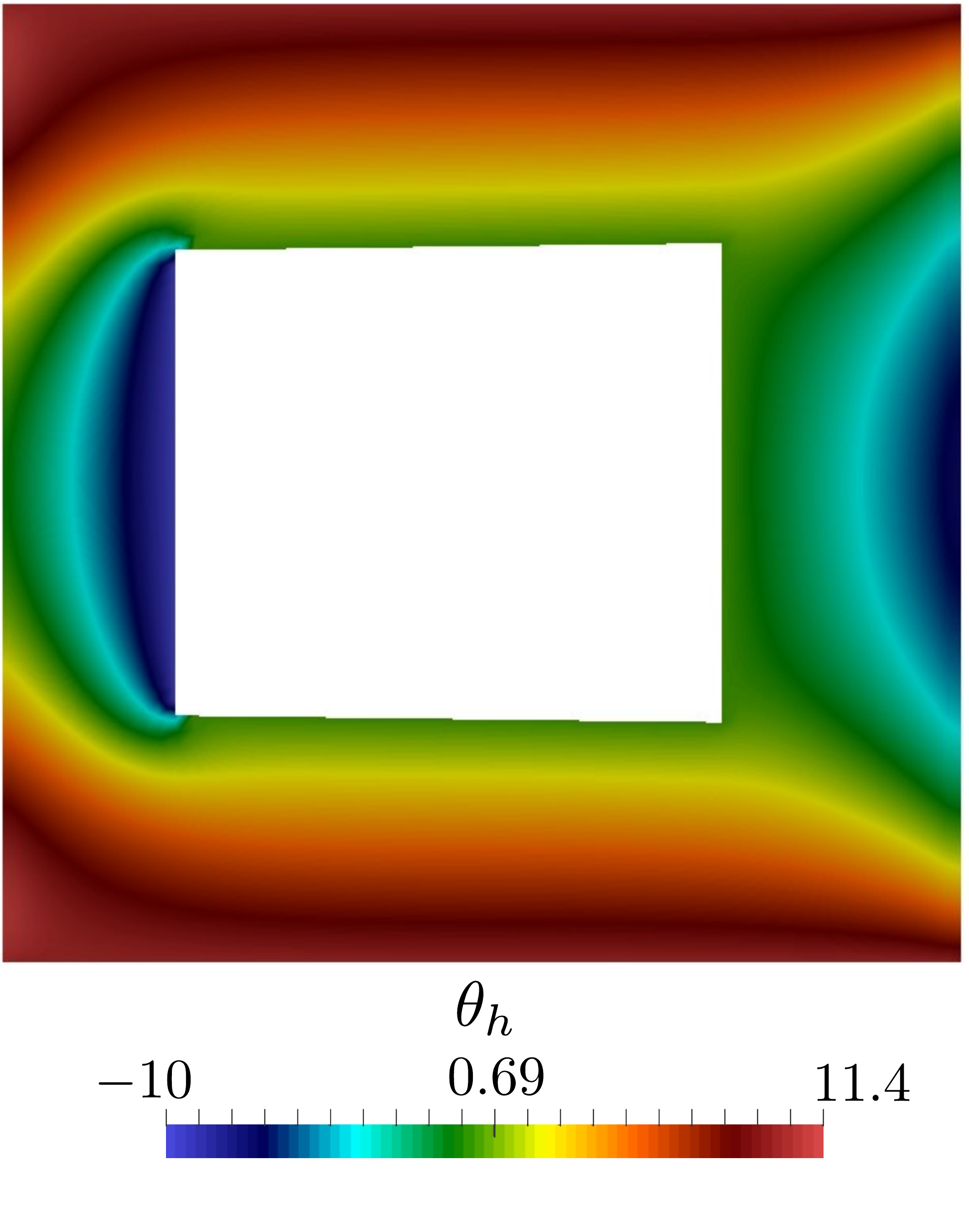}
\end{minipage} 
\begin{minipage}{0.24\textwidth}
\centering
\footnotesize ROM RBF Ortho. \\
\includegraphics[width=\textwidth]{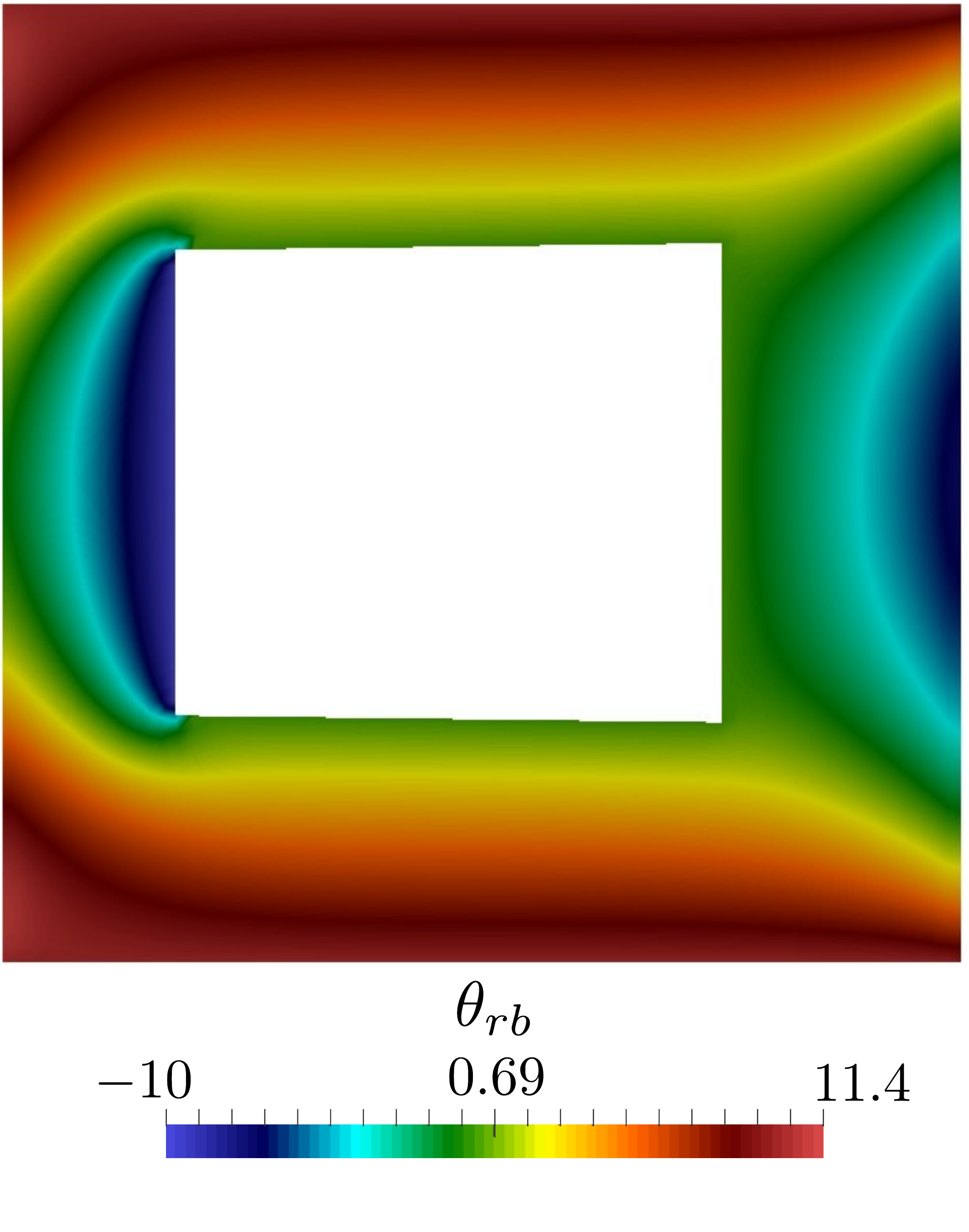}
\end{minipage} 
\begin{minipage}{0.24\textwidth}
\centering
\footnotesize FOM RBF Corr. \\
\includegraphics[width=\textwidth]{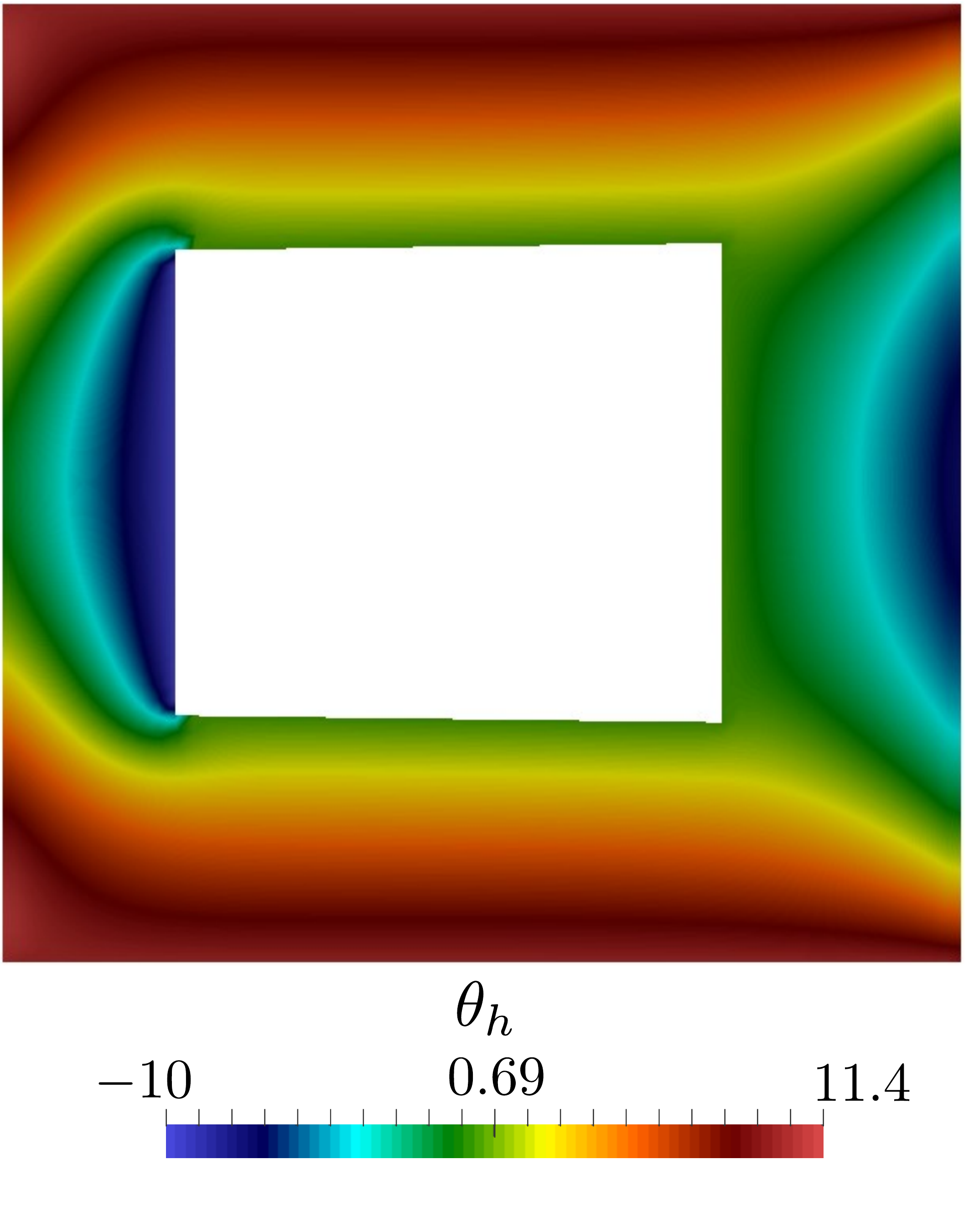}
\end{minipage} 
\begin{minipage}{0.24\textwidth}
\centering
\footnotesize ROM RBF Corr. \\
\includegraphics[width=\textwidth]{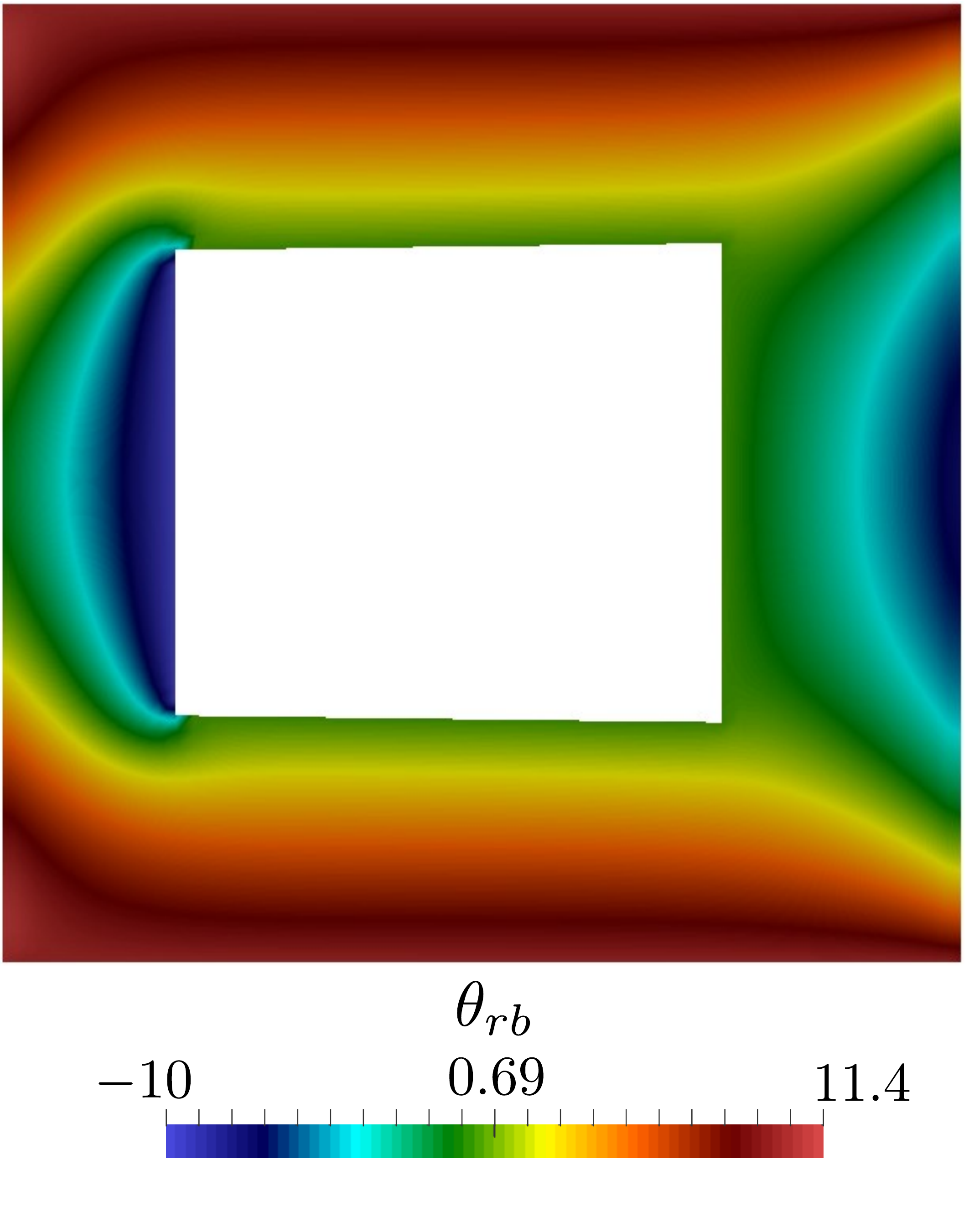}
\end{minipage} 
\end{minipage} 
\caption{{In the plot we report the comparison between full order solution and reduced order solution using the Laplacian smoothing mesh motion approach (first row) and the RBF mesh motion approach (second row). Results are reported without (columns $1-2$) and with (columns $3-4$) orthogonal correction. All the results are for one selected sample point in the parameter space that was giving the worst results in terms of accuracy $\bm{\kappa*} = (-0.213777, 0.0335904)$. } }
\label{fig:DEIM_comp_fig}
\end{figure}

\section{A parametrized incompressible Navier-Stokes problem}\label{sec:num_ex2}
In the previous section we focused the attention on relatively simple heat transfer problems. The choice of such an easy case was justified by the fact that we wanted to reduce as much as impossible additional sources of problems that are generally linked with hyperbolic equations and advection dominated cases in order to provide a better understanding of the influence of the mesh motion strategy in view of model reduction purposes. In this section we focus instead the attention on systems generated by the steady incompressible Navier--Stokes equations and on the difficulties that we have to face in order to extend the methodology introduced in the previous sections to a more complex scenario. 
\subsection{Mathematical Formulation}
The focus is here on a geometrically parametrized problem which is approximated by the steady Navier--Stokes equations. 
Considering an Eulerian frame on a geometrically parametrized domain $\Omega(\bm \mu) \subset \mathbb{R}^d$ with $d=2,3$ the problem consists in finding the vectorial velocity field $\bm{u}(\bm \mu):\Omega(\bm \mu) \to \mathbb{R}^d$ and the scalar pressure field $p(\bm \mu):\Omega(\bm \mu) \to \mathbb{R}$ such that:
\begin{equation}\label{eq:navstokes}
\begin{cases}
\bm{\nabla} \cdot (\bm{u} \otimes \bm{u})- \bm{\nabla} \cdot 2 \nu \bm{\nabla^s} \bm{u}=-\bm{\nabla}p &\mbox{ in } \Omega(\bm \mu),\\
\bm{\nabla} \cdot \bm{u}=\bm{0} &\mbox{ in } \Omega(\bm \mu),\\
\bm{u} (x) = \bm{f}(\bm{x}) &\mbox{ on } \Gamma_{In}(\bm \mu),\\
\bm{u} (x) = \bm{0} &\mbox{ on } \Gamma_{0}(\bm \mu) ,\\ 
(\nu\nabla \bm{u} - p\bm{I})\bm{n} = \bm{0} &\mbox{ on } \Gamma_{Out}(\bm \mu),\\ 
\end{cases}
\end{equation}
where $\Gamma(\bm \mu) = \Gamma_{In}(\bm \mu) \cup \Gamma_{0}(\bm \mu) \cup \Gamma_{Out}(\bm \mu)$ is the boundary of $\Omega(\bm \mu)$ and, is composed by three different parts $\Gamma_{In}(\bm \mu)$, $\Gamma_{Out}(\bm \mu)$ and $\Gamma_0(\bm \mu)$ that indicates, respectively, inlet boundary, outlet boundary and physical walls. The function $\bm{f}(\bm{x})$ represents the boundary conditions. The parameter dependency is given by a parameter vector $\bm \mu\in \mathcal{P}$ with $\mathcal{P}$ denoting the parameter space that describes the deformation of the physical domain. It is moreover assumed that the kinematic viscosity $\nu$ is constant in the spacial domain. For sake of brevity, the parameter dependency of the domain $\Omega$ will be omitted in the next formulations.
\subsection{The discrete formulation}
The above equations are discretized at the full order level using a cell-centered finite volume approach and a SIMPLE algorithm to resolve the pressure velocity coupling. The equations are then written in integral form over each control volume $\Omega_i$ and the divergence terms are transformed into surface integral thanks to the exploitation of the Gauss's theorem and computed numerically as sum of fluxes over the surfaces of each cell $\Omega_i$. In what follows we briefly report the discretization of each terms that appears inside the momentum and continuity equation.
The non-linear convective term is discretized as:
\begin{equation}\label{eq:convect_discr}
\int_{\Omega_i} \bm \nabla \cdot (\bm{u}\otimes\bm{u}) \dif \Omega = \int_{\partial \Omega_i} \bm{n} \cdot (\bm{u}\otimes\bm{u}) \dif \Gamma = \sum_f \bm{S_f} \cdot \bm{u}_f \otimes \bm{u}_f = \sum_f \bm{S_f} \cdot \bm{u}_f \otimes \bm{u}_f = \sum_f F_f \bm{u}_f ,
\end{equation}
where $\bm{u}_f$ indicates the velocity at the centre of the faces. Since the equation is solved using an iterative approach the non-linear term is linearized with the substitution of the term $\bm{S}_f \cdot \bm{u}_f = F_f$, that represents the mass flux over each face, with a previously calculated value of the mass flux that satisfies the continuity equation, for more details about this issue we refer to \cite{Jasak1996}. The diffusive term is discretized as:
\begin{equation}
\int_{\Omega_i} \bm{\nabla} \cdot 2 \nu \bm{\nabla^s} \bm{u} \dif \Omega \int_{\partial \Omega_e} \bm{n} \cdot 2 \nu \bm{\nabla^s}\bm{u} \dif \Gamma = \int_{\partial \Omega_e} \bm{n} \cdot \nu \bm{\nabla}\bm{u} \dif \Gamma = \nu \sum_f \bm{S_f} \cdot (\bm{\nabla u})_f,
\end{equation}
where the first equality follows from the incompressibility constraint and the term  $(\bm{\nabla u})_f$ indicates the gradient of the velocity field at the centre of each face. This is calculated, starting from the values at the centre of the neighbouring cells, using a finite difference scheme that includes a correction in the case of non-orthogonal meshes (see \autoref{eq:nonortho}). For more details on available choices concerning the correction term we refer to \cite{Jasak1996}.

The term originated from the gradient of pressure, is discretized as:
\begin{equation}
\int_{\partial \Omega_e} \bm{n} p  \dif \Gamma = \sum_f \bm{S_f} p_f,
\end{equation}
while the term originated from the divergence of velocity is discretized as:
\begin{equation}\label{eq:NS_con}
\int_{\partial \Omega} \bm{n} \cdot \bm{u} \dif \Gamma = \sum_{f=1}^{N_f} \bm{S_f} \cdot \bm{u_f} = \sum_{f=1}^{N_f} F_f. 
\end{equation}

\subsection{The SIMPLE algorithm}\label{subsec:simple}
Now that the discretization of each term has been introduced we describe the solution strategy used to resolve the coupling between the momentum and continuity equation. The coupling between velocity and pressure is resolved using a SIMPLE strategy \cite{Patankar1972}. The equations can be rewritten in matrix form as:
\begin{equation}
\begin{bmatrix}
[\bm A_u] & [\nabla(\cdot)] \\
[\nabla \cdot (\cdot)]  & [0]
\end{bmatrix} 
\begin{bmatrix}
\bm u \\
p
\end{bmatrix} = 
\begin{bmatrix}
\bm 0 \\
0
\end{bmatrix}
\end{equation}
where $\bm A_u$ is the coefficient matrix coming from the momentum equation, $A_u \bm u = \bm \nabla \cdot (\bm u \otimes \bm u) - \bm \nabla \cdot (2 \nu \bm \nabla^s \bm u)$. The above system matrix has a saddle point structure which is usually not easy to solve using a coupled approach. For this reason we rely on a segregated approach where the momentum equation is solved with a tentative pressure $\tilde p$ and later corrected exploiting the divergence free constraint. 
{\RA The first equation could then, in principle, be solved using a direct solver; however, for efficiency reasons, it is decomposed into a diagonal, upper triangular and lower triangular part $\bm A_u = \bm D_u + \bm L_u + \bm U_u$. } The momentum equation is then rewritten considering this decomposition and using the previous converged value of the velocity to assemble the operator $\bm H(\bm u) = \bm (L_u + \bm U_u) \bm u_{k-1}$. Since now the system matrix is diagonal it is easy to invert and a tentative value of the velocity can be computed as:
\begin{equation}\label{eq:momentum_equation}
\bm u = D_u^{-1}(\bm H (\bm u_{}) - \bm \nabla \tilde p),
\end{equation}
Taking the divergence of the tentative velocity field and exploiting the divergence free constraint one can derive a Poisson equation for pressure of the form:
\begin{equation}\label{eq:pressure_equation}
\bm \nabla \cdot (\bm D_u^{-1} \bm \nabla p) = \bm \nabla \cdot D_u^{-1}(\bm H (\bm u_{k-1})).
\end{equation}
Equations \ref{eq:momentum_equation} and \ref{eq:pressure_equation} are then  solved using an iterative SIMPLE approach \cite{Patankar1972} with an under-relaxation procedure to achieve convergence. The SIMPLE algorithm follows a segregated solution strategy where the momentum equation is firstly solved with a tentative pressure and an intermediate velocity field $\bm u^*$ that is generally not divergence-free is computed. This intermediate velocity is plugged inside the continuity equation to construct a Poisson equation for pressure that is used to obtain a new pressure $p$ that, if inserted inside the momentum equation, delivers a divergence free velocity field. This two steps are repeated until convergence.

\subsection{The Reduced Order Model}
Contrary to what was done in \cite{Stabile2017,stabile_stabilized} the reduced order model is constructed in such a way to be completely consistent with the procedure used at the full order level. In the previously mentioned works, the SIMPLE procedure used at the full order level was transformed into a saddle point problem at the reduced order level. This operation was possible because geometrical parametrization was not considered and therefore the reduced operators could be computed explicitly after the offline stage using the reduced basis functions. 
In this case, since the interest is into geometrical parametrization, and for a new value of the parameters a new mesh is constructed, the reduced order model is based on a reduced version of the SIMPLE algorithm presented in the previous subsection. During the online stage, the reduced operators are constructed on the deformed mesh using the methodology introduced in \autoref{sec:rom}.
For the construction of the reduced basis spaces we used a POD strategy, performed on an average deformed configuration, on the snapshots matrices of the velocity and pressure fields in order to obtain two separate set of reduced basis functions:
\begin{equation}\label{eq:mom_simple}
\bm L_u = [\bm \varphi^u_1, \dots, \bm \varphi^u_{N_r^u}], \quad \bm L_p = [\bm \varphi^p_1, \dots, \bm \varphi^p_{N_r^p}].
\end{equation}
We assume then that the velocity and the pressure fields can be approximated by:
\begin{equation}\label{eq:pres_simple}
\bm u \approx \bm u_r = \sum_{i=1}^{N_r^u} a_i^u \bm \varphi^u_{i}, \quad p \approx p_r = \sum_{i=1}^{N_r^p} a_i^u \varphi^p_{i}.
\end{equation}
As shown in the previous section the reduced basis spaces are generated using a modified version of the mass matrix which is defined as the ensemble average of all the mass matrices obtained during the training stage. In the case of a geometrically parametrized problem, \autoref{eq:mom_simple} and \autoref{eq:pres_simple} can be rewritten as:
\begin{equation}
\bm A_u (\bm u_{k-1}, \bm \mu) \bm u = \bm b_u (\bm \mu, p), \quad \bm A_p(\mu) p = \bm b_p (\bm \mu, \bm u).
\end{equation}
The system matrices and the source terms of both the momentum and pressure equation are therefore parameter dependent. In a reduced setting the SIMPLE algorithm can be rewritten following the procedure of \autoref{alg:simpleRed}.
\begin{algorithm}[t]

\caption{The Reduced SIMPLE algorithm}
\label{alg:simpleRed}
\hspace*{\algorithmicindent} \textbf{Input:} Tentative value of the velocity and pressure coefficient vectors $\bm a^u_0$ and $\bm a^p_0$, $k=1$, tolerance $tol$, $res = tol+1$.
\begin{algorithmic}[1]
\STATE Reconstruct the first attempt full velocity  and pressure fields: $\bm u_0 = {\bm L}_u \bm a^u_0 $, $\bm{p}_0 = \bm{L}_p \bm a^p_0$,  
\WHILE{$res>tol$} 
\STATE assemble reduced momentum equation $\bm A^r_u = \bm{L}_u^T \bm{A}_u(\bm u_{k-1}, \bm \mu) \bm{L}_u$; $\bm b^r_u = \bm{L}_u^T \bm b_u (\bm \mu, p_{k-1})$;  
\STATE compute reduced velocity residual $\to r_u = |\bm A_u^r \bm a^u_{k-1} - \bm b^r_u | $;
\STATE solve $\bm A^r_u \bm a^u_k = \bm b^r_u$ and reconstruct $\bm{u}_{k} = \bm{L}_u \bm a_k^u $;
\STATE assemble reduced pressure equation $\bm A^r_p = \bm{L}_p^T \bm{A}_p(\bm \mu) \bm{L}_p$; $\bm b^r_p = \bm{L}_p^T \bm b_p (\bm \mu, \bm u_{k})$; 
\STATE compute reduced pressure residual $\to r_p = |\bm A_p^r \bm a^p_{k-1} - \bm b^r_p | $;
\STATE solve $\bm A^r_p \bm a^p_k = \bm b^r_p$ and reconstruct $\bm{p}_{k} = \bm{L}_p \bm a_k^p $;
\STATE $res = max(r_u,r_p)$, $k=k+1$;
\ENDWHILE
\end{algorithmic}
\end{algorithm}
For sake of clearness in the algorithm the procedure is presented without additional hyper-reduction but the offline and online part can be properly decoupled relying on the D-EIM procedure shown for the heat transfer case. Therefore, during the online resolution only pointwise evaluation of the discretized differential operators $\bm A_u$, $\bm A_p$ and source terms $\bm b_u$, $\bm b_p$ are required. The same concept is valid for the reconstruction part where it is required to reconstruct the solution $\bm u$ and $p$ only in some points of the domain identified by the magic points of the D-EIM algorithm. The spirit of the reduced algorithm introduced here is similar to the one introduced in \cite{Isoz2019} with the main difference that in the mentioned reference the reduced order model is constructed only for the pressure equation and the velocity is reconstructed a posteriori from the pressure field. In the present case we built instead a reduced order model which considers both the velocity and pressure fields at the same time and that completely mimic the full order SIMPLE algorithm. 

\subsection{A first numerical result on a geometrically parametrized incompressible flow problem}
\begin{figure}[H]
\centering
\begin{minipage}{\textwidth}
\centering
\begin{minipage}{0.49\textwidth}
\centering
\includegraphics[width=\textwidth]{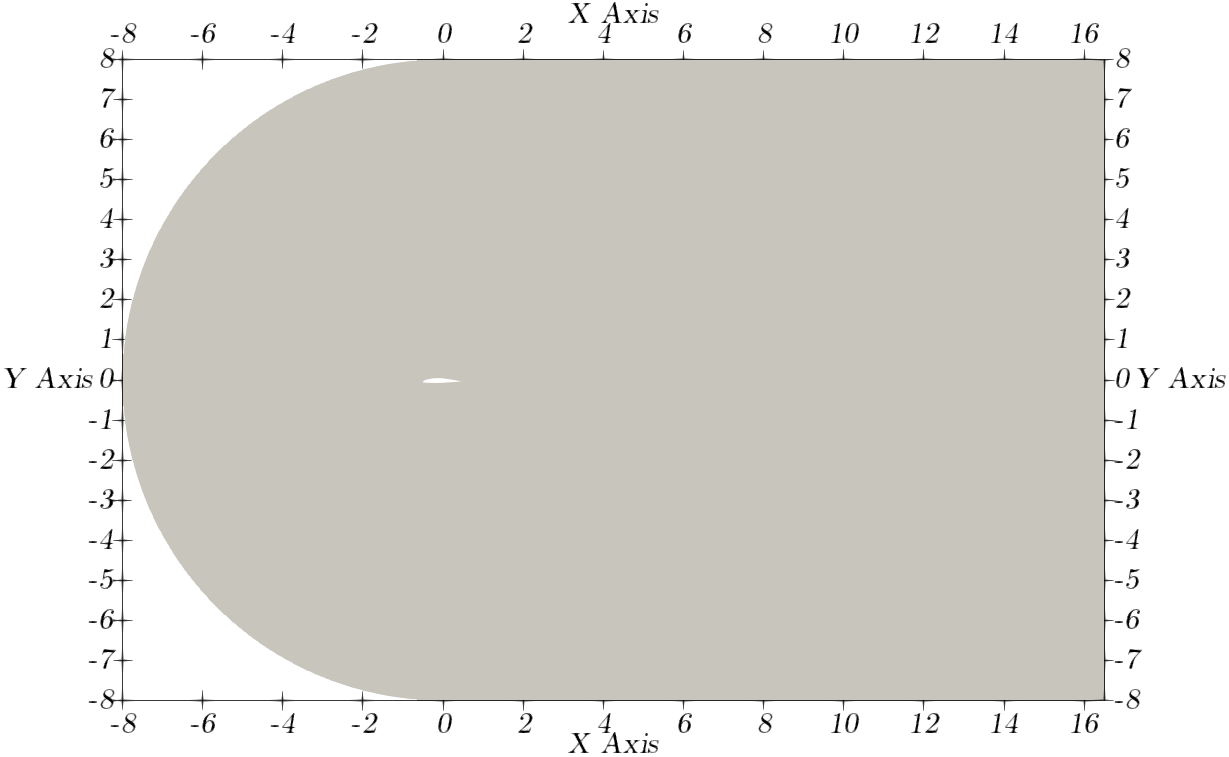}
\end{minipage} 
\begin{minipage}{0.49\textwidth}
\centering
\includegraphics[width=\textwidth]{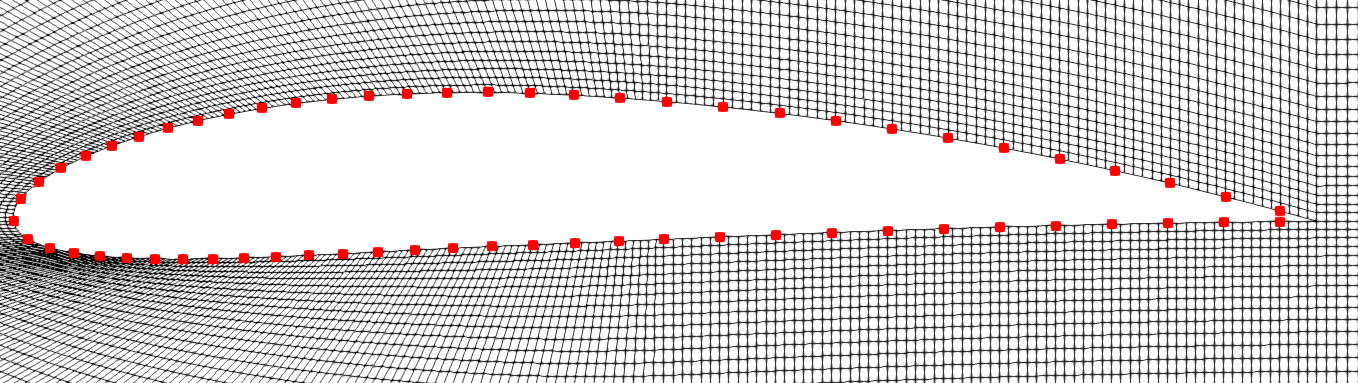}\\
\vspace{1cm}
\begin{tabular}{l | c | c | c}
  & Domain-X & Domain-Y & wing-chord\\
  \hline
  Dimension & $24.5$ & $16$ & $1$ \\
\end{tabular}
\end{minipage} 
\end{minipage} 
\caption{On the left we report a sketch of the mesh for the geometrical parametrized wing problem and on the right we can observe the location of the $64$ control points used to solve the RBF mesh motion problem together with the main dimensions of the domain.}
\label{fig:control_wing}
\end{figure}

\begin{figure}
\centering 
\includegraphics[width=0.5 \textwidth]{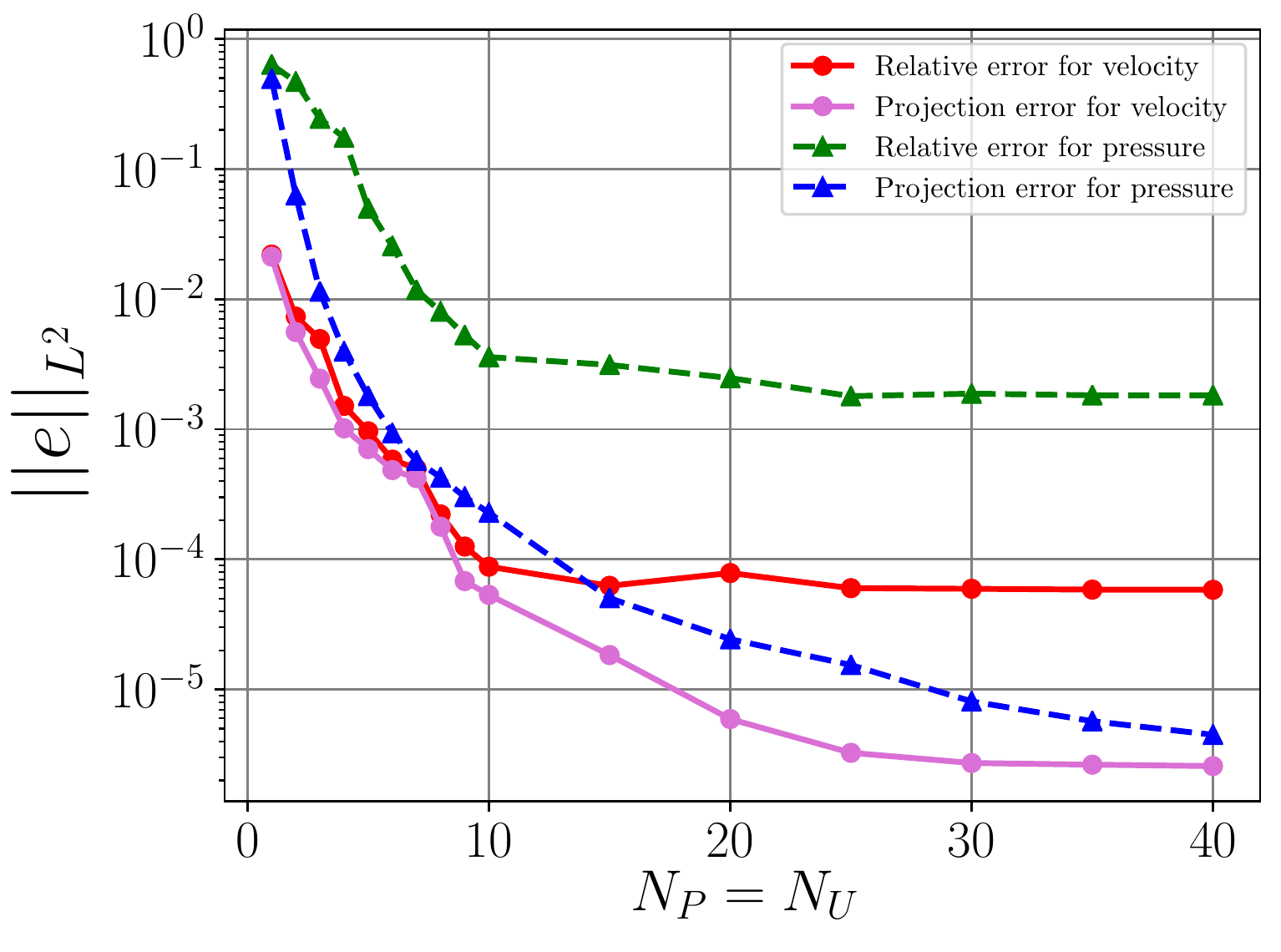}
\caption{\RA The average of the $L^2$ relative error norm of pressure and velocity fields: $\frac{P_{FOM} - P_{ROM}}{P_{FOM}}$ and $\frac{(U_{FOM} - U_{\infty})-(U_{ROM} - U_{\infty})}{U_{FOM} - U_{\infty}}$ where $U_{\infty}$ is the free stream velocity. The $L^2$ norm of the relative error for both pressure and velocity is plotted against the number of modes used for the reconstruction of the online solution in logarithmic $y$ scale. In the same image the reduced errors are compared with the projection errors obtained by the use of the same basis functions and by using the same norm.}
\label{fig:NSerrors}
\end{figure}

\begin{figure}

\begin{minipage}{\textwidth}
\centering
\begin{minipage}{0.24\textwidth}
\includegraphics[width=\textwidth]{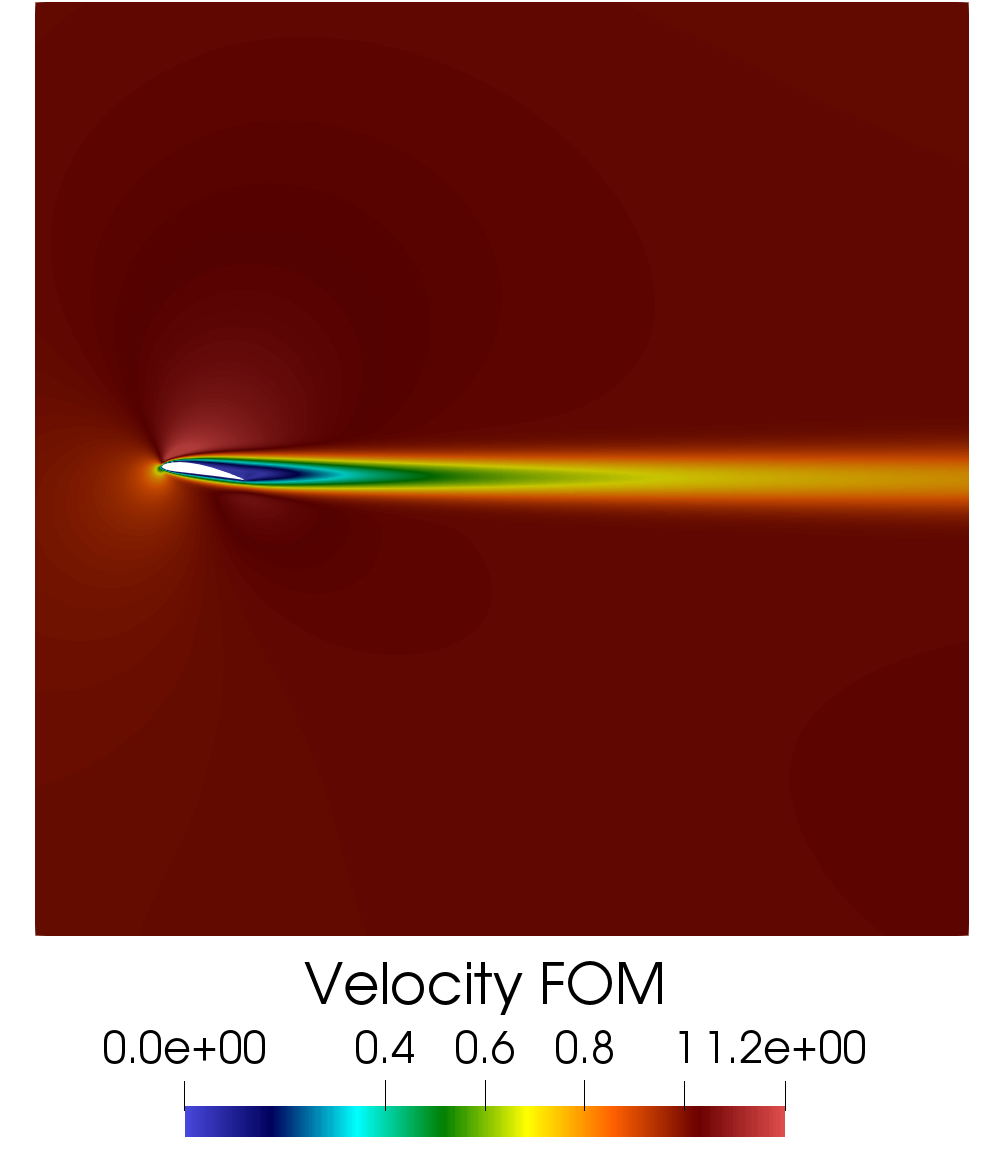}
\end{minipage}
\begin{minipage}{0.24\textwidth}
\includegraphics[width=\textwidth]{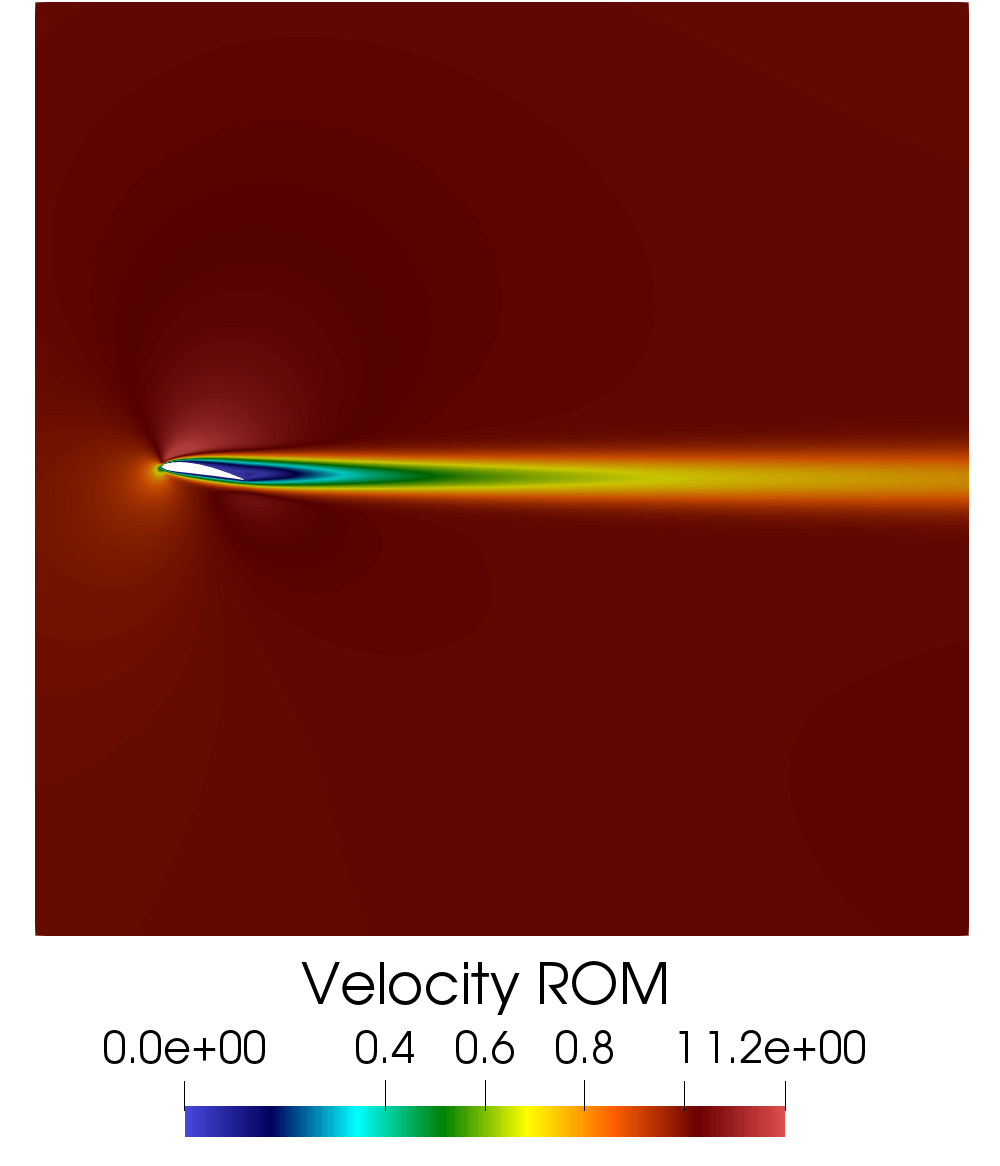}
\end{minipage} 
\begin{minipage}{0.24\textwidth}
\includegraphics[width=\textwidth]{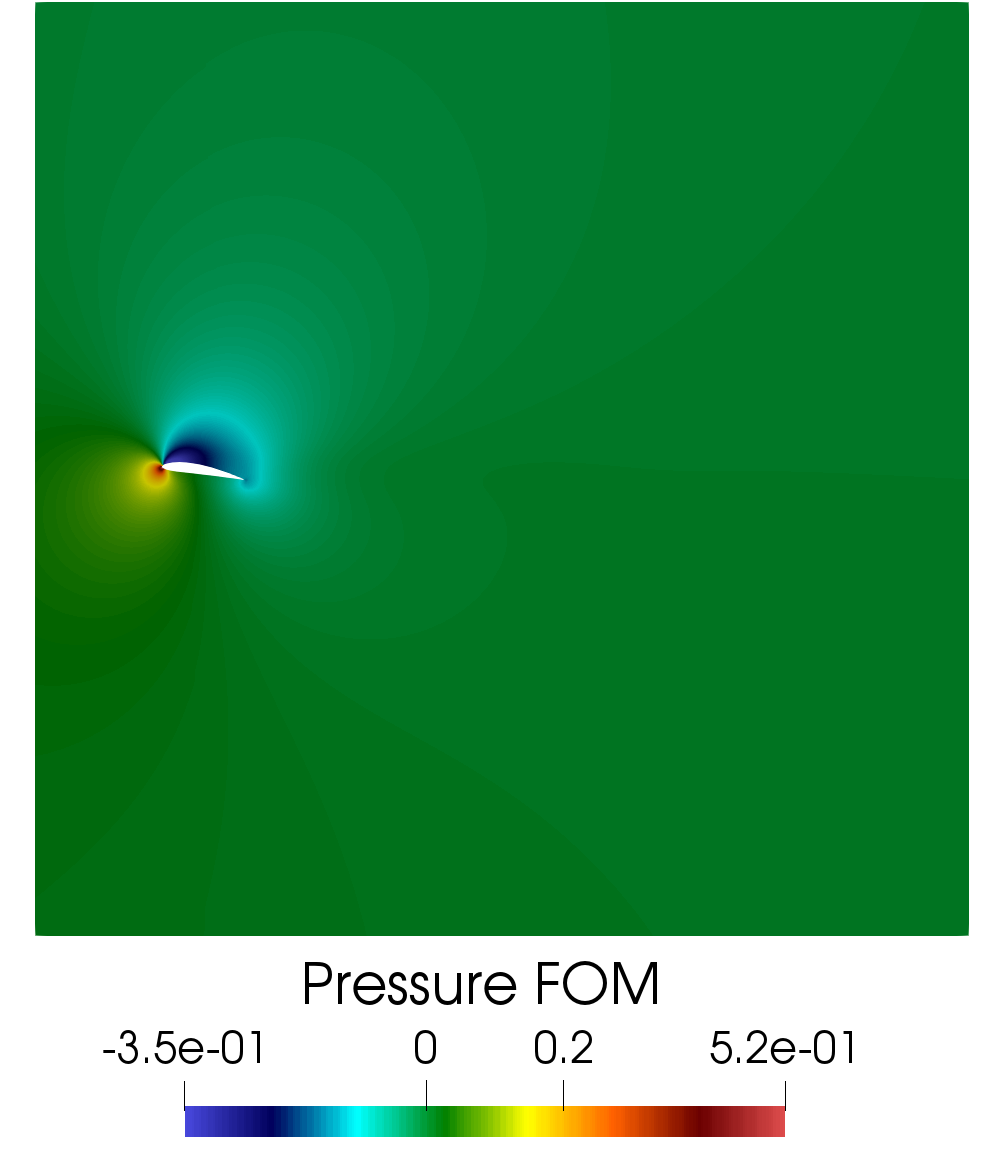}
\end{minipage} 
\begin{minipage}{0.24\textwidth}
\includegraphics[width=\textwidth]{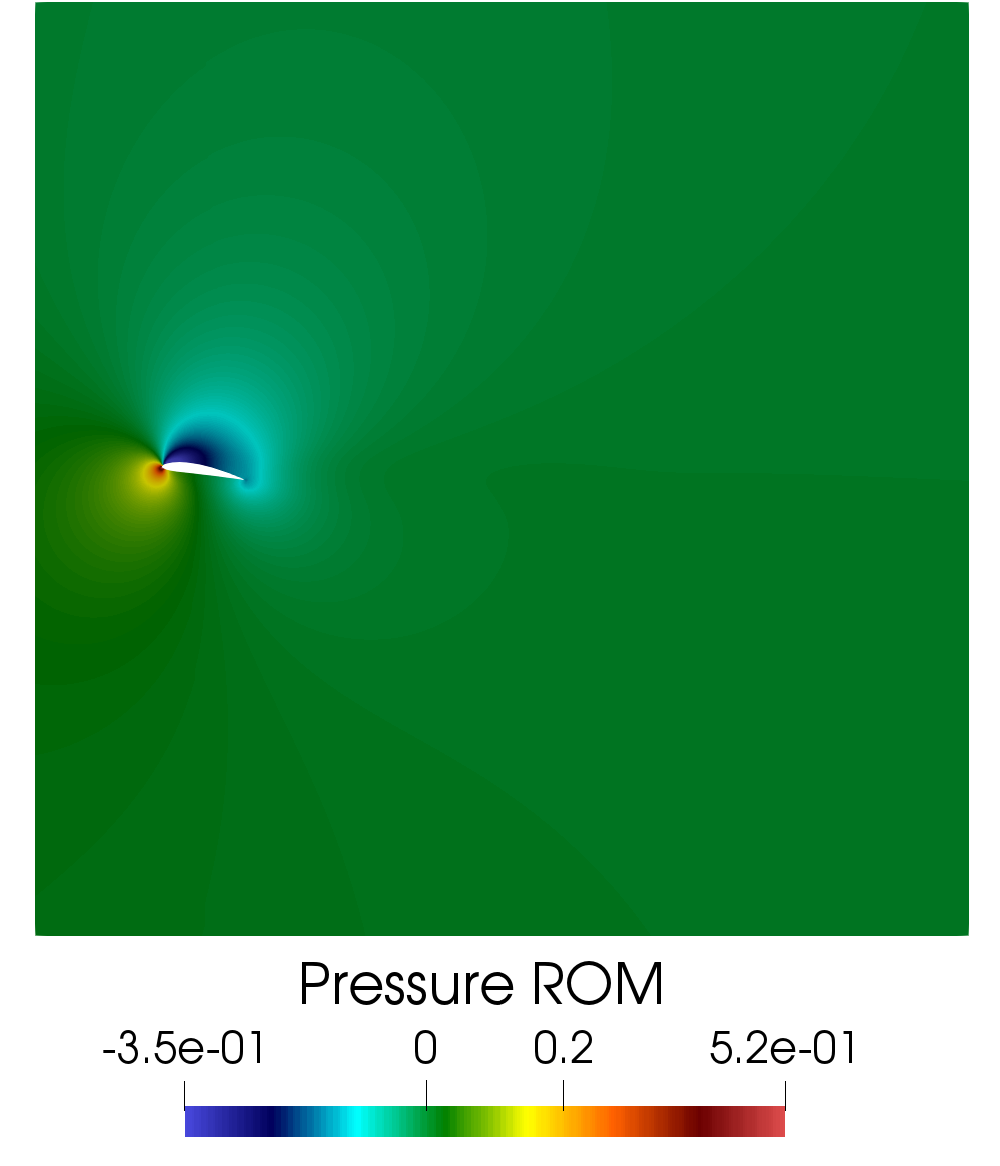}
\end{minipage}
\end{minipage}

\vspace{0.1cm}

\begin{minipage}{\textwidth}
\centering
\begin{minipage}{0.24\textwidth}
\includegraphics[width=\textwidth]{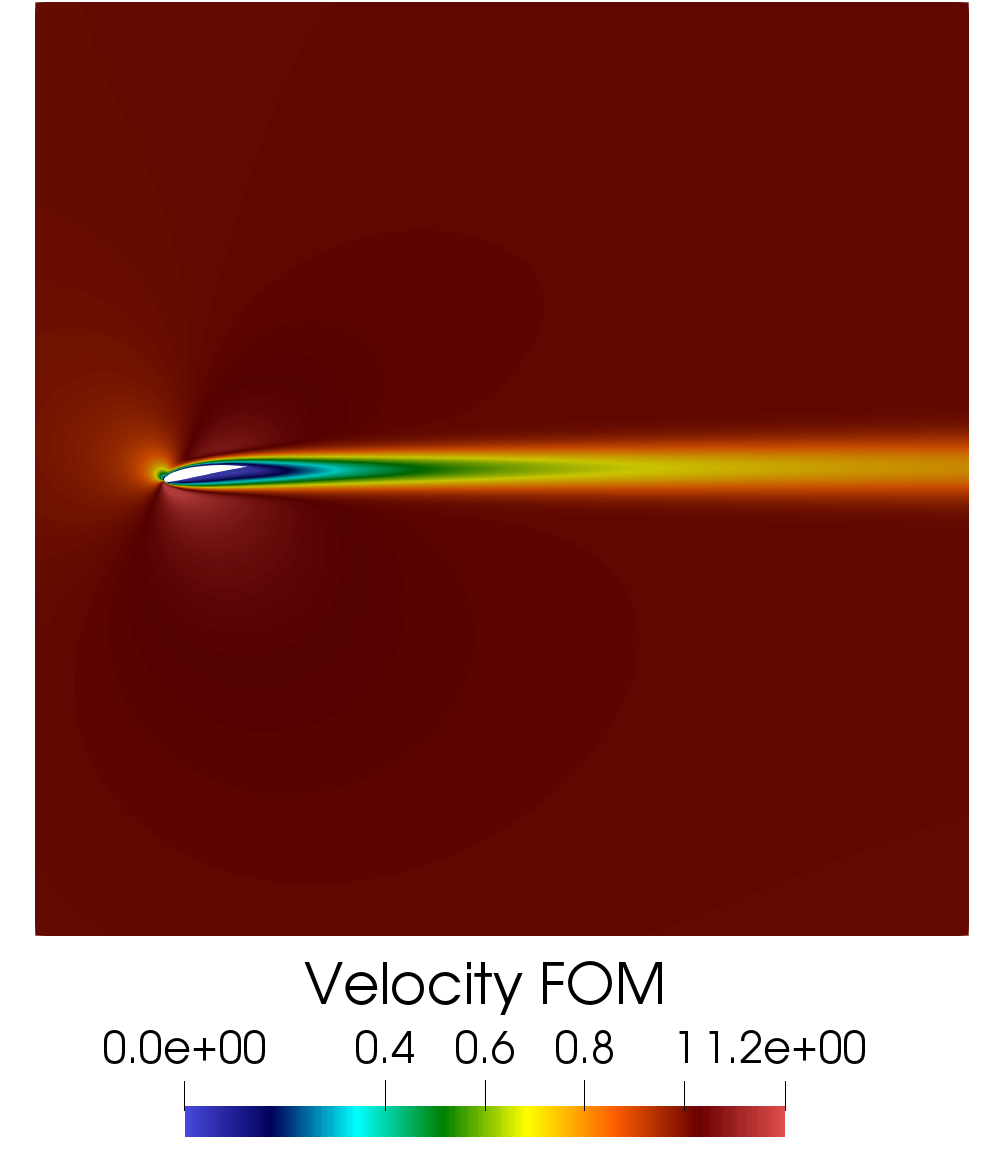}
\end{minipage} 
\begin{minipage}{0.24\textwidth}
\includegraphics[width=\textwidth]{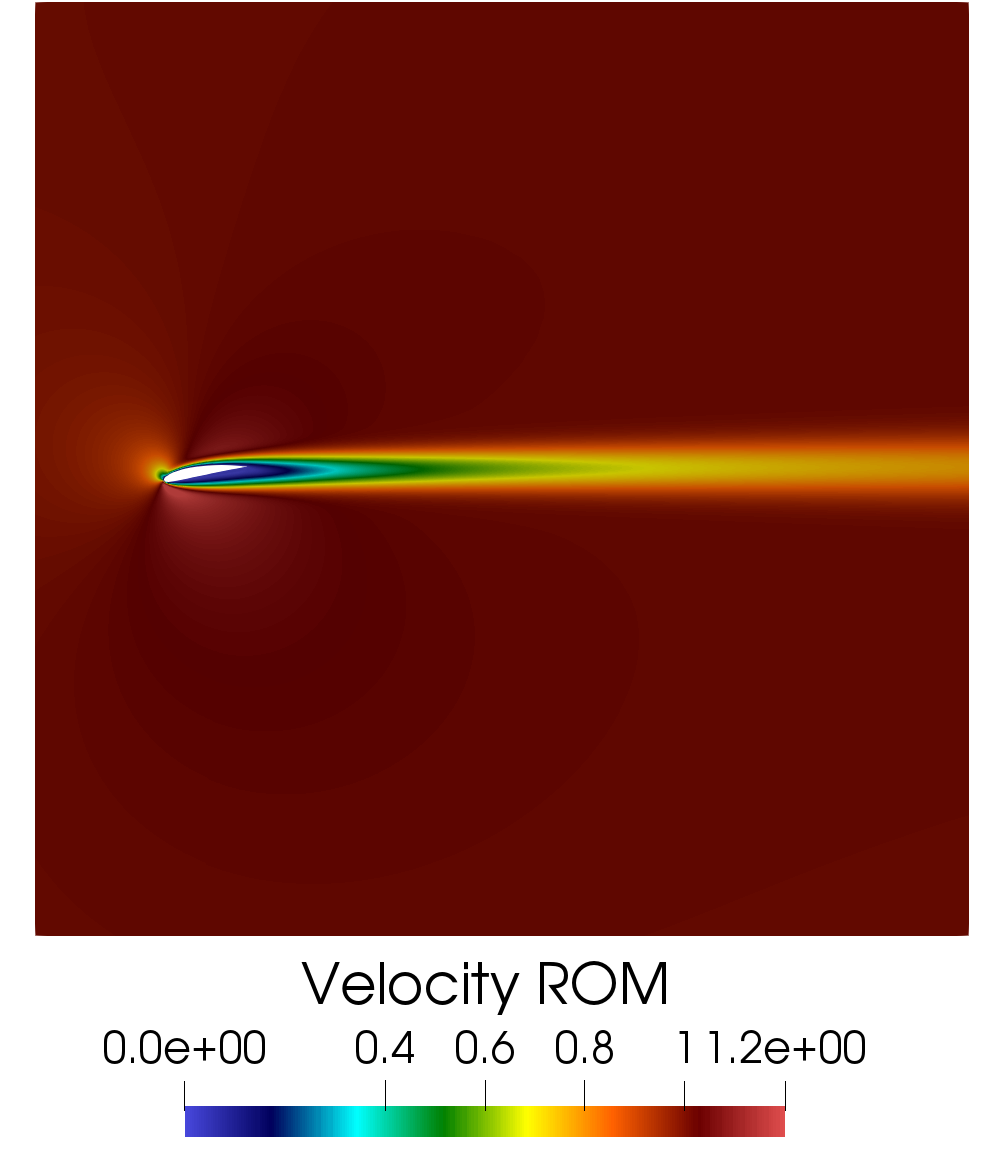}
\end{minipage} 
\begin{minipage}{0.24\textwidth}
\includegraphics[width=\textwidth]{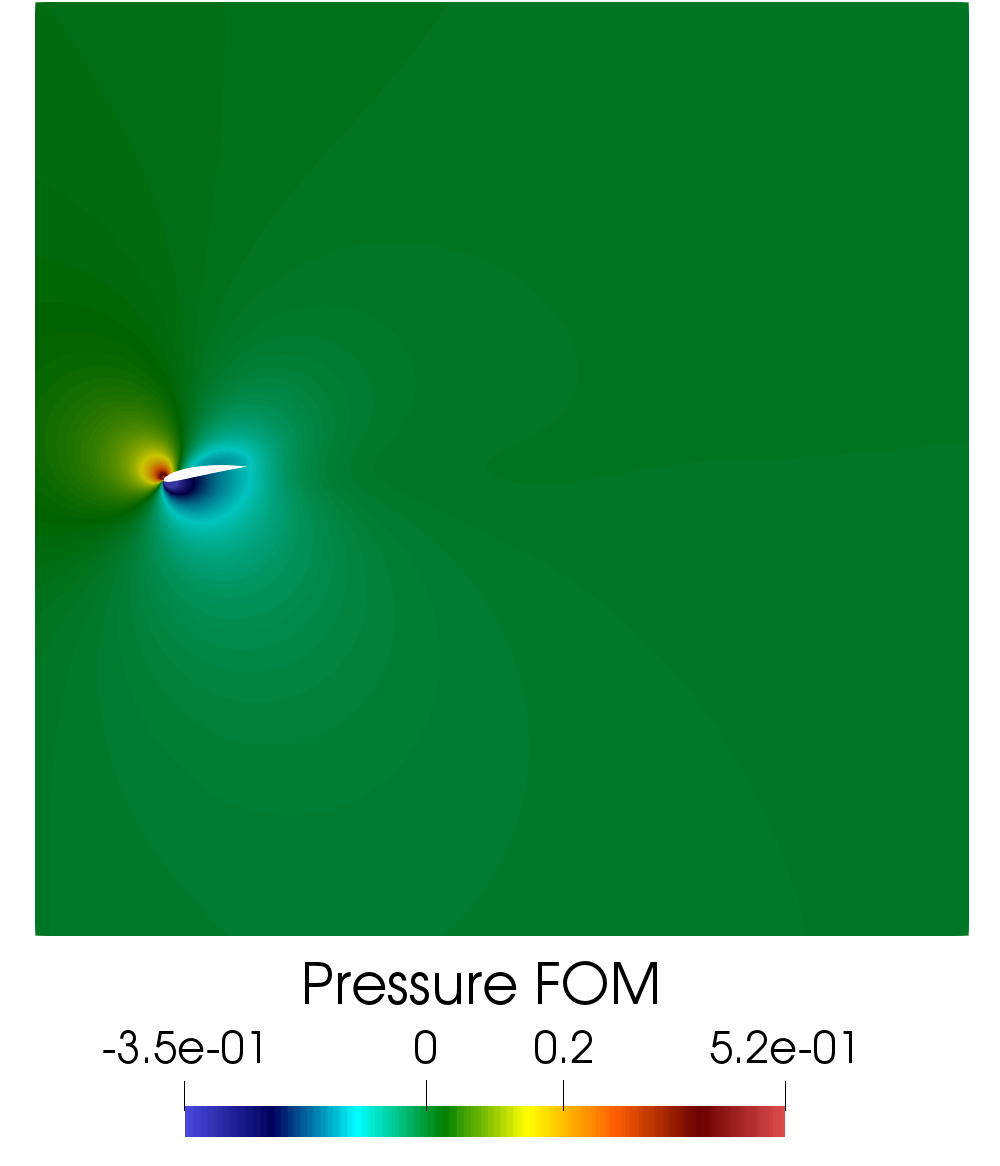}
\end{minipage}
\begin{minipage}{0.24\textwidth}
\includegraphics[width=\textwidth]{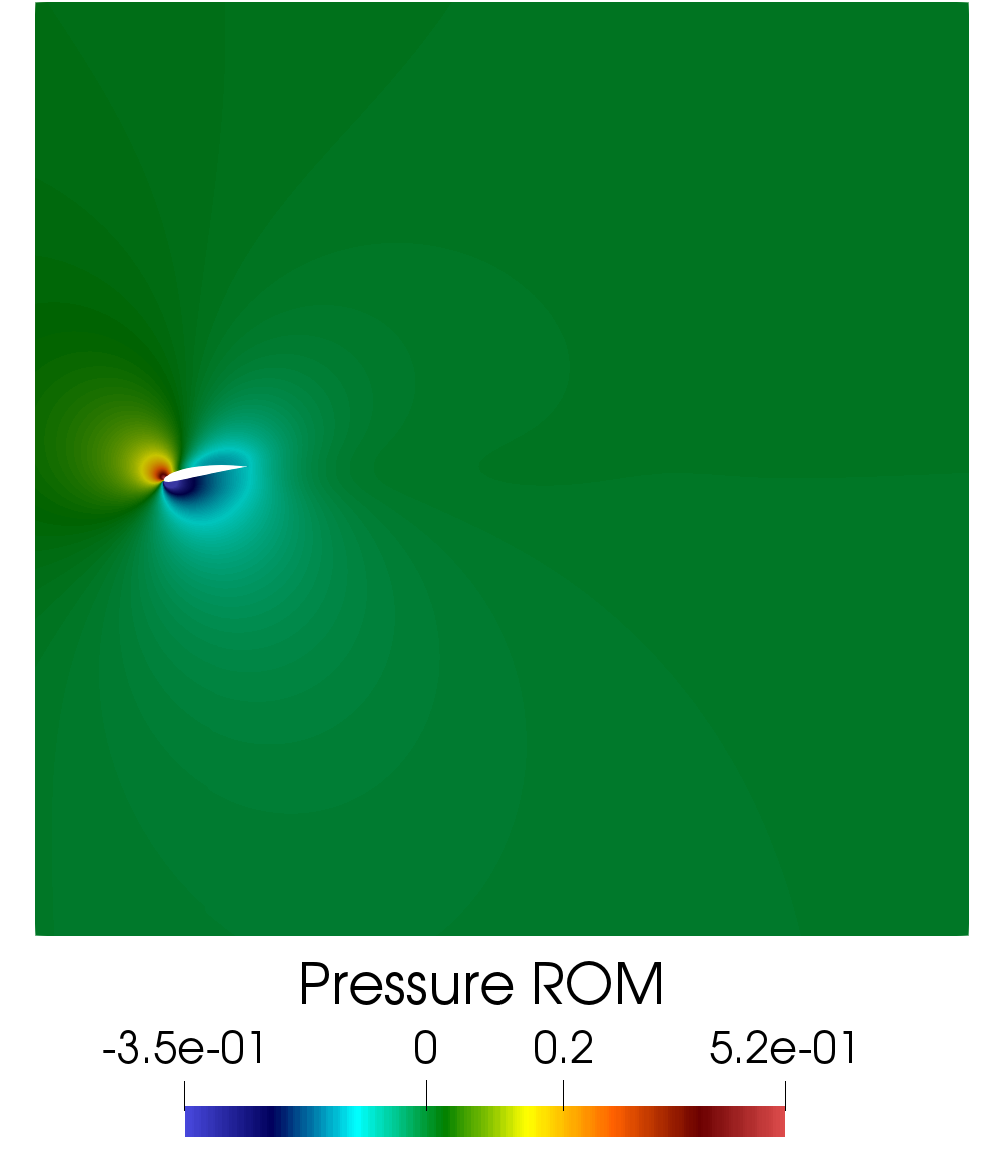}
\end{minipage}
\end{minipage}

\vspace{0.1cm}

\begin{minipage}{\textwidth}
\centering
\begin{minipage}{0.24\textwidth}
\includegraphics[width=\textwidth]{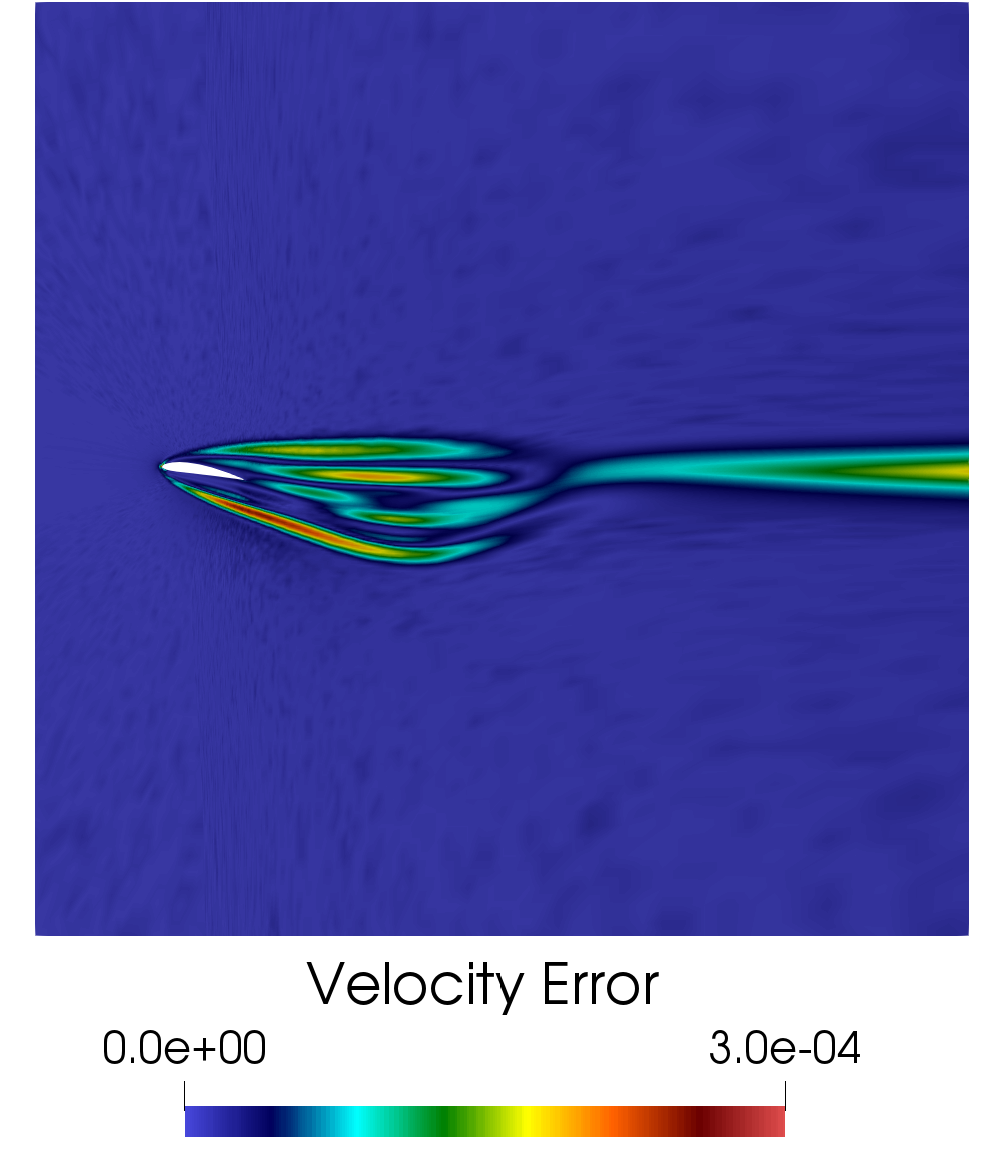}
\end{minipage} 
\begin{minipage}{0.24\textwidth}
\includegraphics[width=\textwidth]{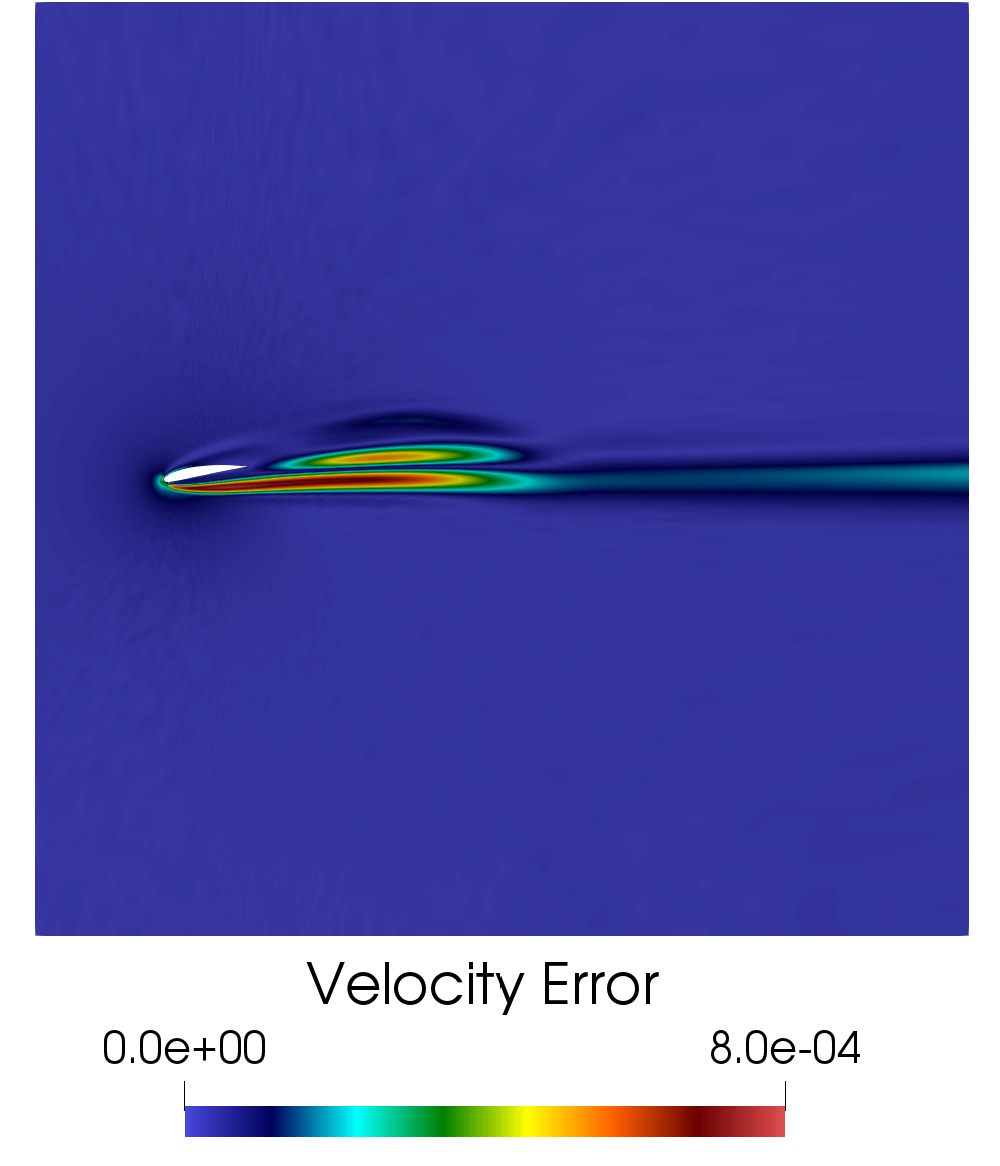}
\end{minipage}
\begin{minipage}{0.24\textwidth}
\includegraphics[width=\textwidth]{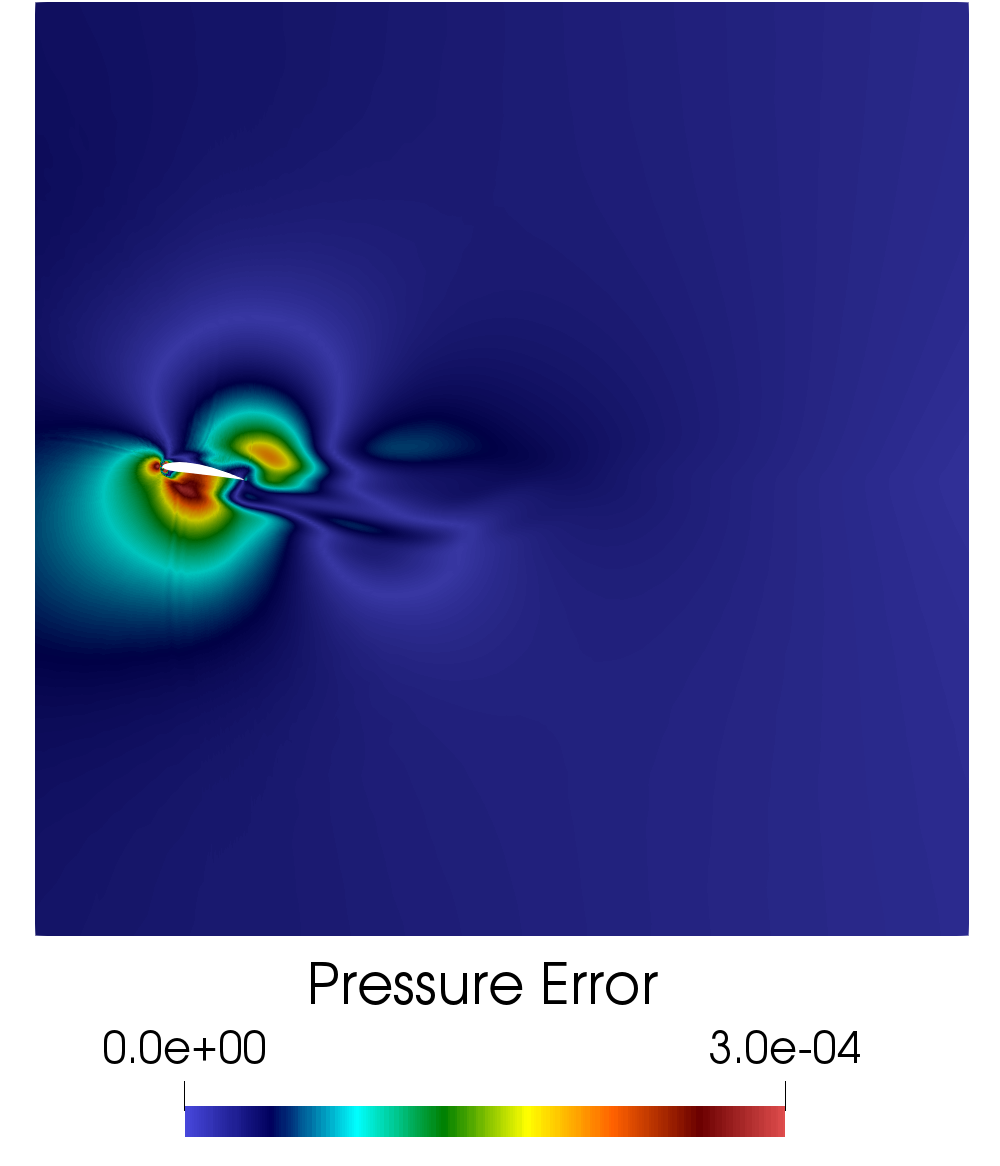}
\end{minipage} 
\begin{minipage}{0.24\textwidth}
\includegraphics[width=\textwidth]{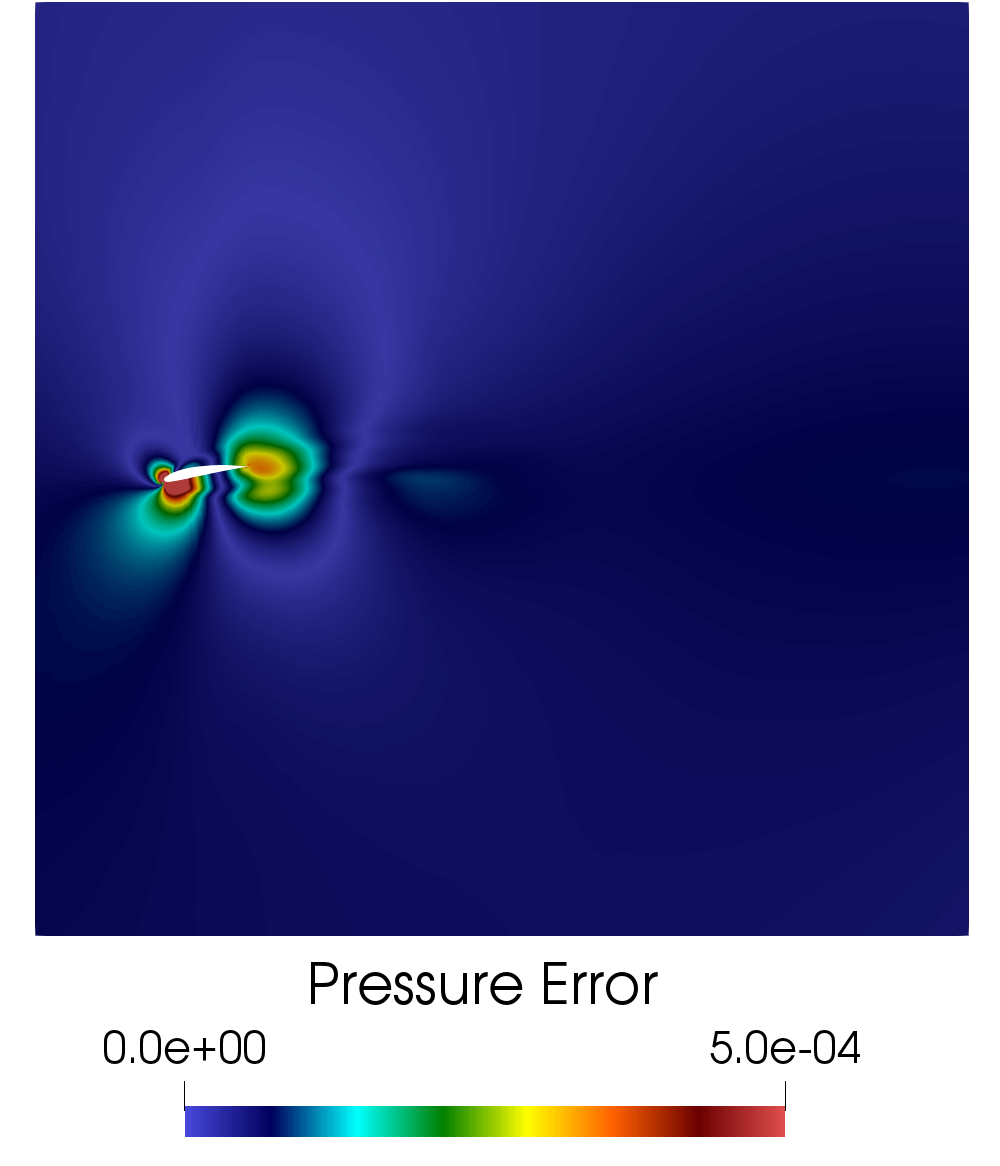}
\end{minipage}     
\end{minipage} 
\caption{Plots of the full order and reduced order velocity and pressure fields and of the error between them for two selected extreme values inside the testing set (${\kappa}^1_\alpha = \ang{9.15}$,  ${\kappa}^2_\alpha = \ang{-9.43}$). In the first and second row, from the left to the right, we report the FOM velocity field, the ROM velocity field, the FOM pressure field and the ROM pressure field for the first and second testing points, respectively. In the third row, from the left to the right, we report the velocity error for the first testing point, the velocity error for the second testing point, the pressure error for the first testing point and the pressure error for the second testing point, respectively.}  
\label{fig:NS_plots}
\end{figure}

In this section we present the numerical results for a geometrically parametrized incompressible Navier--Stokes problem. The case involves the parametrized angle of attack on a wing airfoil NACA $4412$. In this case, since the focus is mainly the methodological development we decided not to introduce additional complexities given by turbulence models and we focused indeed on relatively small Reynolds number without turbulence modeling. The numerical test has been performed using only the RBF mesh motion strategy as introduced in the previous numerical example. In this case the control points have been placed both on the moving boundaries and on the static patches. The locations on the airfoil are depicted in \autoref{fig:control_wing}. The kernel of the RBF interpolation is still given by Gaussian functions with a radius $d_{RBF}= 0.1$. In \autoref{fig:control_wing} we report also the geometry of the problem with the dimensions of the domain. The domain is composed by a quadrilateral whose size is equal to $16.5$ along the $x$ direction, $16$ along the $y$ direction, and by a semicircle attached to the inlet side of the quadrilateral. The wing chord is equal to $1$ and the foil is positioned on the center of the semicircle. The mesh counts $58000$ hexahedral cells. The geometrical parameter is given by the angle of attack that describes the rotation of the wing chord with respect to the inflow velocity. The deformed geometries are obtained rotating the wind around its barycenter. The training set $\mathcal{K}_{train} = \{\bm{\kappa}_{i_{train}}\}_{i=1}^{N_{train}} \in [\ang{-10} ,\ang{10}]$ with $N_{train}=100$ has been generated randomly inside the parameter space. The testing set $\mathcal{K}_{test} = \{\bm{\kappa}_{i_{train}}\}_{i=1}^{N_{test}} \in [\ang{-9.5} ,\ang{9.5}]$ that has been used to verify the accuracy of the ROM counts $N_{test}=50$ samples and has also been generated using a uniform random distribution. In the full order simulations the convective term has been approximated using a second order upwinding scheme while the diffusive terms has been discretized with a linear scheme with non-orthogonal correction. The SIMPLE algorithm runs with under relaxation for both velocity and pressure using the relaxation factors $\alpha_u = 0.7$ and $\alpha_p = 0.3$. As mentioned in the previous section the reduced order model has been constructed to be fully consistent with the SIMPLE procedure employed at the full order level. Therefore the same under-relaxation factors are used also at the reduced order level. The velocity at the inlet is set constant and equal to $1 \si{m}/\si{s}$. The kinematic viscosity is equal to $\nu=3 \times 10^{-3} \si{m}^2/\si{s}$. In \autoref{fig:NSerrors} we show the convergence properties of the ROM changing the number of reduced basis functions used to approximate both velocity and pressure while in \autoref{fig:NS_plots} we show a qualitative comparison between the full order and reduced order velocity and pressure fields for two selected values inside the testing set. The model produces accurate results for both velocity and pressure and more importantly it does not require any additional stabilization which is typical of reduced order models for incompressible flows \cite{stabile_stabilized}. We believe that  this fact is a consequence of the SIMPLE algorithm that has been used also at the reduced order level.

\begin{figure}
\begin{minipage}{0.49\textwidth}
\centering 
\includegraphics[width=\textwidth]{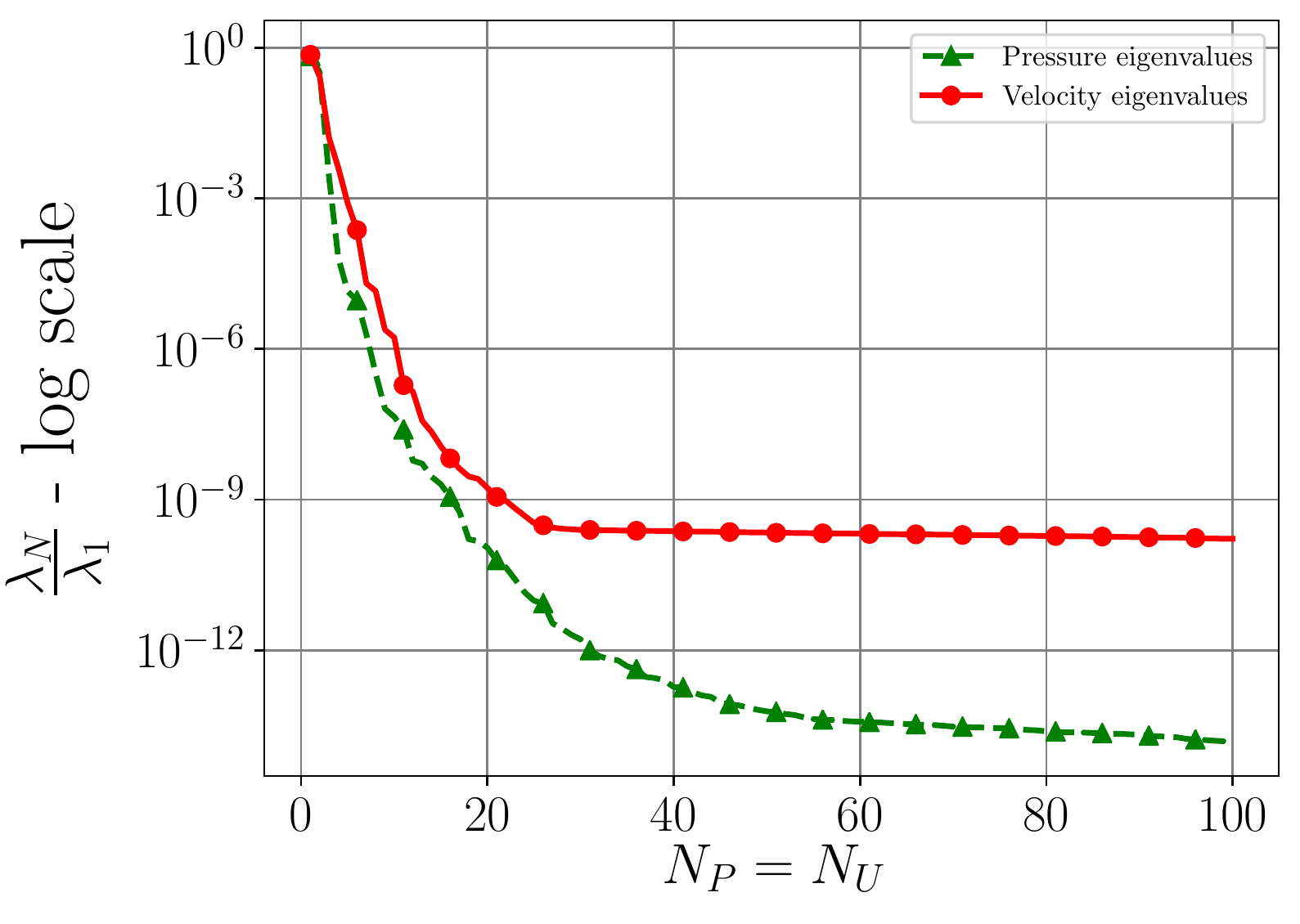}
\end{minipage}
\begin{minipage}{0.49\textwidth}
\centering 
\begin{tabular}{l | c | c}
  & Pressure & Velocity\\
  \hline
  N = 1 & $0.674404466672$ & $0.717187625760$\\
  \hline
  N = 2 & $0.997304641239$ & $0.979072177457$\\
  \hline
  N = 5 & $0.999988310322$ & $0.999728547289$\\
  \hline
  N = 7 & $0.999999521373$ & $0.999981203705$\\
  \hline
  N = 10 & $0.999999956844$ & $0.999999559123$\\
  \hline
  N = 20 & $0.999999999818$ & $0.999999980819$\\
  \hline
  N = 40 & $0.999999999997$ & $0.999999988151$\\
\end{tabular}
\end{minipage}
\caption{Plot reporting the eigenvalue decay relative to the POD procedure used to compute the modes for both pressure and velocity on the left, while on the right the cumulates of the eigenvalues are reported.}
\label{fig:NSeigs}
\end{figure}

{\subsection{A second numerical result on a geometrically parametrized incompressible flow problem}
In this section we show the numerical results obtained exploiting the segregated method explained in the previous sections on a shape deformation test case. In particular, by following what has been done in \cite{LeGresley2001}, we want to tackle the shape deformation of a NACA $4412$ airfoil by the superposition of some bump functions (\autoref{fig:BumF}) to the geometry. The same case and strategy analysed in the previous section are taken into consideration by only changing the geometry of the foil. All the $5$ bump functions are added to the top of the airfoil and subtracted to the bottom so that intersections between the boundaries are avoided, thus the problem is parametrized by the use of $10$ shape coefficients. The angle of attack of the reference foil is always set to zero while solving the problem. The training set $\mathcal{K}_{train} = \{\bm{\kappa}_{i_{train}}\}_{i=1}^{N_{train}} \in [0 , 0.02]^{10}$, where each row contains $10$ different values and $N_{train}=1000$, has been generated randomly inside the parameter space. The testing set $\mathcal{K}_{test} = \{\bm{\kappa}_{i_{train}}\}_{i=1}^{N_{test}} \in [0 , 0.018]^{10}$ that has been used to verify the accuracy of the ROM counts $N_{test}=20$ samples and has also been generated using a uniform random distribution.
In \autoref{fig:NSeigsBubble} the trends of the errors and of the eigenvalues for both velocity and pressure are depicted. \autoref{fig:NS_plots_bubble} shows the comparison between FOM and ROM solutions for both velocity and pressure for this last test case where one of the $20$ online parameters sets have been employed to modify the shape of the airfoil. As one can notice the accordance between the two solutions is high and the accuracy of the method is satisfying.

\begin{figure}
\centering
\includegraphics[width=0.5 \textwidth]{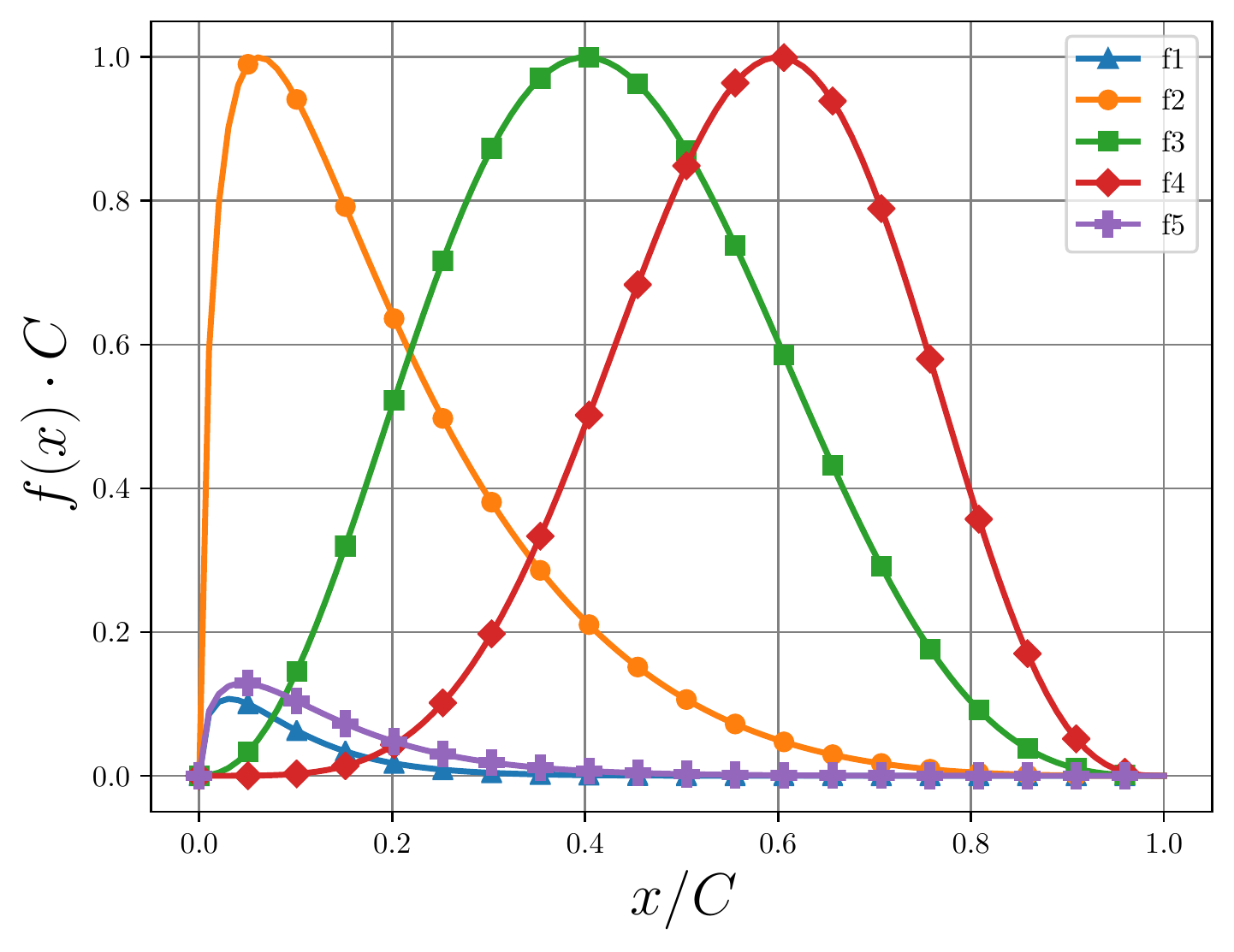}
\caption{Shape of the $5$ different bump functions, where $C$ is the chord of the airfoil.}
\label{fig:BumF}
\end{figure}

\begin{figure}
\begin{minipage}{0.49\textwidth}
\centering 
\includegraphics[width=\textwidth]{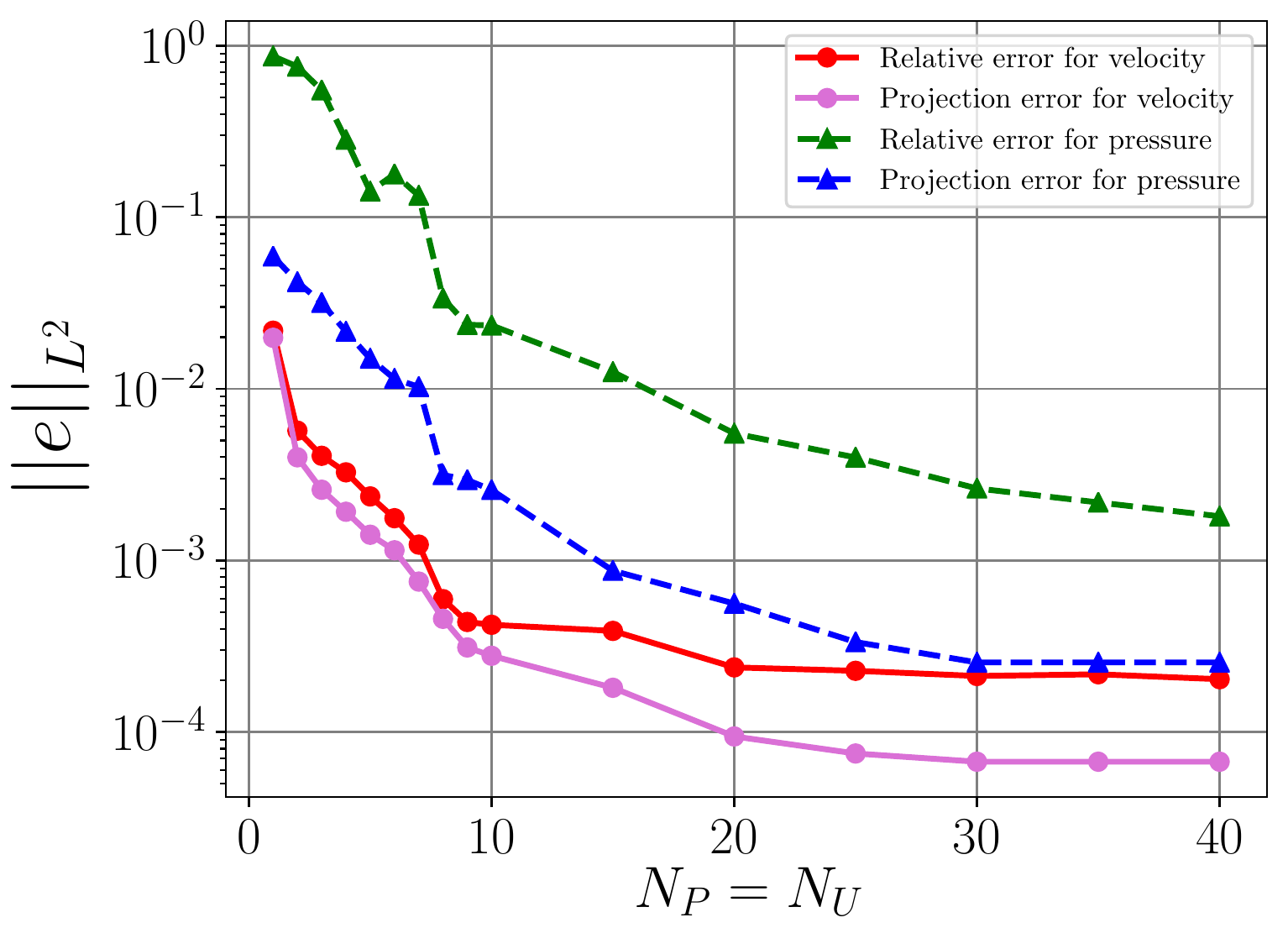}
\end{minipage}
\centering 
\begin{minipage}{0.49\textwidth}
\centering 
\includegraphics[width=\textwidth]{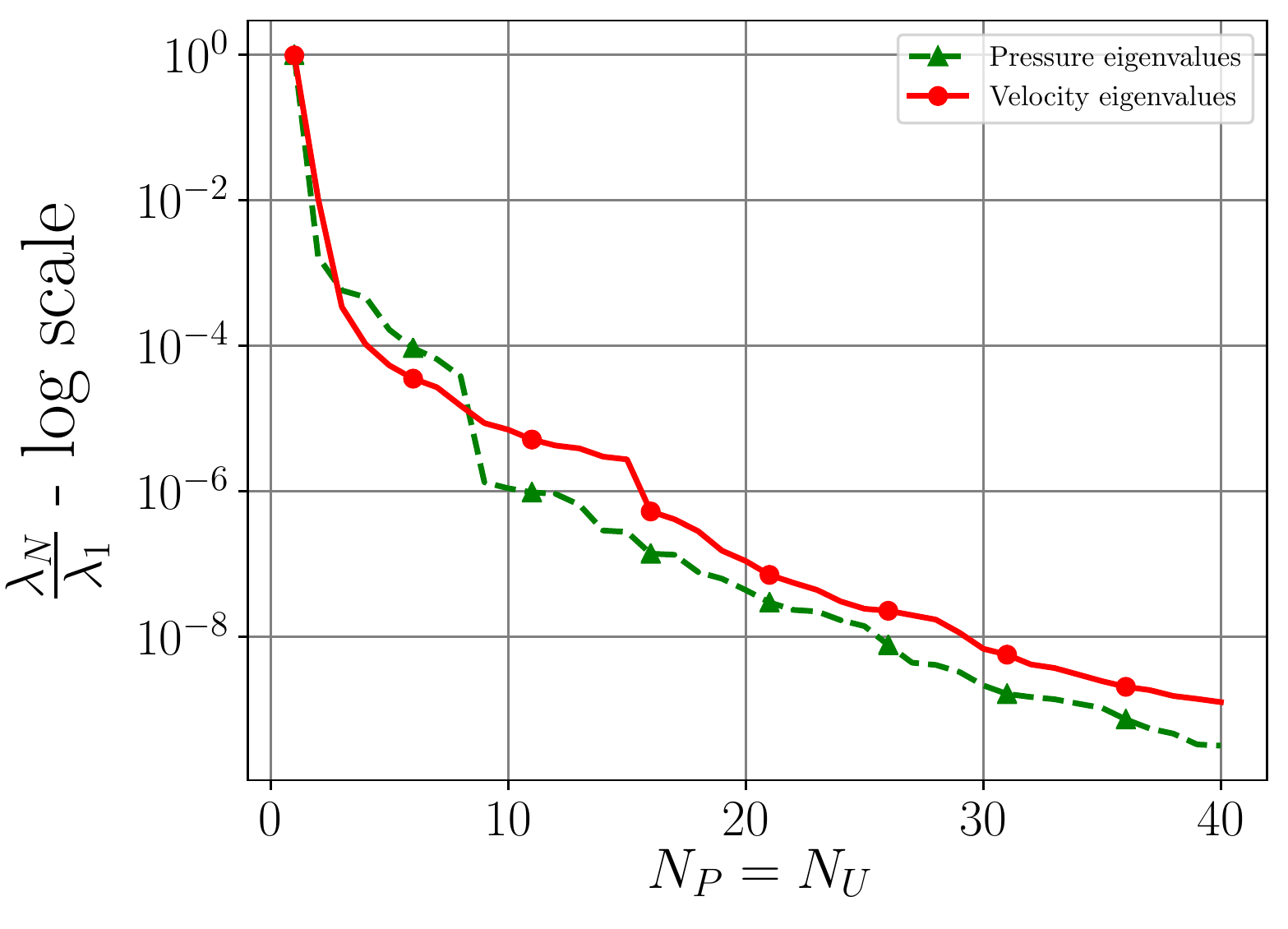}
\end{minipage}
\caption{\RA The average of the $L^2$ relative error norm of pressure and velocity fields has been reported on the left:  $\frac{P_{FOM} - P_{ROM}}{P_{FOM}}$ and $\frac{(U_{FOM} - U_{\infty})-(U_{ROM} - U_{\infty})}{U_{FOM} - U_{\infty}}$ where $U_{\infty}$ is the free stream velocity. The $L^2$ norm of the relative error for both pressure and velocity is plotted against the number of modes used for the reconstruction of the online solution in logarithmic $y$ scale. In the same image the reduced errors are compared with the projection errors obtained by the use of the same basis functions and by using the same norm. On the right the eigenvalue decay relative to the POD procedure used to compute the modes for both pressure and velocity is shown.}
\label{fig:NSeigsBubble}
\end{figure}

\begin{figure}

\begin{minipage}{\textwidth}
\centering
\begin{minipage}{0.24\textwidth}
\includegraphics[width=\textwidth]{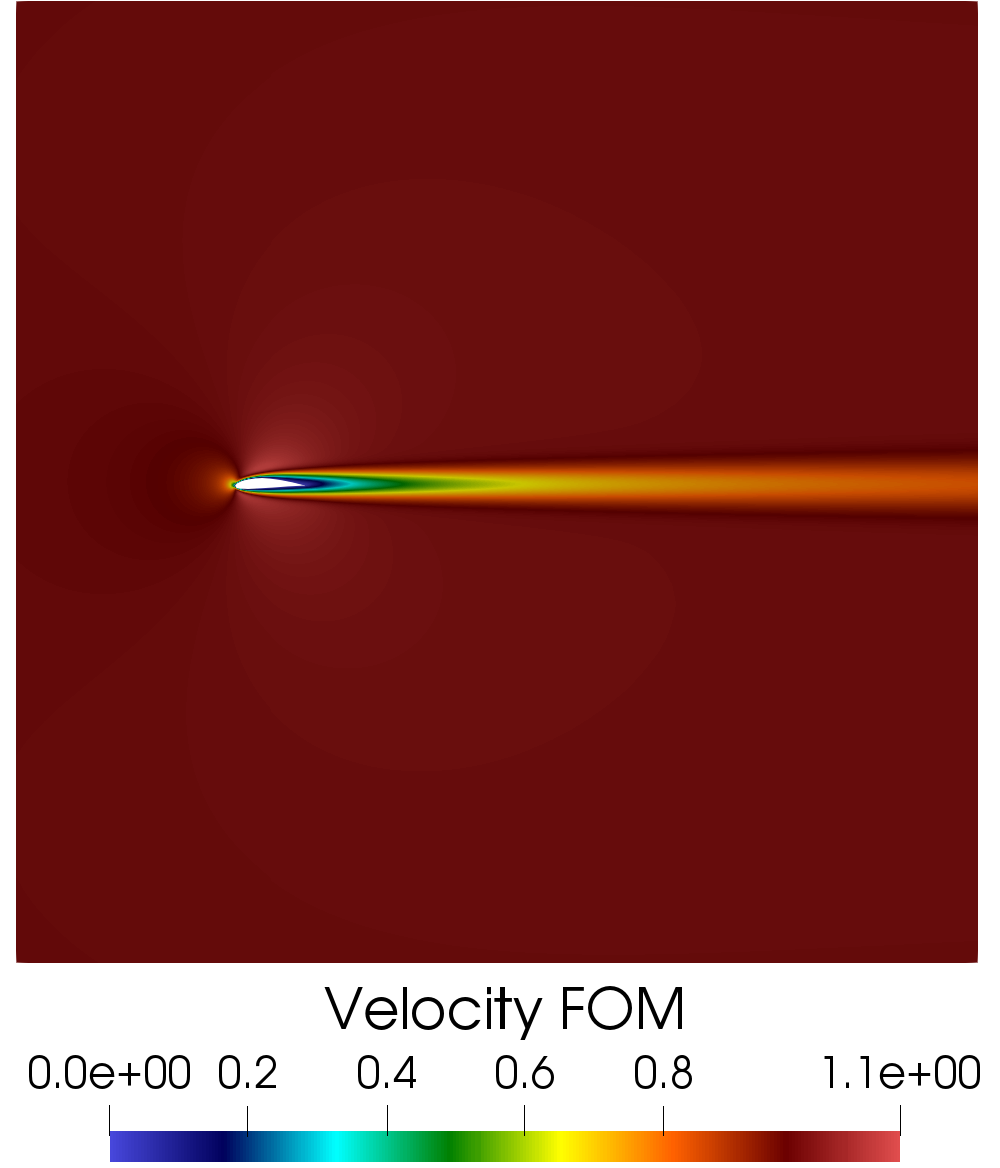}
\end{minipage}
\begin{minipage}{0.24\textwidth}
\includegraphics[width=\textwidth]{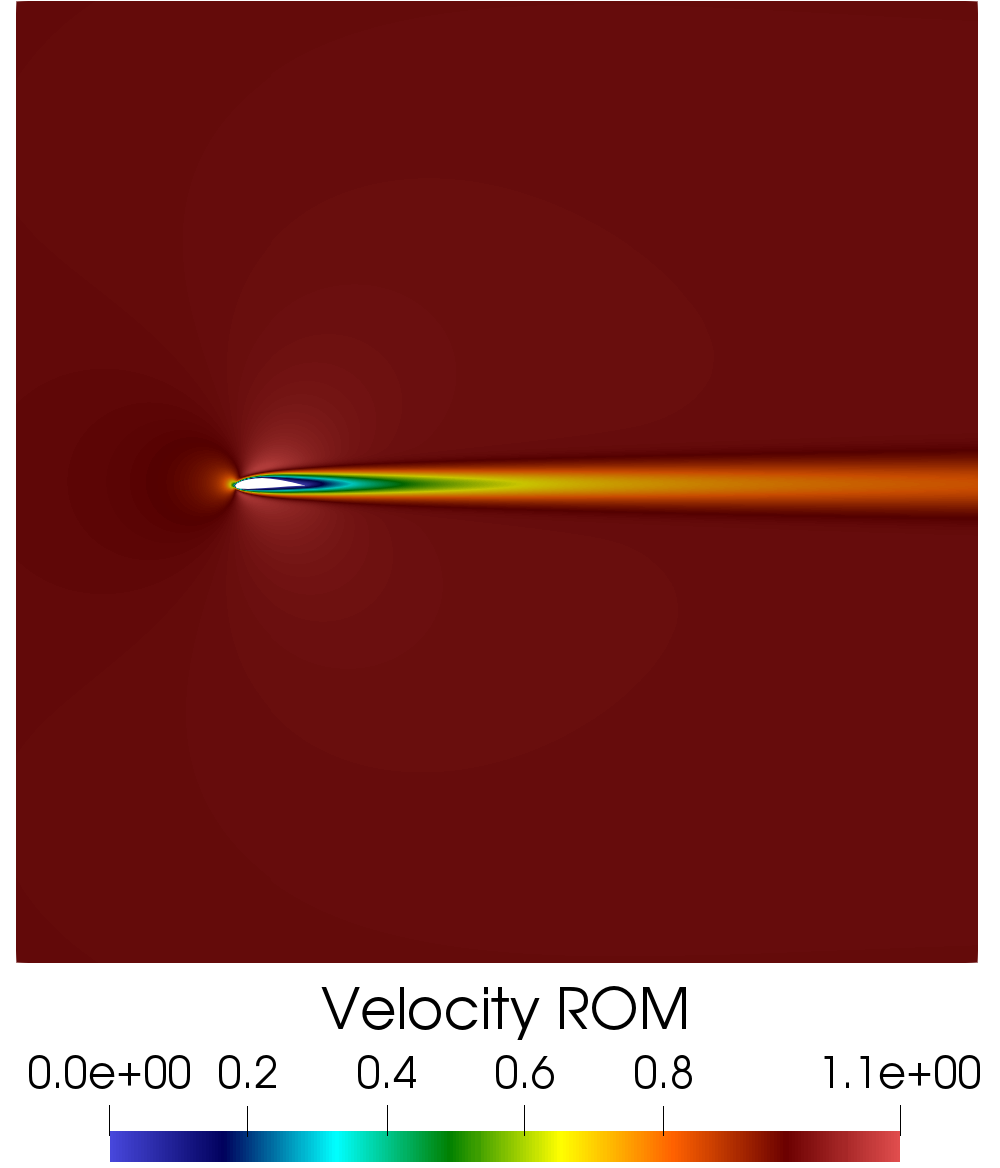}
\end{minipage} 
\begin{minipage}{0.24\textwidth}
\includegraphics[width=\textwidth]{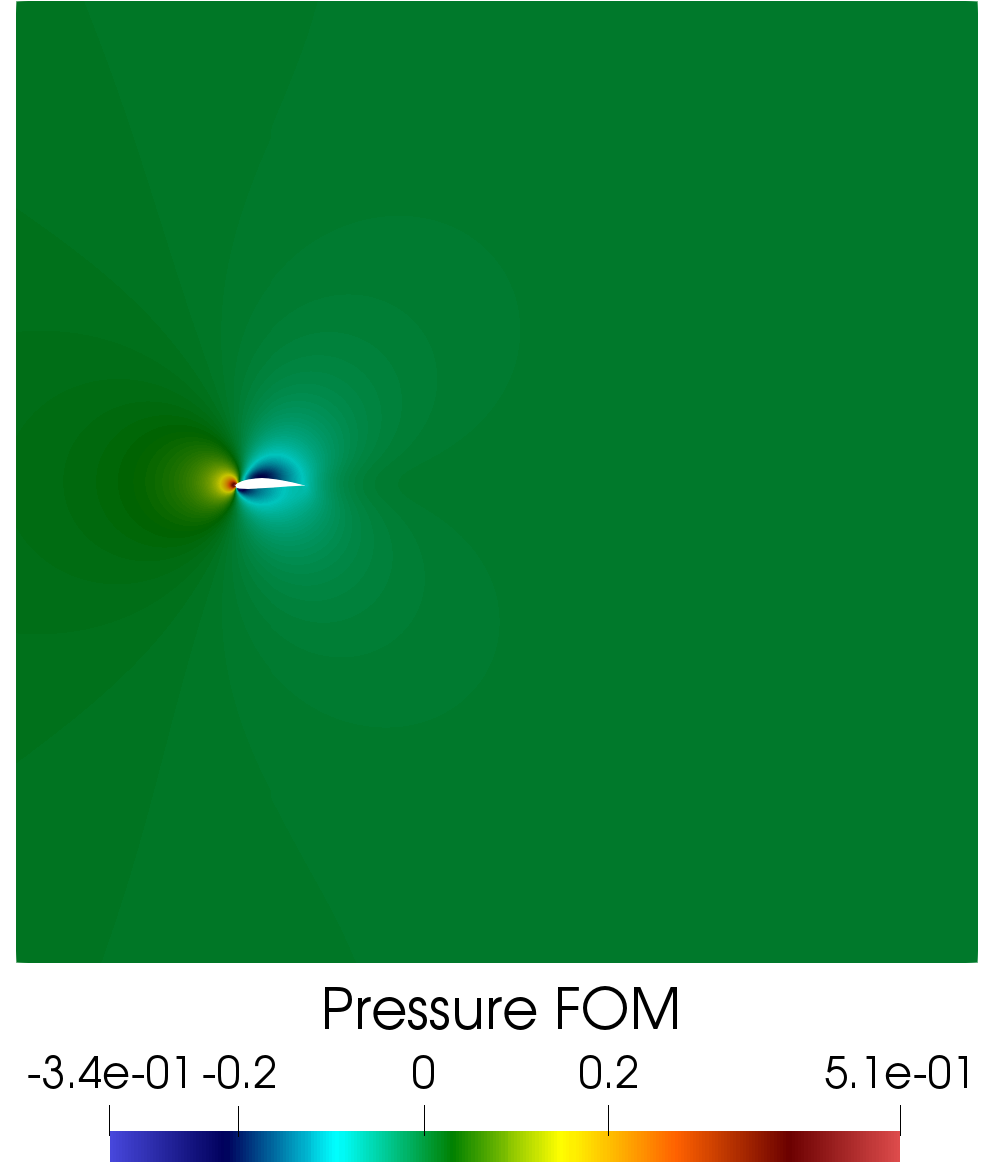}
\end{minipage} 
\begin{minipage}{0.24\textwidth}
\includegraphics[width=\textwidth]{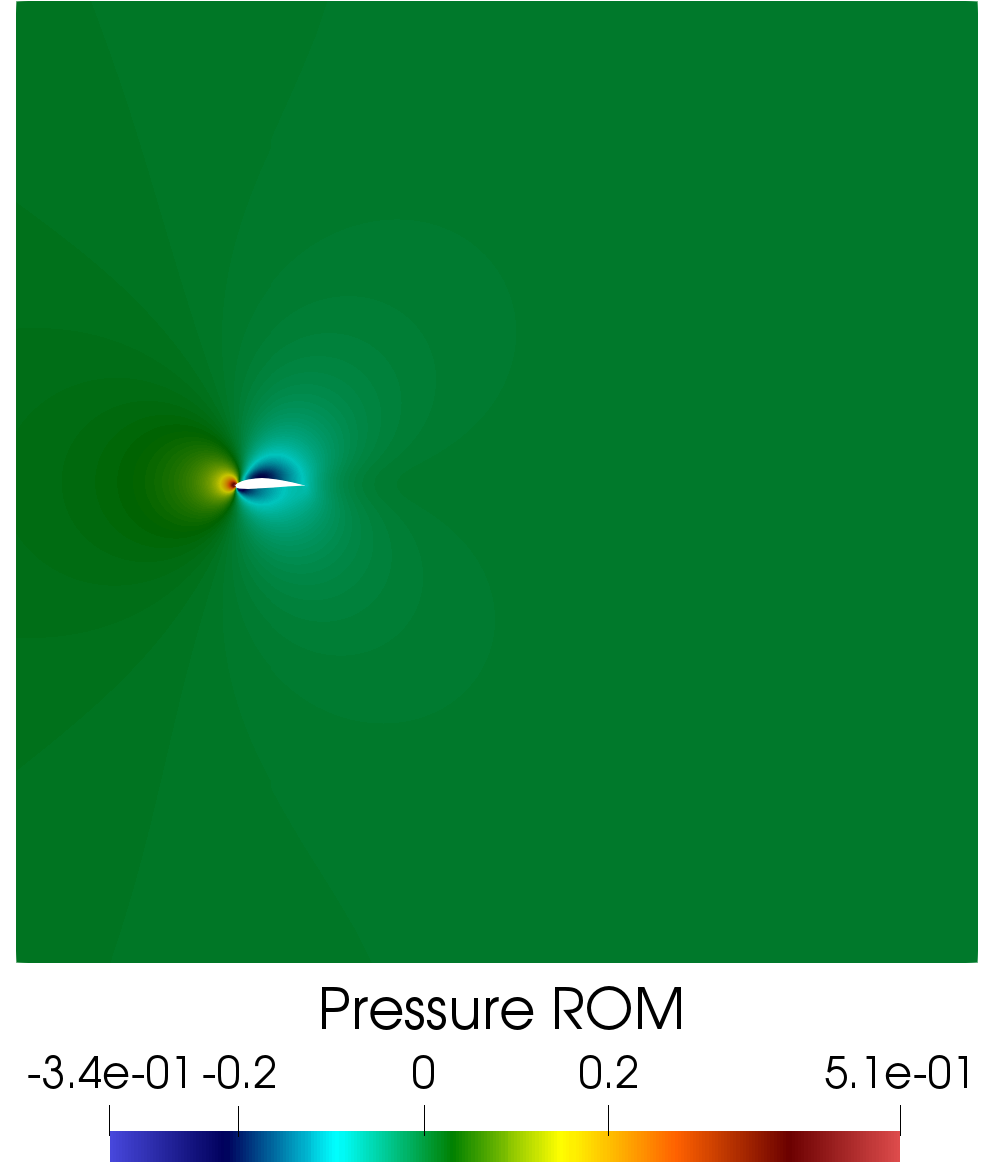}
\end{minipage}
\end{minipage}

\vspace{0.1cm}

\begin{minipage}{\textwidth}
\centering
\begin{minipage}{0.24\textwidth}
\includegraphics[width=\textwidth]{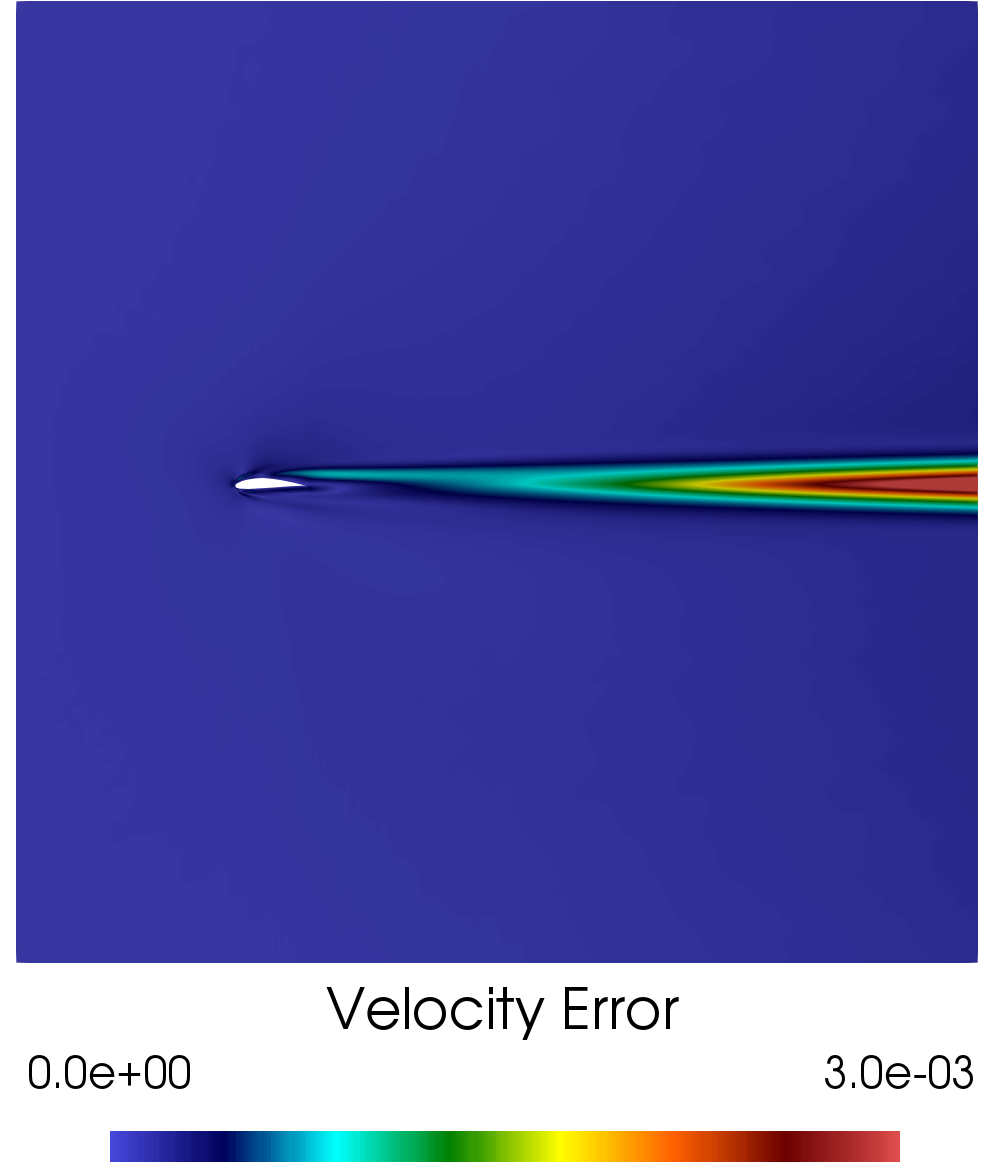}
\end{minipage} 
\begin{minipage}{0.24\textwidth}
\includegraphics[width=\textwidth]{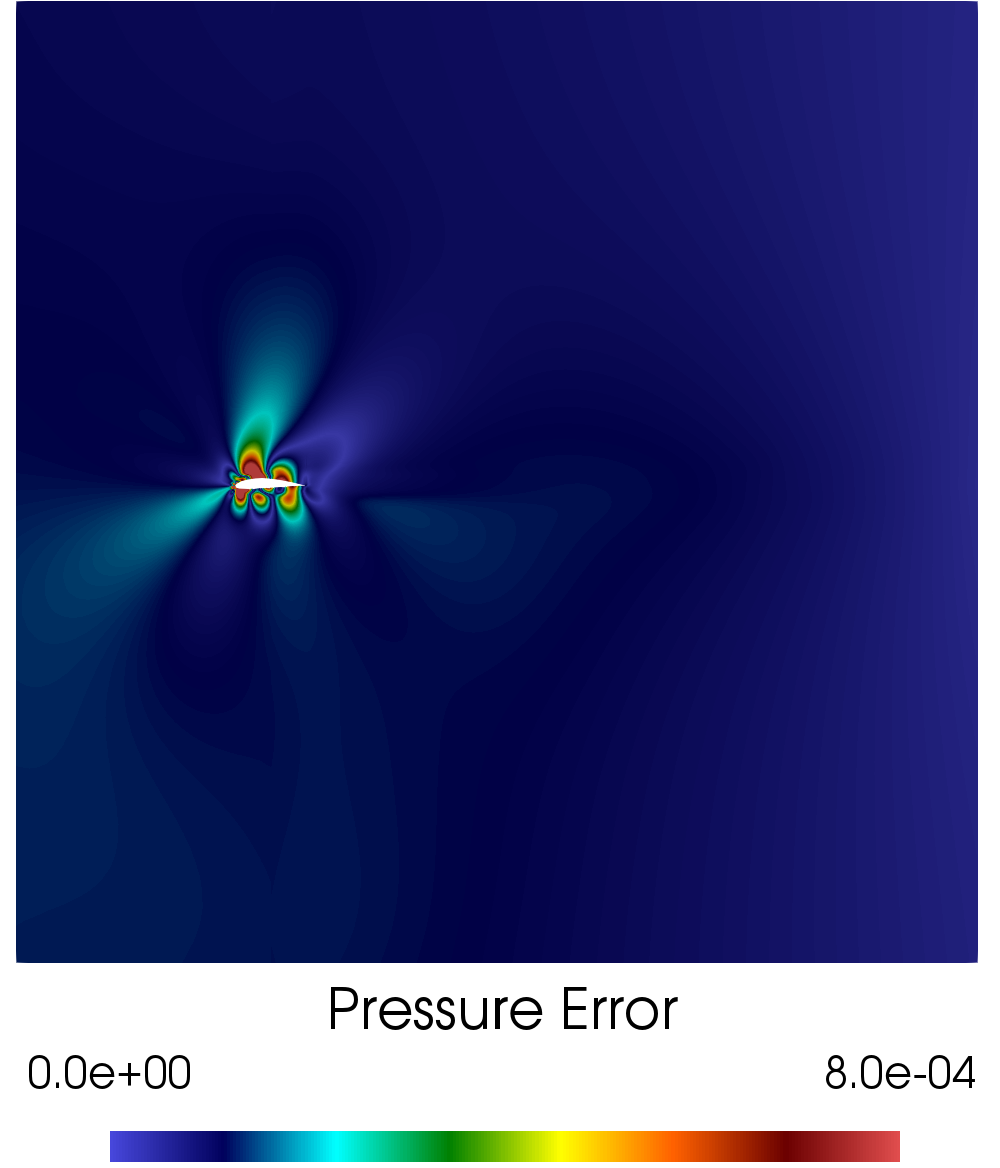}
\end{minipage}     
\end{minipage} 
\caption{Plots of the full order and reduced order velocity and pressure fields and of the error between them for a set of shape parameters inside the testing set ($[1.68124,$ $0.01488,$ $0.98581,$ $0.27690,$ $0.91593,$ $0.30225,$ $0.93038,$ $1.18651,$ $0.00667,$ $1.05945$ $]\times 10^{-2}$). In the first row, from the left to the right, we report the FOM velocity field, the ROM velocity field, the FOM pressure field and the ROM pressure field, respectively. In the second row, from the left to the right, we report the velocity error and the pressure error, respectively. Online solutions have been obtained by the use of $30$ basis functions.}  
\label{fig:NS_plots_bubble}
\end{figure}

\section{Conclusions and future perspectives}\label{sec:outlooks}
In this work we presented a reduced order modeling strategy for geometrical parametrization starting from a full order finite volume discretization. The methodology makes use of an ALE approach to analyze all the possible parametrized configurations and of a modified inner product for the computation of the correlation matrices used during the POD procedure. During the online stage, in order to ensure an efficient Offline-Online decoupling we make extensive use of the discrete empirical interpolation method at both matrix and vector level. 

Different mesh motion strategies have been tested and compared. These are based on a Laplacian smoothing  approach and on a radial basis function approach. Both methodologies have been analyzed with and without non-orthogonality correction and their perspectives in view of model reduction purposes have been deeply analyzed. 

Both methodologies have been tested on a heat transfer problem and on an incompressible steady Navier--Stokes setting providing accurate results and great computational gains. The radial basis function approach proved to be more accurate in presence of large parametric geometrical variations. 

As future perspectives it will be certainly interesting also to use the developed methodology to perform shape optimization in computational fluid dynamic problems at higher Reynolds number \cite{Lombardi_2012} and to couple the present approach to what developed in \cite{1807.11370v1,saddam2018} for turbulent flows. 

Of interest is also to study the applicability of the method to deal with complex geometrical variations arising from free form deformation algorithms.

Moreover, the interest is also into studying FSI transient problem, to deal with parametric interfaces. 

\section{Acknowledgments}
We acknowledge the support provided by the European Research Council Consolidator Grant project Advanced Reduced Order Methods with Applications in Computational Fluid Dynamics - GA 681447, H2020-ERC COG 2015 AROMA-CFD, as well as MIUR FARE-X-AROMA-CFD (Italian Ministry for Education University and Research) and INdAM-GNCS (Istituto Nazionale di Alta Matematica - Gruppo Nazionale di Calcolo Scientifico). The computations in this work have been performed with \emph{ITHACA-FV} \cite{RoSta17}, developed at SISSA mathLab, which is an implementation in OpenFOAM \cite{OF} of several reduced order modeling techniques; we acknowledge developers and contributors to both libraries.
%\section*{Appendix A. List of abbreviations and symbols}
%\begin{multicols}{2}
%\printnomenclature
%\end{multicols}
\bibliographystyle{abbrv}
\bibliography{bib_ok}
\end{document}